\documentclass[10pt,a4paper]{article}

\usepackage{graphicx}

\usepackage{amsmath}
\usepackage{amsthm}
\usepackage{amsfonts}
\usepackage{amssymb}
\usepackage{hyperref}
\usepackage{multirow}
\usepackage{units}
\usepackage{caption}
\usepackage{subcaption}
\usepackage{epstopdf}
\usepackage{cleveref}
\usepackage{tikz}
\usepackage{algorithm}
\floatname{algorithm}{Step}
\usepackage{enumitem}

\usepackage{algpseudocode}

\newtheorem{thm}{Theorem}

\newtheorem{conj}[thm]{Conjecture}
\newtheorem{rmk}{Remark}
\newtheorem{defi}{Definition}
\newtheorem{pf}{Proof}

\begin{document}

\title{High order approximation to non-smooth multivariate functions}

\author{Anat Amir\thanks{School of Mathematical Sciences, Tel Aviv University, Ramat Aviv, Tel Aviv 6997801, Israel} \and David Levin\footnotemark[1]}

\begin{abstract}
Common approximation tools return low-order approximations in the vicinities of singularities.  Most prior works solve this problem for univariate functions. In this work we introduce a method for approximating non-smooth multivariate functions of the form $f = g + r_+$ where $g,r \in C^{M+1}(\mathbb{R}^n)$ and the function $r_+$ is defined by 
\[ r_+(y) = \left\{ \begin{array}{ll}
r(y), & r(y) \geq 0  \\
0, & r(y) < 0 \end{array} \right. \ , \ \forall y \in \mathbb{R}^n \ . \]
Given scattered (or uniform) data points $X \subset \mathbb{R}^n$, we investigate approximation by quasi-interpolation. We design a correction term, such that the corrected approximation achieves full approximation order on the entire domain. We also show that the correction term is the solution to a Moving Least Squares (MLS) problem, and as such can both be easily computed and is smooth. Last, we prove that the suggested method includes a high-order approximation to the locations of the singularities. 
\end{abstract}

\maketitle

\section{Introduction}
\label{intro}
\par
Approximation of non-smooth functions is a complicated problem. Common approximation tools, such as splines or approximations based on Fourier transform, return smooth approximations, thus relying on the smoothness of the original function for the approximation to be correct. 
However, the need to approximate non-smooth functions exists in many applications. For a high-order approximation of non-smooth functions, we need to allow our approximation to be non-smooth. Otherwise, in the vicinities of the singularities, we will get a low-order approximation. In this work we will suggest a method that will allow us to properly approximate non-smooth functions of a given model.
\par 
We will concentrate on functions $f : \mathbb{R}^n \rightarrow \mathbb{R}$ which may be modelled as $f = g + r_+$ where $g,r \in C^{M+1}(\mathbb{R}^n)$ and the function $r_+$ is defined by 
\[ r_+(y) = \left\{ \begin{array}{ll}
r(y), & r(y) \geq 0 \\
0, & r(y) < 0 \end{array} \right. \ , \ \forall y \in \mathbb{R}^n \ . \]
Such functions are obviously continuous, but are non-smooth across the hypersurface 
\[ \Gamma_r := \left\{ z \in \mathbb{R}^n : r(z) = 0 \right\} \ . \]
As an example for such functions, consider shock waves, which are solutions of non-linear hyperbolic PDEs \cite{Acoustics}. 
Another example would be a signed distance function  \cite{LevelSet}, where the distance is measured from a disconnected set. Our goal is to achieve high-order approximations of such functions. To achieve that we will concentrate on a specific family of approximation tools.
\par 
Consider a quasi-interpolation operator $Q$ \cite{wendland2004scattered}. Such an operator receives the values of a function $\phi : \mathbb{R}^n \rightarrow \mathbb{R}$ on a set of data points $X \subset \mathbb{R}^n$. 
The quasi-interpolation operator $Q$ returns an approximation defined by
\[ Q \phi(y) := \sum_{x \in X} q_{x}(y) \phi(x) \ , \ \forall y \in \mathbb{R}^n \ , \] 
where $\left\{ q_{x} \right\} $ are the quasi-interpolation basis functions, each is smooth and has compact support.

Let $h$ be the fill distance of $X$, 
\[ h := \min \left\{ L : B_L(y) \cap X \neq \emptyset \ , \ \forall y \in \mathbb{R}^n \right\} \ , \]
where $B_r(y)$ is the ball of radius $r$ centred at $y$. 
Denote 
\[ \Upsilon_M := \min {\left\{ L > 0 : \forall y \in \mathbb{R}^n \ , \ B_{L h}(y) \cap X \mbox{ is uni-solvent for } \Pi_M(\mathbb{R}^n) \right\}} \ . \]
Here, 
\[ \Pi_M(\mathbb{R}^n) := \left\{ p : \mathbb{R}^n \rightarrow \mathbb{R} : deg(p) \leq M \right\} , \]
and $deg(p)$ is the total degree of the polynomial $p$.
Thus, $h$ is the minimal radius which is guaranteed to contain a data point, and $\Upsilon_M \cdot h$ is the minimal radius that guarantees enough data points to uniquely determine each polynomial in $\Pi_M(\mathbb{R}^n)$.
We will also assume that there exists $N > 0$ such that for all $y \in \mathbb{R}^n$ we have 
\[\dfrac{\#\left( X \cap B_h(y) \right)}{h^n} \leq N \ . \] That is, the data set $X$ has no accumulation points. 
Denote 
\begin{equation} \label{eq:R}
R := min \left\{ \rho > 0 : supp (q_{x}) \subseteq B_{\rho h}(x) \ , \ \forall x \in X \right\} 
\end{equation}
We assume that the operator $Q$ has a bounded Lebesgue constant 
\begin{equation} \label{eq:L_1}
L_1 := \sup \left\{ \sum_{x \in X} \vert q_{x}(y) \vert : y \in \mathbb{R}^n \right\} 
\end{equation}
and reproduces polynomials in $\Pi_M(\mathbb{R}^n)$.
Then, the error in the quasi-interpolation, 
\[E \phi := \phi - Q \phi \]
satisfies for all $\phi \in C^{M+1}(\mathbb{R}^{n})$ and $y \in \mathbb{R}^n$ 
\[ \vert E \phi(y) \vert \leq C_1 \cdot \Vert \phi \Vert_{C^{M+1}} \cdot h^{M+1} \]
where 
\[ C_1 = (1 + L_1) \cdot {R}^{M+1} \] 
and 
\[ \Vert \phi \Vert_{C^{M+1}} := \sum\limits_{\vert \beta \vert = M+1} \frac{\Vert D^{\beta} \phi \Vert_{\infty}}{\beta!} \]
with $\beta$ a multi-index and $\Vert \cdot \Vert_{\infty}$ the maximum norm.
That is, the operator $Q$ has full approximation order for smooth functions \cite{wendland2004scattered}. On the other hand, since the approximation $Q \phi$ is always smooth, the operator gives low-order approximations in the vicinities of singularities.
\par
One example of a quasi-interpolation operator is the MLS approximation \cite{bosMLS,levin1998approximation}. Given a function $\phi : \mathbb{R}^n \rightarrow \mathbb{R}$ and a point $y \in \mathbb{R}^n$ the MLS approximation is defined as $Q \phi(y) := p_{y}(y)$ where 
\begin{equation} \label{eq:MLS}
p_{y} := \operatorname*{arg\,min}\limits_{p \in \Pi_M(\mathbb{R}^n)} \sum\limits_{x \in X} \eta \left( \frac{\Vert y-x \Vert }{h}\right) \cdot \left( p(x) - \phi (x) \right)^2 
\end{equation}
Here $\eta$ is a smooth weight function with compact support. The MLS approximation essentially returns the value at the point $y$ of the $M$-th degree polynomial which gives the best approximation to $\phi$ at the data points 
\[ X \cap supp \left( \omega \left( \frac{ \Vert \cdot - y \Vert}{h} \right) \right) \ . \]
This approximation is especially important in this work, since we base much of our results on our ability to adapt the MLS approximation to our non-smooth scenario.
\par 
Our initial goal for this work was to generalize the work done by Lipman and Levin \cite{Lipman}. In that work, the authors address the problem of approximation of univariate functions of the form 
\[ f(x) = g(x) + \sum_{j=1}^M \frac{\Delta_j}{j!}(x-s)^j_+ \ , \]
where $g \in C^{M+1}(\mathbb{R})$ and 
$(x-s)^j_+ = \left\{ \begin{array}{ll} 
(x-s)^j , & x \geq s \\
0 , & x < s \end{array} \right.$. 
Indeed, such functions are continuous but not smooth. In \cite{Lipman}, the univarite case is solved by modelling the error terms 
of the approximation by a quasi-interpolate $Q$. That is, one searches for variables $s^\ast, \bar{\Delta}^\ast$, such that the errors in the quasi-interpolation approximation of the term 
\[ \tilde{r}_{\bar{\Delta}^\ast, s^\ast}(x) := \sum_{j=1}^M \frac{\Delta^\ast_j}{j!}(x-s^\ast)^j_+ \]
give the best Least-squares approximation to the errors of the function $f$ at the data points. 
It is shown that by adding the error of the approximation of the new term $\tilde{r}_{\bar{\Delta}^\ast, s^\ast}$ to the approximation $Qf$, full approximation order for the function $f$ is achieved.
\par
Another approach to this problem was proposed by Harten \cite{Harten}. The author introduces the essentially non-oscillatory (ENO) and the subcell resolution (SR) schemes. The ENO scheme bases the approximation at each point on only some of the data points in its vicinity. Thus, disregarding points from the other side of the singularity which contaminate the approximation. The SR scheme locates the singularities by intersecting polynomials from supposedly different sides of the singularities. For an examination of these methods for univariate functions with a jump discontinuity in the derivative see \cite{Arandiga}.
\par
Archibald et al (\cite{Gelb2009}, \cite{Gelb2008}) suggest using polynomial annihilation to locate the singularity $\Gamma_r$. Of-course, once the singularity is known we can approximate each connected component of $\mathbb{R}^n \setminus \Gamma_r$ independent of the values in the other connected components. 
\par 
Other approaches were suggested by Markakis and Barack \cite{markakis2014high}, where the authors revise the Lagrange interpolation formula to approximate univariate discontinuous functions, and by Plaskota et al (\cite{Plaskota2009power}, \cite{Plaskota2013adaptive}), where the authors suggest using adaptive methods for this approximation. Batenkov et al (\cite{batenkov2012complete}, \cite{Batenkov2011algebraic}, \cite{Batenkov2012algebraic}) address a similar problem, the reconstruction of a piecewise smooth function from its integral measurements. One disadvantage of the methods mentioned above, is that they do not easily adapt to multivariate singular functions.
\par 
Thus, the main advantage of the method we suggest in this paper is its ability to deal with the multivariate case. Indeed, our method enables us to approximate multivariate functions which have non-continuous derivatives across smooth hyper-surfaces,  $\Gamma_r$. Note that while the dimension $n$ of the domain of the function $f$ affects the required number of data points in $X$ and the dimension of $\Pi_M(\mathbb{R}^n)$, the correction procedure is not otherwise affected by the higher dimension.

\section{Main results} \label{sec:results}
As described in the introduction, our goal is to fix the approximation of the function $f = g + r_+ $, where $g , r \in C^{M+1}(\mathbb{R}^n)$.
\begin{rmk}
The decomposition $f = g + r_+$ is not unique. Indeed, 
\[ f = g + r_+ = g + r + (-r)_+ \ . \]
However, we only correct the approximation error, for which we do have uniqueness.
\end{rmk}
We will achieve this, following the main idea in \cite{Lipman}, by investigating the error terms 
\[E f(y) = \underbrace{E g (y)}_{O(h^{M+1})} + E r_+ (y) = E r_+(y) + O(h^{M+1}) \ . \] 

\begin{defi}[$\lambda$-neighbourhood of $\Gamma_r$]
For $\lambda > 0$ define 
\[ \mathcal{G}(\lambda) := \bigcup_{z \in \Gamma_r} B_{\lambda}(z) \ . \]
\end{defi}

Note that if the point $y$ is far enough from the singularity, the restriction of $f$ to the $Rh$ neighbourhood of $y$ is a smooth function, hence in this case there is no need to fix the approximation. Thus, we need only correct the approximation for points in the set  $\mathcal{G} (R h)$. 
In the following we suggest an algorithm for bounding the set $\mathcal{G}(Rh)$. We begin by estimating whether the function $r$ returns a positive or negative value at each data point $x \in X$. For this we will need the following definition: 
\begin{defi}[Partition of the data points with respect to sign] \label{defi:P}
For a set $X \subset \mathbb{R}^n$ with fill-distance $h$ and a function $r$ we will say that the set $P \subset X$ partitions the data points in $X$ with respect to the sign of the function $r$ if 
\begin{enumerate}
\item $\forall x \in P$ either $r(x) > 0$ or $r(x) = O(h^{M+1})$.
\item $\forall x \in X \setminus P$, either $r(x) < 0$ or $r(x) = O(h^{M+1})$.
\end{enumerate}
\end{defi}
In \cref{sec:signs} we will introduce an algorithm that partitions the data points in $X$ with respect to the sign of the function $r$.
For now, let us assume that a set $\mathcal{P}$ which partitions the data points in $X$ with respect to the sign of the function $r$ is known. 
\begin{rmk}
Apparently, once we have a set $\mathcal{P}$ which partitions the data points in $X$ with respect to the sign of the function $r$, we can approximate a point where $r$ has a positive value using only the data points in $\mathcal{P}$, and a point where $r$ has negative value using only the data points in $X \setminus \mathcal{P}$. However, we predict whether $r$ has a positive or negative value only on the data points $X$, and not on the entire domain. Specifically, for a point close to the singularity location, $\Gamma_r$, we can not tell whether we should approximate based on $\mathcal{P}$ or on $X \setminus \mathcal{P}$. Hence we can not rely only on the set $\mathcal{P}$ to fix the approximation.
\end{rmk}

\begin{defi}[The set $\mathcal{G}$] \label{defi:G}
Denote 
\[ \mathcal{B}(y) := B_{(R + 2 \Upsilon_M + 2)h}(y) \ , \] 
and define 
\[ \mathcal{G} := \left\{ y \in \mathbb{R}^n : \mathcal{B}(y) \cap \mathcal{P} \mbox{ and } \mathcal{B}(y) \cap (X \setminus \mathcal{P}) \mbox{ are uni-solvent for } \Pi_M(\mathbb{R}^n) \right\} \ . \]
\end{defi}

Of-course, 
\[\mathcal{G} (R h) \neq \mathcal{G} \ , \]
however, we can prove that 
\[ \mathcal{G} (Rh) \subset \mathcal{G} ((R+1) h) \subset \mathcal{G} \ . \]
This gives us a bound on the region in which we need to fix the approximations.
\begin{thm}[The domain of the correction] \label{thm:G_bound} If $\nabla r(z) \neq 0$ for all $z \in \Gamma_r$, then 
\[ \mathcal{G} ((R+1) h) \subseteq \mathcal{G}\ . \]
\end{thm}

\subsection{The corrected approximation}
We may now describe the corrected approximation of the function $f$.
Pick a point $y \in \mathbb{R}^n$. 
As explained above, there is no need to fix the approximation outside the set $\mathcal{G}$. Hence, if $y \in \mathbb{R}^n \setminus \mathcal{G}$ define the corrected approximation as 
\[ \widehat{Q} f(y) = Q f(y) \ . \]
Otherwise assume that $y \in \mathcal{G}$. To fix the approximation of $f$ we need a polynomial $p^y$ which approximates the function $r$ locally. In section \ref{sec:p^y} we will introduce two methods that will allow us to construct the approximation $p^y$.
\begin{defi} [Smooth $M$-th order approximation] \label{def:p^y}
We will say that a mapping $\Psi : \mathcal{G} \rightarrow \Pi_M(\mathbb{R}^n)$ is a smooth $M$-th order approximation to the function $r : \mathbb{R}^n \rightarrow \mathbb{R}$ if $\Psi$ satisfies the following conditions :
\begin{itemize}
\item There exists a constant $C_2 > 0$, independent of $h$, such that for all $y \in \mathcal{G}$, 
\[ \vert (\Psi(y))(u) - r(u) \vert \leq C_2 \cdot \left( \Vert g \Vert_{C^{M+1}} + \Vert r \Vert_{C^{M+1}}\right) \cdot h^{M+1} \ , \ \forall u \in B_{R h}(y) \ . \]
\item Let $\left\{ p_\alpha \right\}_{\alpha \in I}$ be a polynomial basis of $\Pi_M(\mathbb{R}^n)$ and write 
\[ \Psi(y) = \sum_{\alpha \in I} \lambda_\alpha (y) \cdot p_\alpha \ . \]
Then, the mappings $y \mapsto \lambda_{\alpha}(y)$ are infinitely smooth.
\end{itemize}
\end{defi}
If the mappings $y \mapsto p^y$ are a smooth $M$-th order approximation to $r$, then we may define the corrected approximation as 
\[ \widehat{Q} f(y) = Q f(y) + E \left( (p^y)_+ \right) (y) \ . \]
Hence, the corrected approximation $\widehat{Q}$ is defined as follows : 
\begin{defi} [Corrected approximation]
\[ \widehat{Q} f(y) := Q f(y) + \left\{ \begin{array}{ll}
E \left((p^y)_+ \right) (y), & y \in \mathcal{G} \\
0, & \mbox{otherwise} \end{array} \right. \ . \]
\end{defi}
Thus we get 
\begin{thm}[Corrected approximation errors]\label{thm:good_correction}
Let $y \mapsto p^y$ be a smooth $M$-th order approximation to $r$.
Then there exists a constant $C_3 > 0$ such that 
\[ \vert \widehat{E} f( y) \vert = \vert f(y) - \widehat{Q} f(y) \vert \leq C_3 \cdot \left( \Vert g \Vert_{C^{M+1}} + \Vert r \Vert_{C^{M+1}} \right) \cdot h^{M+1} \ . \]
\end{thm}

\begin{defi}[The function $\widehat{r}$] \label{defi:widehat_r}
Given a smooth $M$-th order approximation to $r$, $y \mapsto p^y$, define $\widehat{r} : \mathcal{G} \rightarrow \mathbb{R}$ by 
\[ \widehat{r}(y) = p^y(y) \ . \]
\end{defi}

\begin{thm}[Smoothness of the approximation] \label{thm:smooth}
Assume that 
\[ \nabla r (z) \neq 0 \quad , \quad \forall z \in \Gamma_r \ . \]
For small enough $h$, the corrected approximation term $\widehat{Q} f $ is a smooth function on 
\[ \mathbb{R}^n \setminus \left\{ \widehat{r} = 0 \right\} \ . \]
\end{thm}

\begin{thm}[Approximation of the singularity location $\Gamma_r$]\label{thm:curves_diff}
Assume that 
\[ \nabla r (z) \neq 0 \quad , \quad \forall z \in \Gamma_r \ . \]
Denote  
\[ \left\{ \begin{array}{l} 
\mathcal{C}_1 := \left\{ y \in \mathcal{G} : r(y) = 0 \right\} \\
\mathcal{C}_2 := \left\{ y \in \mathcal{G} : \widehat{r}(y) = 0 \right\} \end{array} \right. \ , \]
then 
\[ d_{H}(\mathcal{C}_1 , \mathcal{C}_2) = O(h^{M+1}) , \]
where $d_{H}$ is the Hausdorff distance of the two sets.
\end{thm}

\begin{rmk} \label{rmk:exp_decay}
Although the results in this section were proven for quasi-interpolations with basis functions of finite support, they may also be proven for quasi-interpolations with basis functions of exponential decay. In the Numerical results (Section \ref{sec:numeric}), we have used MLS quasi-interpolation with weight function of exponential decay.
\end{rmk}

\section{Proofs} \label{sec:proofs}

\subsection{Proof of \cref{thm:G_bound}}
\begin{pf}
Pick \[y \in \mathcal{G}((R+1) h) = \bigcup_{z \in \Gamma_r} B_{(R+1) h}(z) \ . \]
Then, there exists $z \in \Gamma_r$ with $\Vert y - z \Vert < (R+1) h$.
By our assumption $\nabla r(z) \neq 0$. Using Taylor's approximation we can see that for any point 
\[ u \in B_{\Upsilon_M h} \left(z + (\Upsilon_M +1) h \cdot \dfrac{\nabla r (z)}{\Vert \nabla r(z) \Vert} \right) \] we have  
\[ r(u) = r(z) + \langle \nabla r (z), u-z \rangle + O(h^2) \geq h \Vert \nabla r (z) \Vert + O(h^2) \ . \]
Hence, 
\[ B_{\Upsilon_M h} \left(z + (\Upsilon_M +1) h \cdot \dfrac{\nabla r (z)}{\Vert \nabla r(z) \Vert} \right) \cap X \subset \mathcal{P} \ . \]
Similarly we may show that 
\[ B_{\Upsilon_M h} \left(z - (\Upsilon_M +1) h \cdot \dfrac{\nabla r (z)}{\Vert \nabla r(z) \Vert} \right) \cap X \subset X \setminus \mathcal{P} \ . \]
The intersections of each of the balls $B_{\Upsilon_M h} \left(z \pm (\Upsilon_M +1) h \cdot \dfrac{\nabla r (z)}{\Vert \nabla r(z) \Vert} \right)$ with $X$ must be uni-solvent for $\Pi_M(\mathbb{R}^n)$, thus both 
\[ B_{(R + 2 \Upsilon_M + 2 )h} (y) \cap \mathcal{P} \quad \mbox {and} \quad B_{(R + 2 \Upsilon_M + 2)h} (y) \cap (X \setminus \mathcal{P}) \]
are uni-solvent for $\Pi_M(\mathbb{R}^n)$ and $y \in \mathcal{G}$.
\end{pf}

\subsection{Proof of \cref{thm:good_correction}}
\begin{pf}
If $y \in \mathbb{R}^n \setminus \mathcal{G}$ then 
\[B_{(R+1) h}(y) \cap \Gamma_r = \emptyset \ , \]
thus the restriction of $f$ to $B_{R h}(y)$ is smooth and 
\[ \left| \widehat{E}f (y) \right| = \left| E f(y) \right| \leq C_1 \cdot \Vert f \Vert_{C^{M+1}} \cdot h^{M+1} \leq C_1 \cdot ( \Vert g \Vert_{C^{M+1}} + \Vert r \Vert_{C^{M+1}} ) \cdot h^{M+1} \ . \]
Otherwise, if $y \in \mathcal{G}$, then by our assumptions, 
\[ \vert p^{y}(u) - r(u) \vert \leq C_2 \cdot \left( \Vert g \Vert_{C^{M+1}} + \Vert r \Vert_{C^{M+1}} \right) \cdot h^{M+1} \ , \ \forall u \in B_{R h}(y) \ . \]
Consequently, we get, 
\[ \vert E (r_+)(y) - E (p^y)_+ (y) \vert \leq (1+L_1) \cdot C_2 \cdot \left( \Vert g \Vert_{C^{M+1}} + \Vert r \Vert_{C^{M+1}} \right) \cdot h^{M+1} \ . \]
Therefore,  
\begin{align*}
\vert \widehat{E} f( y) \vert & = \vert f(y) - \widehat{Q} f(y) \vert \\
& = \vert \underbrace{f(y) - Q f(y)}_{E f(y) } - E \left( (p^y)_+ \right) (y) \vert \\
& = \vert E f(y)- E \left( (p^y)_+ \right) (y) \vert \\ 
& \leq \vert E r_+(y) - E \left((p^y)_+\right)(y) \vert + \vert E g (y) \vert \\
& \leq  (1+L_1) \cdot C_2 \cdot \left( \Vert g \Vert_{C^{M+1}} + \Vert r \Vert_{C^{M+1}} \right) \cdot h^{M+1} + C_1 \Vert g \Vert_{C^{M+1}} h^{M+1} \ .
\end{align*} \qed
\end{pf}

\subsection{Proof of \cref{thm:smooth}}
\begin{pf}
For $y \in int(\mathbb{R}^n \setminus \mathcal{G})$, there exists $\delta > 0$ such that within the ball 
\[B_\delta (y) \subset \mathbb{R}^n \setminus \mathcal{G} \] the correction term is equal to 
\[ \widehat{Q}f = Qf \ , \]
which is a smooth function. 

Similarly, for $y \in int(\mathcal{G}) \setminus \left\{ \widehat{r} = 0 \right\}$, there exists $\delta > 0$ such that 
\[ \forall v \in B_{\delta} (y) \subset \mathcal{G} \setminus \left\{ \widehat{r} = 0 \right\} \ , \] the correction term is equal to 
\[ \widehat{Q} f(v) = Qf(v) + E((p^v)_+)(v) = Qf(v) + (p^v)_+(v) - Q((p^v)_+)(v) \ . \]
The function $Qf$ is obviously smooth.
Likewise, the smoothness of $Q((p^v)_+)(v)$ follows from the smoothness of the mapping $v \mapsto p^v$.
Last, since 
\[B_{\delta} (y) \cap \left\{ \widehat{r} = 0 \right\} = \emptyset \ , \] then $(p^v)_+(v)$ is also smooth.

We still have to show that the corrected term is smooth for $y \in \partial(\mathcal{G})$.
Since $\mathcal{G}((R+1)h) \subset \mathcal{G}$, we must have 
\[ B_{(R+1)h} (y) \cap \Gamma_r = \emptyset \ . \]
Specifically, for $v \in B_{\frac{1}{2}h} (y)$ we get 
\[ B_{(R+\frac{1}{2})h}(v) \cap \Gamma_r = \emptyset \ . \]
By our assumptions, 
\[ \nabla r(z) \neq 0 \ , \ \forall z \in \Gamma_r \ , \]
hence, since $y \mapsto p^y$ is a smooth $M$-th order approximation of $r$, we get for small enough $h$ that 
\[ B_{Rh}(v) \cap \left\{ p^v = 0 \right\} = \emptyset \ . \]
That is, $p^v$ does not change sign in $B_{Rh} (v)$, and $E((p^v)_+)(v) = 0$, which gives us 
\[\widehat{Q}f(v) = Qf(v) \ , \ \forall v \in B_{\frac{1}{2} h}(y) \ . \] 
Thus the correction term is smooth at $y$. \qed

\end{pf}

\subsection{Proof of \cref{thm:curves_diff}}
\begin{pf}
Pick $y \in \mathcal{C}_1$, then $r(y) = 0$ and $\nabla r(y) \neq 0$.
Using Taylor's approximation we have 
\begin{equation} \label{eq:r(y)+eps}
r(y \pm \epsilon \nabla r (y)) = \pm \epsilon \Vert \nabla r(y) \Vert^2 + O(\epsilon^2 \cdot \Vert \nabla r (y) \Vert^2) \ . 
\end{equation}
Hence, \[ r(y-\epsilon \nabla r(y)) < 0 < r(y+\epsilon \nabla r(y)) \ .\]
Since $y \mapsto p^y$ is a smooth $M$-th order approximation to $r$, the function $\widehat{r}(y) = p^y(y)$ must be smooth and satisfy 
\begin{equation} \label{eq:wide_r-r} 
\widehat{r}(y) - r(y)  = O(h^{M+1}) \ . 
\end{equation}
hence 
\[ \widehat{r}(y-\epsilon \nabla r(y)) < 0 < \widehat{r}(y+\epsilon \nabla r(y)) \ .\]
Then there must exist $u = y + \lambda \nabla r(y)$ with $\vert \lambda \vert < \epsilon$ such that $\widehat{r}(u) = 0$.
However, from (\ref{eq:r(y)+eps}) and (\ref{eq:wide_r-r}) we get 
\[ O(h^{M+1}) = r(u) = \lambda \Vert \nabla r (y) \Vert^2 + O( \lambda^2 \Vert \nabla r(y) \Vert^2 ) \ . \] 
Hence, 
\[ \lambda \Vert \nabla r (y) \Vert^2 = O(h^{M+1}) \ , \]
and consequently
\[ \Vert y - u \Vert = \vert \lambda \vert \Vert \nabla r (y) \Vert = O(h^{M+1}) \ . \]
That is, there must exist $u$ with $\Vert y - u \Vert = O(h^{M+1})$ and $\widehat{r}(u) = 0$, hence, 
\[ \sup_{y \in \mathcal{C}_1} \inf_{u \in \mathcal{C}_2} d(y, u) = O(h^{M+1}) \ . \]
Similarly we show that 
\[ \sup_{y \in \mathcal{C}_2} \inf_{u \in \mathcal{C}_1} d(y, u) = O(h^{M+1}) \ , \]
and we have \[d_H(\mathcal{C}_1, \mathcal{C}_2) = O(h^{M+1}) \ . \] \qed
\end{pf}

\section{Partitioning the data points in $X$ with respect to the sign of $r$} \label{sec:signs}
Our method relies upon our ability to correctly identify a set $\mathcal{P} \subset X$ which partitions the data points in $X$ with respect to the signs of the function $r$. That is, 
\begin{enumerate}
\item For all $x \in \mathcal{P}$ either $r(x) > 0$ or $r(x) = O(h^{M+1})$.
\item For all $x \in X \setminus \mathcal{P}$ either $r(x) < 0$ or $r(x) = O(h^{M+1})$.
\end{enumerate}
We propose an algorithm for building the set $\mathcal{P}$ according to the following steps: 
\begin{enumerate}[label={Step (\arabic*)}]
\item \label{step:S} Find a set $\mathcal{S} \subset X$ satisfying $\Gamma_r \subset \bigcup\limits_{x \in \mathcal{S}} B_{3h} (x) \ . $
\item \label{step:components} Denote by $A_1, \ldots, A_k$ the connected components of 
\[ \mathbb{R}^n \setminus \left( \bigcup\limits_{x \in \mathcal{S}} B_{3h} (x) \right) \subset \mathbb{R}^n \setminus \Gamma_r \ . \] 
Note that for all $1 \leq i \leq k$ and $x_1, x_2 \in A_i$ we have $r(x_1) \cdot r(x_2) > 0$.
\item \label{step:a} Define a function $a : X \rightarrow \left\{ 1, \ldots , k \right\}$ satisfying 
\[ \forall x \in X : a(x) = i \quad \Longrightarrow \quad \begin{array}{ll} \mbox{either} & r(x) = O(h^{M+1}) \\ \mbox{or} & r(x) \cdot r(u) > 0 \ , \ \forall u \in A_i \cap X \end{array} \ . \]
For $1 \leq i \leq k$ set $\bar{A}_i := \left\{ x \in X : a(x) = i \right\} \ . $
\item \label{step:sigma} Define a function $\sigma : \left\{ 1 , \ldots , k \right\} \rightarrow \left\{ 1 , 2 \right\}$ satisfying 
\[ \sigma(i_1) = \sigma(i_2) \quad \Longrightarrow \quad \forall x_1 \in A_{i_1} \ , \ \forall x_2 \in A_{i_2} : r(x_1) \cdot r(x_2) > 0  \ . \]
\item \label{step:P} Set $\mathcal{P} := \bigcup\limits_{\sigma(i) = 1} \bar{A}_i \ . $

Refine the set $\mathcal{P}$.
\end{enumerate}

\begin{algorithm}[H]
\caption{Build the set $\mathcal{S}$ satisfying $\Gamma_r \subset \bigcup\limits_{x \in \mathcal{S}} B_{3h} (x)$}
\label{algo:S}
\begin{algorithmic}[1]
\State Define $\Psi : \mathbb{R}^n \rightarrow \Pi_M(\mathbb{R}^n)$ by 
\[ \Psi(y) := \operatorname*{arg\,min}\limits_{p \in \Pi_M(\mathbb{R}^n)} \sum\limits_{x \in X} \omega \left( \frac{\Vert x-y \Vert }{h}\right) \cdot \left( p(x) - f (x) \right)^2 \ . \]
Here $\omega$ is a smooth weight function with compact support 
\[ supp(\omega) = [ 0 , \rho ] \supset [0, \Upsilon_M] \ . \]
\State For $u \in X$ define $\epsilon_u := \max \left\{ \vert (\Psi(u))(x) - f(x) \vert : x \in B_{\rho h}(u) \cap X \right\}$.
\State Set \[ \mathcal{S} := \left\{ x \in X \left| \exists u \in B_{(\rho + 5)h} (x) \cap X : \vert (\Psi(u))(x) - f(x) \vert > \epsilon_u \cdot \left( \frac{\rho + 5 }{\rho} \right)^{M+1} \right. \right\} \ . \]
\end{algorithmic}
\end{algorithm}

\begin{algorithm}[H]
\caption{Find the connected components of $\mathbb{R}^n \setminus \left( \bigcup\limits_{x \in \mathcal{S}} B_{3h} (x) \right)$}
\label{algo:components}
\begin{algorithmic}[1]
\State Build a graph $G = (V,E)$ with vertices 
\[ V := X \setminus \left( \bigcup\limits_{x \in \mathcal{S}} B_{3h}(x) \right) \ , \]
and edges 
\[ E := \left\{ (u,v) \in V \left| d(u,v) < 2h \right. \right\} \ . \]
\State Set $A_1 , A_2 , \ldots , A_k$ to be the connected components of the graph $G$.
\end{algorithmic}
\end{algorithm}

\begin{algorithm}[H]
\caption{Compute the function $a : X \rightarrow \left\{ 1 , \ldots , k \right\}$ satisfying : \\ $ a(x) = i \quad \Rightarrow \quad $
either $ r(x) = O(h^{M+1}) $ or $ r(x) \cdot r(u) > 0 \ , \ \forall u \in A_i \cap X $ .}
\label{algo:a}
\begin{algorithmic}[1]
\State Let $\omega$ be a smooth weight function with compact support.
\State For $1 \leq i \leq k$ denote 
\[ \Omega_i = \left\{ y \in \mathbb{R}^n : supp\left( \omega \left( \frac{ \Vert \cdot - y \Vert }{h} \right) \right) \cap A_i \cap X \mbox{ is uni-solvent for } \Pi_M(\mathbb{R}^n) \right\} \ , \]
and define $\Theta_i : \Omega_i \rightarrow \Pi_M(\mathbb{R}^n)$ by 
\[ \Theta_i(y) := \operatorname*{arg\,min}\limits_{p \in \Pi_M(\mathbb{R}^n)} \sum\limits_{x \in A_i \cap X} \omega \left( \frac{\Vert x-y \Vert }{h}\right) \cdot \left( p(x) - f (x) \right)^2  \ . \]
\State For $x \in X$ define 
\[ a(x) := \left\{ \begin{array}{ll} i & , x \in A_i \\ 
\operatorname*{arg\,min}\limits_{\left\{ 1 \leq i \leq k : x \in \Omega_i \right\} } \left| (\Theta_i(x))(x) - f(x) \right| & , x \in X \setminus \left( \bigcup\limits_{i=1}^k A_i \right) \end{array} \right. \ . \]
\end{algorithmic}
\end{algorithm}

\begin{algorithm}[H]
\caption{Compute the function $\sigma : \left\{ 1 , \ldots , k \right\} \rightarrow \left\{ 1 , 2 \right\} $ satisfying : \\ $ \sigma(i_1) = \sigma(i_2) \quad \Rightarrow \quad \forall x_1 \in A_{i_1} \ , \ \forall x_2 \in A_{i_2} : r(x_1) \cdot r(x_2) > 0  $ .}
\label{algo:sigma}
\begin{algorithmic}[1]
\State Define $\sigma^0 : \left\{ 1 , \ldots , k \right\} \rightarrow \left\{ 1 , \ldots , k \right\}$ by $\sigma^0(i) = i$.
\State For $1 \leq l \leq k$ define $X^0_l = \bar{A}_l$. 
\For{$j = 0 , \ldots , k-3$}
\State For each $1 \leq l \leq k-j$ denote $X^{j+1}_l := \bigcup\limits_{\sigma^j(i) = l} X^j_i $.
\State For each $1 \leq l \leq k-j$ define $\Phi_l : \mathbb{R}^n \rightarrow \Pi_M(\mathbb{R}^n)$ by 
\[ \Phi_l(y) := \operatorname*{arg\,min}\limits_{p \in \Pi_M(\mathbb{R}^n)} \sum\limits_{x \in X^{j+1}_l} \omega \left( \frac{\Vert x-y \Vert }{h}\right) \cdot \left( p(x) - f (x) \right)^2 \ . \]
Here $\omega$ is a smooth weight function with compact support $supp(\omega) = [0, \rho]$. 
\For{$1 \leq l_1 , l_2 \leq k-j$}
\State $N \leftarrow \left\{ x \in X^{j+1}_{l_1} : B_{\rho h}(x) \cap X^{j+1}_{l_2} \neq \emptyset \right\}$
\If{$\# N == 0$}
\State $D_{l_1 , l_2} \leftarrow \infty$
\Else
\State $D_{l_1 , l_2} \leftarrow \max \left\{ \left| (\Phi_{l_2} (x))(x) - f(x) \right| : x \in N \right\}$
\EndIf
\EndFor
\State Pick $1 \leq l_1 <  l_2 \leq k-j$ for which $(D+D^T)_{l_1 , l_2}$ is minimal.
\State Define $\sigma^{j+1} : \left\{ 1 , \ldots , k-j \right\} \rightarrow \left\{ 1 , \ldots , k-j-1 \right\}$ by 
\[ \sigma^{j+1}(i) = \left\{ \begin{array}{ll} 
i , & 1 \leq i \leq l_2 - 1 \\
l_1, & i = l_2 \\
i-1, & l_2 < i \leq k-j \end{array} \right. \ . \]
\EndFor
\State Set $\sigma = \sigma^{k-2} \circ \sigma^{k-3} \circ \ldots \circ \sigma^0 \ . $
\end{algorithmic}
\end{algorithm}

\begin{algorithm}[H]
\caption{Refinement of the set $\mathcal{P}$}
\label{algo:refine}
\begin{algorithmic}[1]
\State $B \leftarrow \left\{ x \in X : \exists x_1 \in \mathcal{P} \ , \ x_2 \in X \setminus \mathcal{P} \mbox{ s.t. } d(x, x_1) , d(x, x_2) < 2h \right\}$
\State $O_1 \leftarrow \mathcal{P} \setminus B$
\State $O_2 \leftarrow X \setminus  (\mathcal{P} \cup B)$
\State For $k = 1,2$ define $\Xi_k : \mathbb{R}^n \rightarrow \Pi_M(\mathbb{R}^n)$ by 
\[ \Xi_k(y) := \operatorname*{arg\,min}\limits_{p \in \Pi_M(\mathbb{R}^n)} \sum\limits_{x \in O_k} \omega \left( \frac{\Vert x-y \Vert }{h}\right) \cdot \left( p(x) - f (x) \right)^2 \ . \]
Here $\omega$ is a smooth weight function with compact support.
\State $\mathcal{P} \leftarrow O_1 \cup \left\{ x \in B : \vert (\Xi_1 (x))(x) - f(x) \vert < \vert (\Xi_2 (x))(x) - f(x) \vert \right\}$.
\end{algorithmic}
\end{algorithm}

To show that the algorithm indeed generates the set $\mathcal{P}$, as defined in Definition \ref{defi:P}, we observe the following: 
\begin{enumerate}
\item \label{enum:S} In Step \ref{algo:S}, $\Psi$ is an MLS polynomial approximation of $f$, hence if 
\[B_{\rho h}(u) \cap \Gamma_r = \emptyset \] then the restriction of $f$ to $B_{\rho h}(u)$ is either $g$ or $g+r$. 
W.l.o.g. we will assume that 
\[ r(y) < 0 \ , \ \forall y \in B_{\rho h}(u) \ , \]
hence the restriction of $f$ to $B_{\rho h}(u)$ is the smooth function $g$ and 
\begin{equation} \label{eq:Psi} 
 \vert (\Psi(u))(y) - g(y) \vert = O(h^{M+1}) \ , \ \forall y \in B_{(\rho + 5)h}(u) 
\end{equation}
Specifically, for $x \in B_{\rho h}(u) \cap X$ we have $f(x) = g(x)$, and 
\begin{equation} \label{eq:e_u}
\epsilon_u = \max \left\{ \vert (\Psi(u))(x) - f(x) \vert : x \in B_{\rho h}(u) \cap X \right\} = O(h^{M+1})
\end{equation}
However, if 
\[\exists z \in B_{(\rho + 2)h}(u) \cap \Gamma_r \quad \mbox{with} \quad \nabla r(z) \neq 0 \] then there must exist 
\[x \in B_h \left( z + 2 h \cdot \frac{\nabla r (z) }{\Vert \nabla r (z) \Vert} \right) \cap X \subset B_{(\rho + 5)h}(u) \cap X  \]
for which  
\[r(x) > h \cdot \Vert \nabla r (z) \Vert + O(h^2) > 0 \ . \]
From (\ref{eq:Psi}) we get 
\begin{eqnarray*}
\vert (\Psi(u))(x) - f(x) \vert & = & \vert (\Psi(u))(x) - g(x) - r(x) \vert \\
& > & \vert r(x) \vert - \vert (\Psi(u))(x) - g(x) \vert \\
& = & h \cdot \Vert \nabla r(z) \Vert + O(h^2) + O(h^{M+1}) 
\end{eqnarray*}
For small enough $h$ we get by (\ref{eq:e_u})
\[  \vert (\Psi(u))(x) - f(x) \vert > h \cdot \Vert \nabla r(z) \Vert + O(h^2) > \epsilon_u \cdot \left( \frac{\rho + 5}{\rho} \right)^{M+1} \ . \]
Thus $x \in \mathcal{S}$ and $z \in \bigcup\limits_{x \in \mathcal{S}} B_{3h} (x)$.
\item \label{enum:a} In Step \ref{algo:a}, each $\Theta_i$ is an MLS approximation of $f$ based only upon the data points in $A_i \cap X$. Since $A_i \cap \Gamma_r = \emptyset$, $f$ is smooth on each $A_i$ and $\Theta_i$ is either a polynomial approximation of $g$ or of $g+r$. 

Pick $x \in \mathcal{S}$ and assume that both $x$ and $A_i$ belong to the same connected component of $\mathbb{R}^n \setminus \Gamma_r$, hence 
\[ \left| (\Theta_i (x))(x) - f(x) \right| = O(h^{M+1}) \ . \]
W.l.o.g. assume that $r(x) > 0$. 
If $r(u) < 0$ for $u \in A_{a(x)}$ then $\Theta_{a(x)} (x)$ is a polynomial approximation of $g$, which gives us 
\begin{eqnarray*} 
O(h^{M+1}) &=& \left| (\Theta_i(x)) (x) - f(x) \right| \\
& \geq & \left| (\Theta_{a(x)}(x)) (x) - f(x) \right| \\
& \geq & \left| r(x) \right| - \left| (\Theta_{a(x)})(x) - g(x) \right| = \vert r(x) \vert + O(h^{M+1}) 
\end{eqnarray*}
Thus in this case $r(x) = O(h^{M+1})$. 
\item \label{enum:sigma} In Step \ref{algo:sigma}, the operator $\Phi_l$ returns MLS approximations of either $g$ or $g+r$. At each step of the for-loop we merge two subsets $X^{j+1}_{l_1}$ and $X^{j+1}_{l_2}$ on which the approximations $\Phi_{l_1}$ and $\Phi_{l_2}$ are close. At the end we will have two subsets of $X$, such the the restriction of $f$ to one subset would be $g$ and the restriction of $f$ to the other subset would be $g+r$.
\item \label{enum:refine} In Step \ref{algo:refine} we propose to refine the initial set $\mathcal{P}$. We refine this set by removing from $\mathcal{P}$ and from $X \setminus \mathcal{P}$ data points that are close to the boundary. Then, we add to $\mathcal{P}$ only the boundary data points for which the MLS approximation based upon $\mathcal{P}$ has smaller errors than the MLS approximation based upon $X \setminus \mathcal{P}$. 
\end{enumerate}
The above observations prove that the set $\mathcal{P}$ returned by our proposed algorithm is a partition of the data points in $X$ with respect for the sign of the function $r$. 
\begin{rmk} \label{rem:chice}
In Step \ref{step:P} of the algorithm we arbitrarily choose the set $\mathcal{P}$, thus we might choose the set on which $r$ returns negative values. However, this choice has no effect on the approximation algorithm. Indeed, 
\[E r_+ = E (r-r_-) = \underbrace{E r}_{O(h^{M+1})} + E \underbrace{-(r_-)}_{(-r)_+} = O(h^{M+1}) + E(-r)_+ \ , \]
where $E$ is the approximation error of the quasi-interpolation operator.
Hence, the initial choice is insignificant.
\end{rmk}

\begin{figure*}
\begin{subfigure}{0.5\textwidth}
\includegraphics[width=\textwidth]{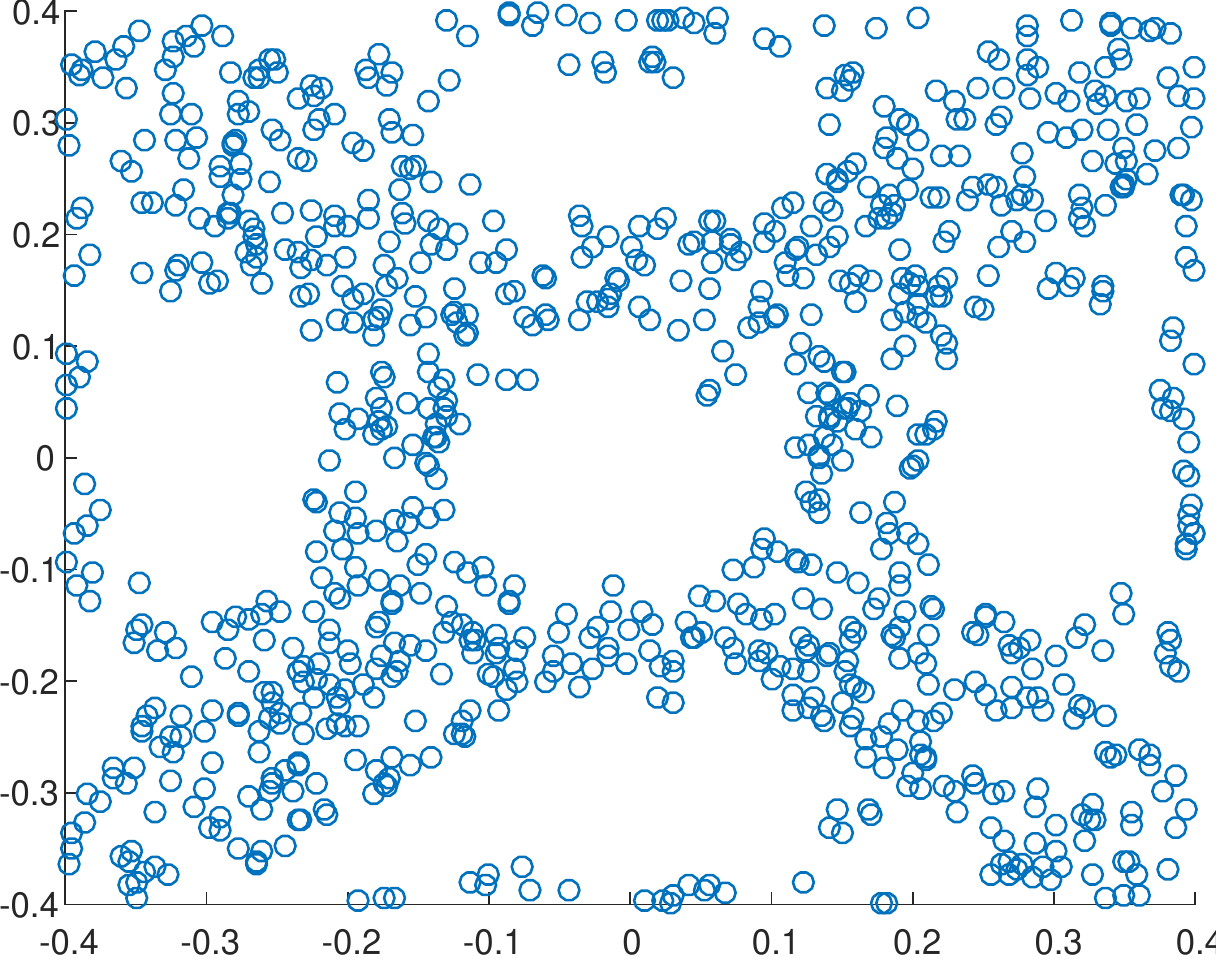}
\caption{$\mathcal{S}$}
\label{fig:S_0}
\end{subfigure}
\begin{subfigure}{0.5\textwidth}
\includegraphics[width=\textwidth]{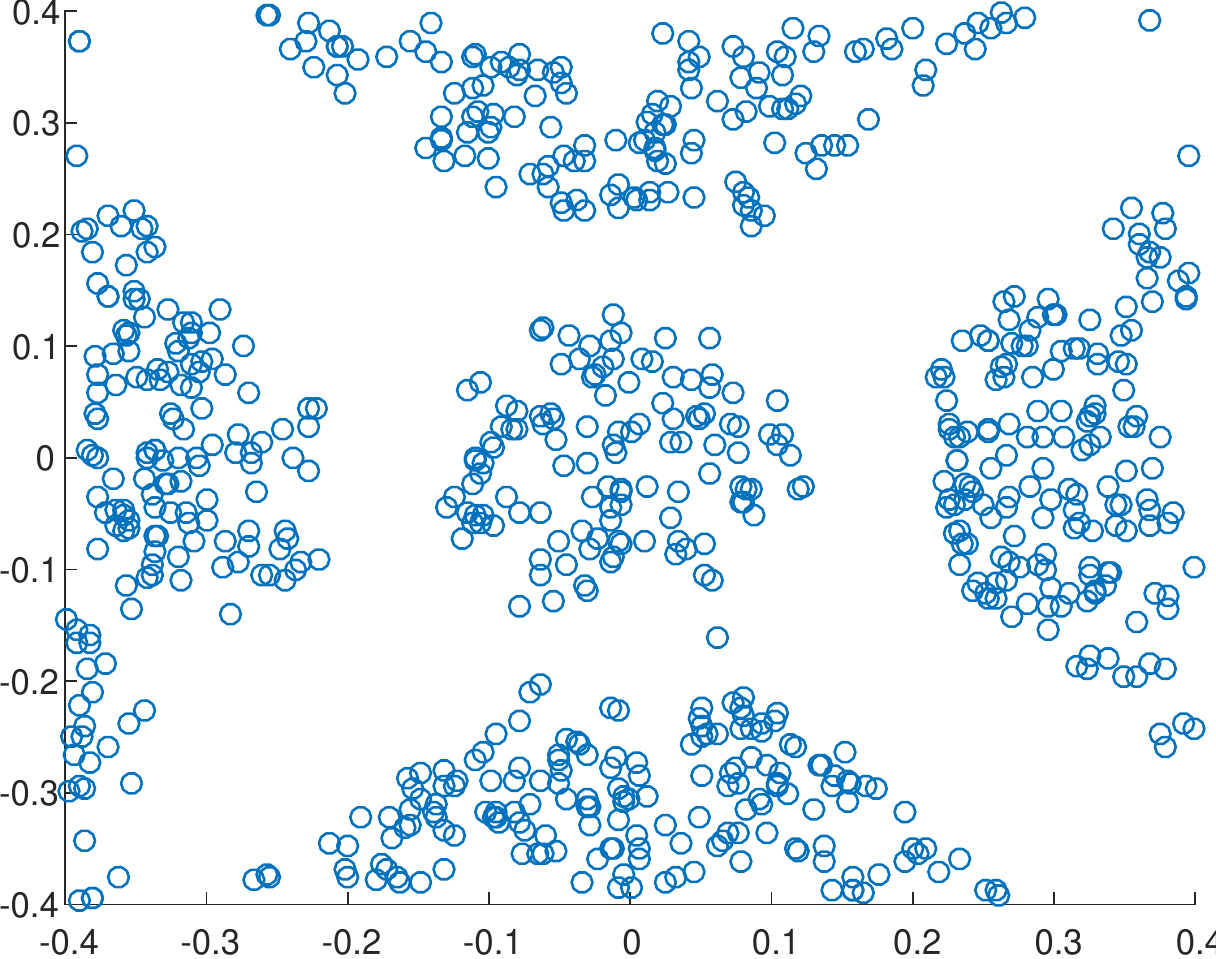}
\caption{$X \setminus \mathcal{S}$}
\label{fig:X-S_0}
\end{subfigure}
\caption{The sets $\mathcal{S}$ and $X \setminus \mathcal{S}$.} 
\label{fig:S_X-S}
\end{figure*}

\begin{figure*}
\begin{subfigure}{0.5\textwidth}
\includegraphics[width=\textwidth]{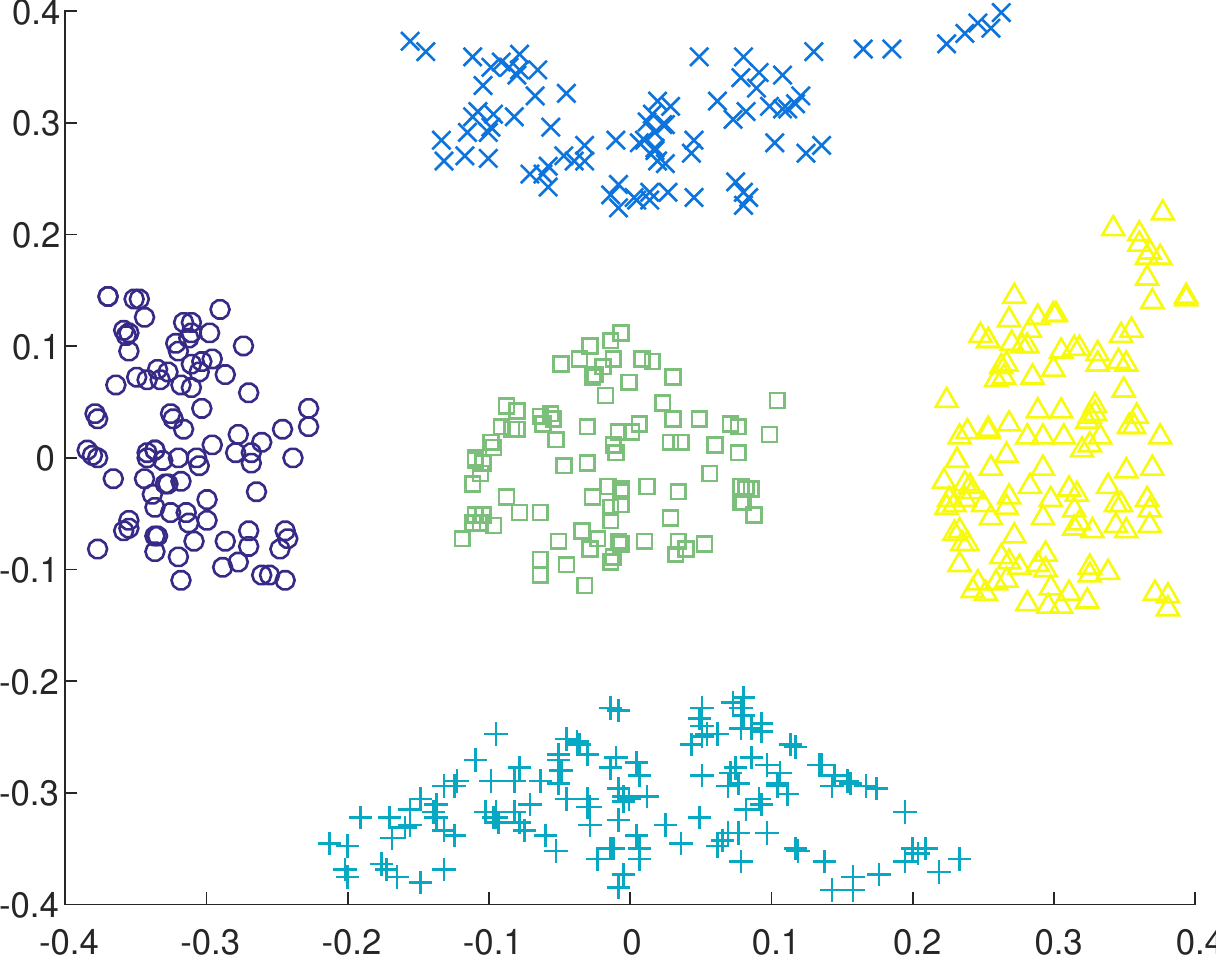}
\caption{The connected components of $X \setminus \mathcal{S}$}
\label{fig:components_X-S}
\end{subfigure}
\begin{subfigure}{0.5\textwidth}
\includegraphics[width=\textwidth]{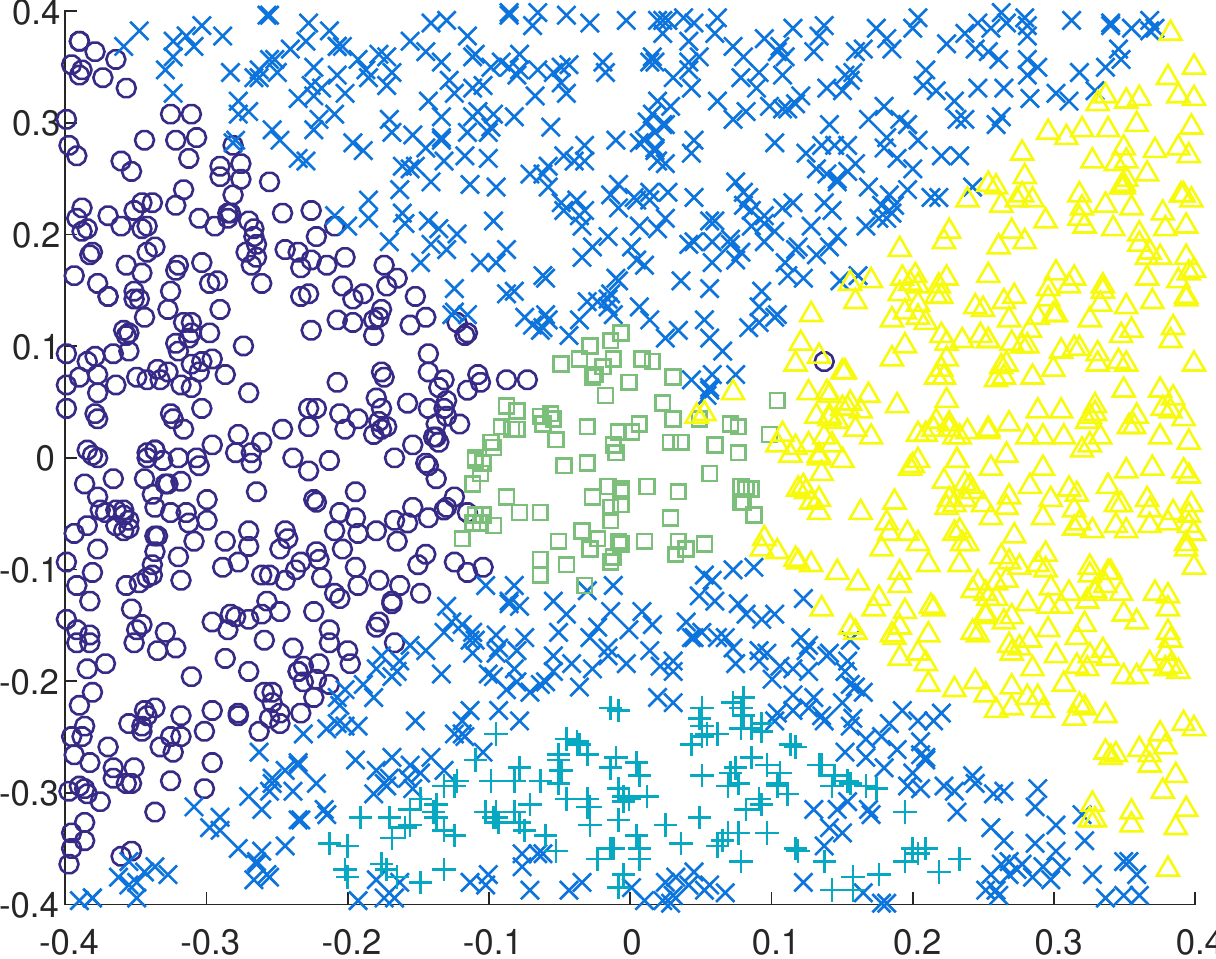}
\caption{The connected components of $X \setminus \Gamma_r $}
\label{fig:components_X-Gamma}
\end{subfigure}
\caption{The connected components of $X \setminus \mathcal{S}$ and of $X \setminus \Gamma_r$.} 
\label{fig:con_components}
\end{figure*}

\begin{figure*}
\begin{subfigure}{0.5\textwidth}
\includegraphics[width=\textwidth]{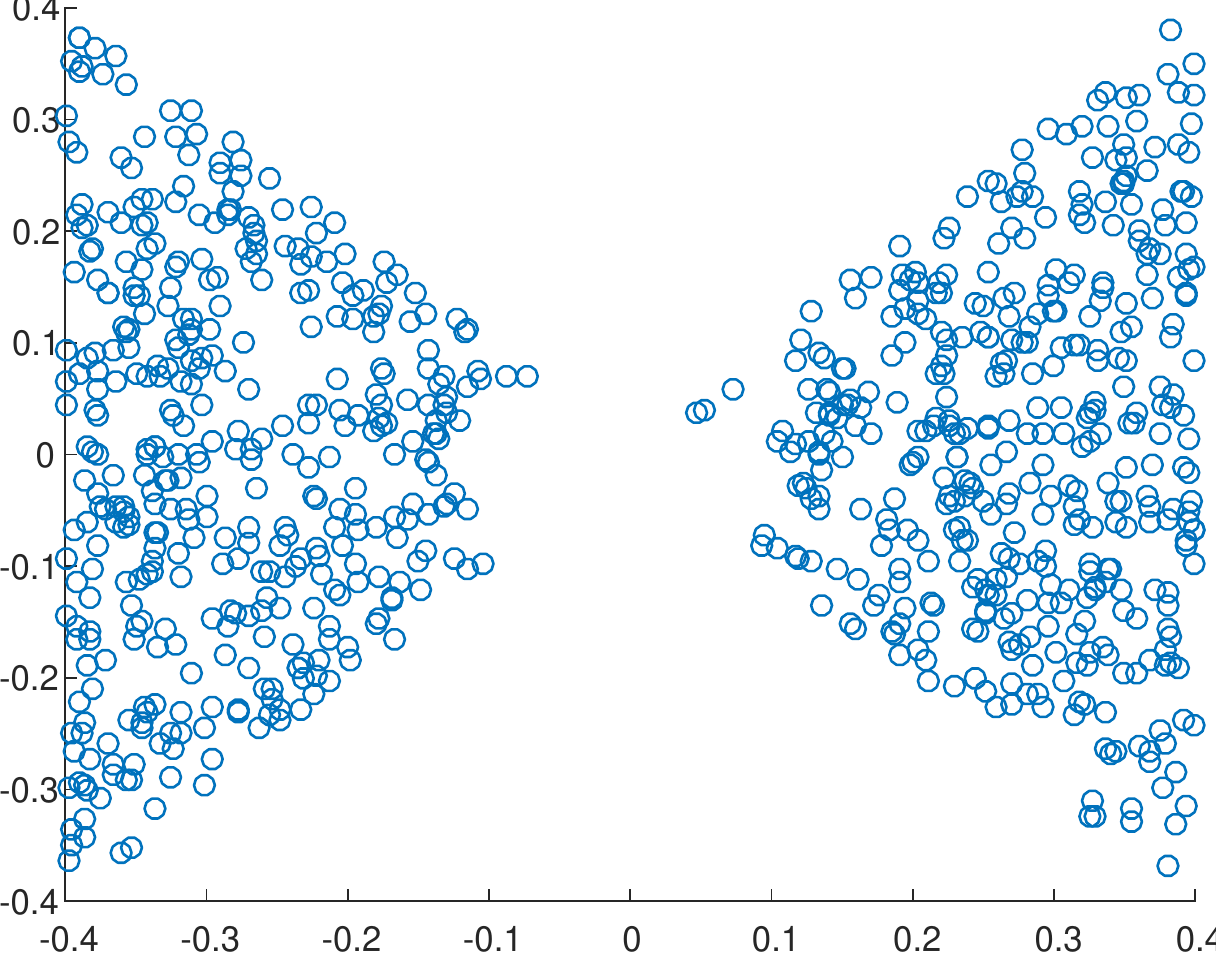}
\caption{Initial $\mathcal{P}$}
\label{fig:initial_P}
\end{subfigure}
\begin{subfigure}{0.5\textwidth}
\includegraphics[width=\textwidth]{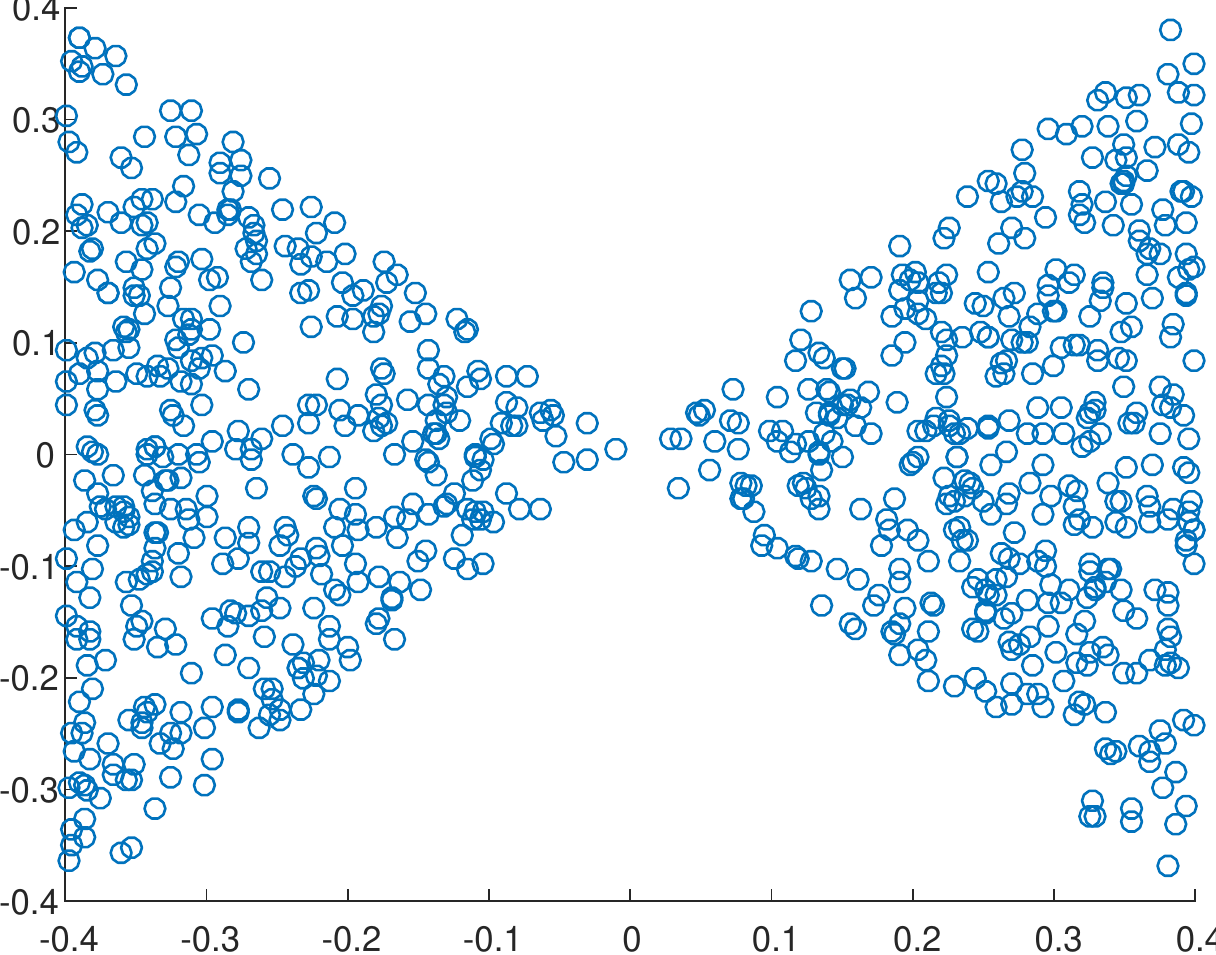}
\caption{Final $\mathcal{P}$}
\label{fig:final_P}
\end{subfigure}
\caption{The set $\mathcal{P}$ before and after the refinement.} 
\label{fig:P}
\end{figure*}

For example, we ran the partitioning algorithm on the function 
\[f(x,y) = ((x+y) \cdot (x-y))_+ \ . \]
In Figure \ref{fig:S_X-S} one can see the sets $\mathcal{S}$ and $X \setminus \mathcal{S}$ (see Step \ref{algo:S}). 
In Figure \ref{fig:con_components} one can see the connected components of $X \setminus \mathcal{S}$ (Step \ref{algo:components}), and the connected components of $X \setminus \Gamma_r$ (see Step \ref{algo:a}).
In Figure \ref{fig:P} one can see the initial set $\mathcal{P}$ and the final $\mathcal{P}$ after the refinement (see Step \ref{algo:refine}).

\begin{rmk}\label{rem:domain}
Although the sign determination algorithm, is described for $f$ defined on $\mathbb{R}^n$, the algorithm can also be applied to a compact  domain. Moreover, since the computation of the corrected approximation $\widehat{Q}f$, is a local procedure, it might be computationally preferable to break the domain into smaller compact subsets. 
\end{rmk}

\section{Approximation of the signed function $r$} \label{sec:p^y}
In Section \ref{sec:results} we have introduced the notion of a smooth $M$-th order approximation to the function $r$. In this section we will introduce two methods we may use to find a smooth $M$-th order approximation to $r$.
Recall that a smooth $M$-th order approximation is a mapping $y \mapsto p^y$ satisfying : 
\begin{itemize}
\item There exists $C_2 > 0$ such that for all $y \in \mathcal{G}$, 
\[ \vert p^y(u) - r(u) \vert \leq C_2 \cdot \left( \Vert g \Vert_{C^{M+1}} + \Vert r \Vert_{C^{M+1}}\right) \cdot h^{M+1} \ , \ \forall u \in B_{R h}(y) \ . \]
\item Let $\left\{ p_\alpha \right\}_{\alpha \in I}$ be a polynomial basis of $\Pi_M(\mathbb{R}^n)$ and write 
\[ p^y = \sum_{\alpha \in I} \lambda_\alpha (y) \cdot p_\alpha \ . \]
Then, the mappings $y \mapsto \lambda_{\alpha}(y)$ are infinitely smooth.
\end{itemize}
The methods that we introduce are both based upon MLS, hence the smoothness is trivial. As for the first condition, while the first method gives slightly better approximations to $r$ (as can be seen in section \ref{sec:numeric}), it also depends upon a non-singularity conjecture (which can be verified numerically). The first method, described in subsection \ref{sub:p^y_1}, is in-fact a generalization of the method suggested by Lipman and Levin \cite{Lipman} for the univariate case. Here also we analyse the quasi-interpolation errors, and find a polynomial $p^y_1$ which gives the best simulation for these errors. 
The second method, described in subsection \ref{sub:p^y_2} utilizes two MLS approximations and is independent upon a non-singularity condition. One MLS approximation is based upon the data set $\mathcal{P}$ and the other is based upon its complement $X \setminus \mathcal{P}$. The difference between the two approximations would be the polynomial approximation to $r$.

\subsection{Approximation by error analysis} \label{sub:p^y_1}
This method derives the polynomial approximation to $r$ from an analysis of the quasi-interpolation errors. 
\begin{defi}[The approximant $p^y_1$] \label{def:p^y_1}
For $y \in \mathbb{R}^n$ define $p^{y}_1 \in \Pi_M(\mathbb{R}^n)$ by 
\[ p^y_1 := \operatorname*{arg\,min}\limits_{p \in \Pi_M(\mathbb{R}^n)} \sum\limits_{x \in X} \omega \left( \dfrac{\Vert y-x \Vert}{h} \right) \cdot \left( E (\chi_{\mathcal{P}} \cdot p)(x) - E f (x) \right)^2 \ . \]
Here $\omega : \mathbb{R}_+ \rightarrow \mathbb{R}_+$ is an infinitely smooth positive weight function with compact support, 
\[ supp(\omega) \supset [0, (R + 2 \Upsilon_M + 2)h] \] 
and $\chi_{\mathcal{P}}$ is the indicator function defined by 
\[ \chi_{\mathcal{P}} (z) = \left\{ \begin{array}{ll} 
1, & z \in \mathcal{P} \\
0, & z \notin \mathcal{P} \end{array} \right. \ , \ \forall z \in \mathbb{R}^n \ . \]
\end{defi}
\begin{rmk}
A more natural choice for the polynomial $p^y_1$ would have been the polynomial minimizing the sum
\[\sum\limits_{x \in X} \omega \left( \dfrac{\Vert y-x \Vert}{h} \right) \cdot \left( E (p_+)(x) - E f (x) \right)^2 \ . \]
However, this computation is not linear.
\end{rmk}

\begin{thm}[$p^y_1$ is a smooth $M$-th order approximation to $r$] \label{thm:p^y_1}
If the vectors 
\[\left\{ E (\chi_{\mathcal{P}} \cdot p_{\alpha})(x) : x \in \Omega(y) \cap X \right\}_{\alpha \in I}\]
are linearly independent for all $y \in \mathcal{G}$, with 
\[ \Omega(y) = \left\{ y + z : \dfrac{\Vert z \Vert}{h} \in supp(\omega) \right\} \ , \]
then the mapping $y \mapsto p^y_1$ is a smooth $M$-th order approximation to $r$.
\end{thm}

\begin{pf}
Let us define an MLS operator $\tilde{Q}$ for a function $\phi : \mathbb{R}^n \rightarrow \mathbb{R}$ and $y \in \mathbb{R}^n$ by 
\[ \tilde{Q} \phi (y)  =  p_{\phi,y}(y) \]
where 
\[p_{\phi , y} := \operatorname*{arg\,min}\limits_{p \in \Pi_M(\mathbb{R}^n)} \sum\limits_{x \in X} \omega \left( \dfrac{\Vert y-x \Vert }{h} \right) \cdot \left( E (\chi_{\mathcal{P}} \cdot p)(x) - E (\chi_{\mathcal{P}} \cdot \phi) (x) \right)^2 \ . \]
Note the difference between the above definition and the original MLS definition (\ref{eq:MLS}). In the original MLS definition we sum the squares of the difference between the values of $\phi$ and of the approximating polynomial $p$ at the data points, while in this definition we sum the squares of the differences between the approximation errors of $\chi_{\mathcal{P}} \cdot \phi$ and $\chi_{\mathcal{P}} \cdot p$ at the data points.
Write \[ p = \sum_{\alpha \in I} \delta_\alpha p_\alpha \quad \mbox{and} \quad p_{\phi, y} = \sum_{\alpha \in I} (\delta_{\phi, y})_{\alpha} p_\alpha \ , \]
then from the linearity of $E$ we have
\[\pmb{\delta_{\phi , y}} = \operatorname*{arg\,min}\limits_{\pmb{\delta}} \sum\limits_{x \in X} 
\omega \left( \dfrac{ \Vert y-x \Vert }{h} \right) \cdot 
\left( \sum\limits_{\alpha} \delta_\alpha E (\chi_{\mathcal{P}} \cdot p_\alpha )(x) - E (\chi_{\mathcal{P}} \cdot \phi) (x) \right)^2 \ . \]
To solve this problem we follow \cite{levin1998approximation}, from which we know that if the vectors 
\[ \left\{ E\left(\chi_{\mathcal{P}} \cdot p_\alpha \right)(x) : x \in \Omega(y) \cap X \right\}_{\alpha \in I} \] are linearly independent then the solution to this problem is the vector $\pmb{\delta_{\phi , y}}$ defined by 
\begin{equation} \label{eq:delta}
\pmb{\delta_{\phi , y}} = (ADA^T)^{-1} AD \pmb{\phi} 
\end{equation}
where $A$ is the matrix 
\[A = \left( E \left( \chi_{\mathcal{P}} \cdot p_{\alpha} \right) (x) \right)_{\alpha, x} \ , \ \forall \alpha \in I \ , \ \forall x \in X \]
$D$ is the diagonal matrix with values
\[ D_{x , x} = \omega \left( \frac{\Vert y - x \Vert}{h} \right) \ , \forall x \in X \]
and $\pmb{\phi}$ is the vector 
\[\phi_{x} = E (\chi_{\mathcal{P}} \cdot \phi)(x) \ , \ \forall x \in X \ . \]
Moreover, the operator $\tilde{Q}$ clearly reproduces polynomials in $\Pi_M(\mathbb{R}^n)$ and as such is a quasi-interpolation operator
\begin{equation} \label{eq:tildeQ}
\tilde{Q} \phi (y) = \sum_{x \in X} \tilde{q}_{x} (y) E (\chi_{\mathcal{P}} \cdot \phi) (x) 
\end{equation}
with basis functions 
\begin{eqnarray} \label{eq:basis}
\tilde{q}_{x}(y) := \sum\limits_{\alpha \in I} p_{\alpha}\left(y \right) \cdot \left( (ADA^T)^{-1} AD \right)_{\alpha,x} 
\end{eqnarray}
of compact support
\[ R_2 := \inf \left\{ \rho : supp(\tilde{q}_{x}) \subseteq B_{\rho h}(x) \ , \ \forall x \in X \right\} \ . \]
Note that 
\begin{equation} \label{eq:p^y_1}
p^y_1 (z) = \sum\limits_{x \in X} \tilde{q}_x (z) \cdot Ef(x) \ . 
\end{equation}
Let us show that $\tilde{Q}$ has a bounded Lebesgue constant, 
\[ L_2 := \sup \left\{ \sum_{x \in X} \vert \tilde{q}_{x}(z) \vert : z \in \mathbb{R}^n \right\} \ . \]
Assume that the polynomials $\left\{ p_{\alpha} \right\}$ are each of the form 
\[ p_\alpha(u) = (u - z)^\alpha \ , \]
with $\alpha$ a multi-index. 
Then we may write 
\[ A = H \cdot F \] where $H$ is the diagonal matrix with values 
\[ H_{\alpha, \alpha} = h^{\vert \alpha \vert} \ , \ \forall \alpha \in I \] 
and $F$ is the matrix 
\[F = \left( E \left( \chi_{\mathcal{P}} \cdot \left( \frac{\cdot - z }{h} \right)^\alpha \right) (x) \right)_{\alpha, x} \ , \ \forall \alpha \in I \ , \ \forall x \in X \ . \]
From (\ref{eq:R}) and (\ref{eq:L_1}) we get 
\[ \vert F_{\alpha, x} \vert = \left| E \left(\chi_{\mathcal{P}} \cdot \left( \frac{\cdot - z }{h} \right)^\alpha \right)(x) \right| \leq (1+L_1) \cdot R^M \ . \]
Then, 
\begin{eqnarray*}
&\vert \tilde{q}_{x}(z) \vert := 
& \left| \sum\limits_{\alpha \in I} p_{\alpha}\left(z \right) \cdot \left( (ADA^T)^{-1} AD \right)_{\alpha,x} \right| = \\
&& \left| \sum\limits_{\alpha \in I} p_{\alpha}\left(z \right) \cdot \left( ((HF)D (HF)^T)^{-1} (HF) D \right)_{\alpha,x} \right| = \\
&& \left| \sum\limits_{\alpha \in I} p_{\alpha}\left(z \right) \cdot \left( H^{-1} (F D F^T )^{-1} H^{-1} H F D \right)_{\alpha,x} \right| = \\
&& \left| \sum\limits_{\alpha \in I} p_{\alpha}\left(z \right) \cdot \left( H^{-1} (F D F^T )^{-1} F D \right)_{\alpha,x} \right| = \\
&& \left| \sum\limits_{\alpha \in I} \left( \frac{z-z}{h} \right)^\alpha \cdot \left( (F D F^T )^{-1} F D \right)_{\alpha,x} \right| \leq \\
&& \Vert (FDF^T)^{-1} \Vert \cdot (1+L_1) \cdot R^M \cdot \omega \left( \frac{\vert z - x \vert}{h} \right)
\end{eqnarray*}
with 
\[ \Vert (FDF^T)^{-1} \Vert = \sup\left\{ \frac{ \Vert (FDF^T)^{-1} \pmb{v} \Vert}{\Vert \pmb{v} \Vert} : \pmb{v} \neq \pmb{0} \right\} \ , \]
Note that while the term $\Vert (FDF^T)^{-1} \Vert$ is independent upon the value of $h$, it is dependent upon the distribution of the data points.
Thus, 
\[ L_2 \leq \Vert (FDF^T)^{-1} \Vert \cdot (1+L_1) \cdot R^M \cdot K \ , \]
where 
\[ K := \sup \left\{ \sum_{x \in X} \omega \left( \frac{z - x}{h} \right) : z \in \mathbb{R}^n \right\} \ . \]
That is, the operator $\tilde{Q}$ has a bounded Lebesgue constant. Note that this constant is independent of the choice of the polynomial basis $\left\{ p_\alpha \right\}_{\alpha \in I}$.

From (\ref{eq:p^y_1}) we have
\[ p^y_1 (z) = \sum\limits_{x \in X} \tilde{q}_x (z) \cdot Ef(x) \ . \] 
Also, since the operator $\tilde{Q}$ reproduces polynomials we have that 
\[ p^{Taylor}_{y} (z) = \sum\limits_{x \in X} \tilde{q}_x(z) \cdot E (\chi_{\mathcal{P}} \cdot p^{Taylor}_{y}) (x) \ , \] 
where $p^{Taylor}_{y}$ is the Taylor approximation of the function $r$ at the point $y$.

Pick $z \in B_{Rh}(y)$, then 
\begin{eqnarray*}
&& \vert p^y_1 (z) - p^{Taylor}_{y} (z) \vert = \\
&& \left| \sum\limits_{x \in X} \tilde{q}_x (z) \cdot \left( Ef(x) - E (\chi_{\mathcal{P}} \cdot p^{Taylor}_{y}) (x) \right) \right| \leq \\
&& L_2 \cdot \max_{x \in X \cap B_{R_2 h}(z)} \left\{ \vert Ef(x) - E (\chi_{\mathcal{P}} \cdot p^{Taylor}_{y}) (x) \vert \right\}
\end{eqnarray*}
Compute, 
\begin{eqnarray*}
&&\vert Ef(x) - E (\chi_{\mathcal{P}} \cdot p^{Taylor}_{y}) (x) \vert = \\
&&\vert Eg(x) + Er_+(x) - E (\chi_{\mathcal{P}} \cdot p^{Taylor}_{y}) (x) \vert \leq \\
&&\vert Eg(x) \vert + \vert Er_+(x) - E (\chi_{\mathcal{P}} \cdot p^{Taylor}_{y}) (x) \vert = \\
&&\vert Eg(x) \vert + \vert E(\chi_{\mathcal{P}} \cdot r)(x) - E (\chi_{\mathcal{P}} \cdot p^{Taylor}_{y}) (x) \vert = \\
&&\vert Eg(x) \vert + \vert E (\chi_{\mathcal{P}} \cdot (r - p^{Taylor}_{y})) (x) \vert \leq \\
&&C_1 \Vert g \Vert_{C^{M+1}} h^{M+1} + (1+L_1) \Vert r \Vert_{C^{M+1}} \Vert x - y \Vert^{M+1}
\end{eqnarray*}
Hence for all $z \in B_{R h}(y)$ we have 
\begin{eqnarray*}
&&\vert p^y_1(z) - r(z) \vert \leq \\
&&\vert p^y_1(z) - p^{Taylor}_{y}(z) \vert  + \vert p^{Taylor}_{y}(z) - r(z) \vert \leq \\
&& L_2 \left( C_1 \Vert g \Vert_{C^{M+1}} + (1+L_1) \Vert r \Vert_{C^{M+1}} {(R+R_2)}^{M+1} \right) h^{M+1} + \\
&&{R}^{M+1} \Vert r \Vert_{C^{M+1}} h^{M+1}
\end{eqnarray*}
Clearly for 
\begin{eqnarray*}
C_2 = \max{\left\{ L_2 C_1 \quad , \quad L_2 (1+L_1) (R+R_2)^{M+1} + R^{M+1} \right\}}
\end{eqnarray*}
we get for all $z \in B_{R h}(y)$ that
\[ \vert p^y_1(z) - r(z) \vert \leq C_2 \cdot \left( \Vert g \Vert_{C^{M+1}} + \Vert r \Vert_{C^{M+1}} \right) \cdot h^{M+1} \ . \]
Moreover, from (\ref{eq:p^y_1}) and from (\ref{eq:basis}) we get the representation 
\[ p^y_1 = \sum\limits_{\alpha \in I} p_{\alpha} \cdot \underbrace{\sum\limits_{x \in X} \left( (ADA^T)^{-1} AD \right)_{\alpha, x} Ef(x)}_{\lambda_{\alpha}(y)} \ , \] 
from which it is clear that the mappings of the coefficients $y \mapsto \lambda_{\alpha}(y)$ are smooth.  \qed
\end{pf}

\begin{rmk} \label{rem:independence}
In Theorem \ref{thm:p^y_1} we have assumed that the vectors 
\[\left\{ E (\chi_{\mathcal{P}} \cdot p_{\alpha})(x) : x \in \Omega(y) \cap X \right\}_{\alpha \in I}\]
are linearly independent. 
While it seems intuitive that a truncated polynomial can not be approximated by polynomials, we did not succeed in proving this. Thus we leave it as a conjecture.
\end{rmk}

\begin{conj}[Linear independence] \label{conj:independence}
For any point $y \in \mathcal{G}$, denote 
\[ \Omega(y) = \left\{ y + z : \dfrac{\Vert z \Vert}{h} \in supp(\omega) \right\} \ . \]
Then, the vectors 
\[\left\{ E (\chi_{\mathcal{P}} \cdot p_{\alpha})(x) : x \in \Omega(y) \cap X \right\}_{\alpha \in I}\]
are linearly independent.
\end{conj}

We also suggest a second method for approximating $r$ which does not depend upon the above conjecture. 

\subsection{Approximation by partitioned MLS} \label{sub:p^y_2}

This method utilizes two MLS approximations, the first is based only upon the data points in $\mathcal{P}$, while the second is based only upon the data points in $X \setminus \mathcal{P}$. 

\begin{defi}[The approximant $p^y_2$] \label{def:p^y_2} 
For a point $y \in \mathcal{G}$ define $p^y_2 := T_1(f) - T_2(f)$ where 
\[ T_1(f) := \operatorname*{arg\,min}\limits_{p \in \Pi_M(\mathbb{R}^n)} \sum\limits_{x \in \mathcal{P}} \omega \left( \frac{\Vert y-x \Vert }{h}\right) \cdot \left( p(x) - f(x) \right)^2 \ , \] 
and 
\[ T_2(f) := \operatorname*{arg\,min}\limits_{p \in \Pi_M(\mathbb{R}^n)} \sum\limits_{x \in X \setminus \mathcal{P}} \omega \left( \frac{\Vert y-x \Vert }{h}\right) \cdot \left( p(x) - f(x) \right)^2 \ . \] 
Here $\omega$ is a weight function satisfying $supp(\omega) \supset [0, (R+2\Upsilon_M+2)h]$.
\end{defi}
\begin{thm}[$p^y_2$ is a smooth $M$-th order approximation to $r$] \label{thm:p^y_2}
The mapping $y \mapsto p^y_2$ is a smooth $M$-th order approximation to $r$.
\end{thm}

\begin{pf}
Pick $y \in \mathcal{G}$, by Definition \ref{defi:G}, both 
\[ B_{(R + 2 \Upsilon_M + 2)h} (y) \cap \mathcal{P} \quad \mbox{and} \quad B_{(R + 2 \Upsilon_M + 2)h}(y) \cap (X \setminus \mathcal{P}) \]
are uni-solvent for $\Pi_M(\mathbb{R}^n)$. 

The operator $T_1 : \mathbb{R}^X \rightarrow \Pi_M(\mathbb{R}^n)$ is an MLS operator, and as such there exists a constant $c_1$ satisfying 
\[ \vert T_1(\phi)(u) - \phi(u) \vert \leq c_1 \cdot \Vert \phi \Vert_{C^{M+1}} \cdot h^{M+1} \ , \ \forall u \in B_{R h}(y) \ , \ \phi \in C^{M+1}(\mathbb{R}^n) \ . \]
Note that the restriction of $f = g+r_+$ to $\mathcal{P}$ satisfies 
\[ f(x) - g(x) - r(x) = O(h^{M+1}) \ , \ \forall x \in \mathcal{P} \ . \]
Hence 
\[  T_1(f)(u) - T_1(g+r)(u) = O(h^{M+1}) \ , \ \forall u \in B_{R h}(y) \ . \]
Moreover, since $g, r \in C^{M+1}(\mathbb{R}^n)$, 
\[ T_1(f)(u) - g(u) - r(u) = O(h^{M+1}) \ , \ \forall u \in B_{R h} (y) \ . \]
Similarly we may show that 
\[ T_2(f)(u) - g(u) = O(h^{M+1})  \ , \ \forall u \in B_{R h} (y) \ . \]
Therefore, 
\begin{eqnarray*}
\vert p^y_2(u) - r(u) \vert & = & \vert T_1(f)(u) - T_2(f)(u) - r(u) \vert \\
& \leq &  \vert T_1(f)(u) - g(u) - r(u) \vert + \vert g(u) - T_2(f)(u) \vert = O(h^{M+1}) \ . 
\end{eqnarray*}

Last note that since $p^y_2$ is defined as the difference between two MLS solutions, the mappings from $y$ to the coefficients of the polynomial $p^y_2$ must be smooth.
\end{pf}

\section{Computational Complexity}
The correction of the approximation at a point $y \in \mathbb{R}^n$ consists of three steps.
\begin{enumerate}
\item We partition the data points in $X \cap B_{(R + 2 \Upsilon_M + 2)h}(y)$ with respect to the sign of the function $r$ (see Section \ref{sec:signs}). This partitioning is based upon several MLS approximations, hence its complexity is $O(N \cdot (d^2 \cdot N + d^3 \cdot j))$ where 
\[ N := \max \left\{k \in \mathbb{N} : \forall y \in \mathbb{R}^n \ , \ \# B_{(R + 2 \Upsilon_M + 2)h}(y) \cap X \leq k \right\} \ , \]
$j$ is the number of connected components of $B_{(R + 2 \Upsilon_M + 2)h}(y) \setminus \Gamma_r$ and 
\[ d := dim (\Pi_M(\mathbb{R}^n)) \ . \]

\item We find the approximating polynomial to the function $r$, $p^y$. In Section \ref{sec:p^y} we have proposed two methods for finding such a polynomial. Since both methods are essentially based upon MLS, they have the same complexity of $O(d^2  \cdot  N + d^3)$.

\item Last, we set the corrected approximation to be 
\[ \widehat{Q} f(y) = Q f(y) + E ((p^y)_+) (y) \ . \]
This step has a complexity of $O(N)$.
\end{enumerate}

\section{Numerical results} \label{sec:numeric}
\begin{figure*}
\begin{subfigure}{0.5\textwidth}
\includegraphics[width=\textwidth]{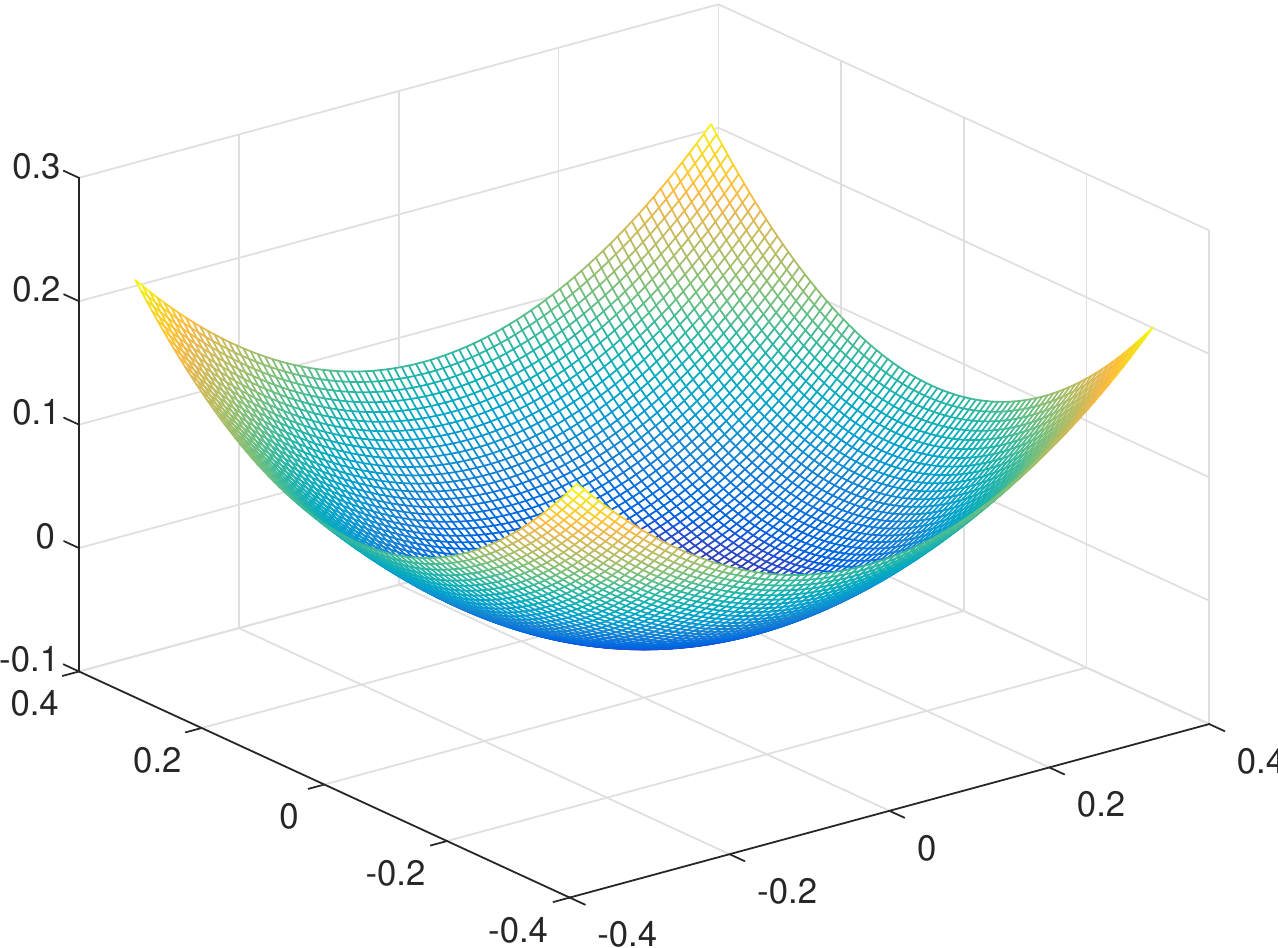}
\caption{$r_1$}
\label{fig:r_1}
\end{subfigure}
\begin{subfigure}{0.5\textwidth}
\includegraphics[width=\textwidth]{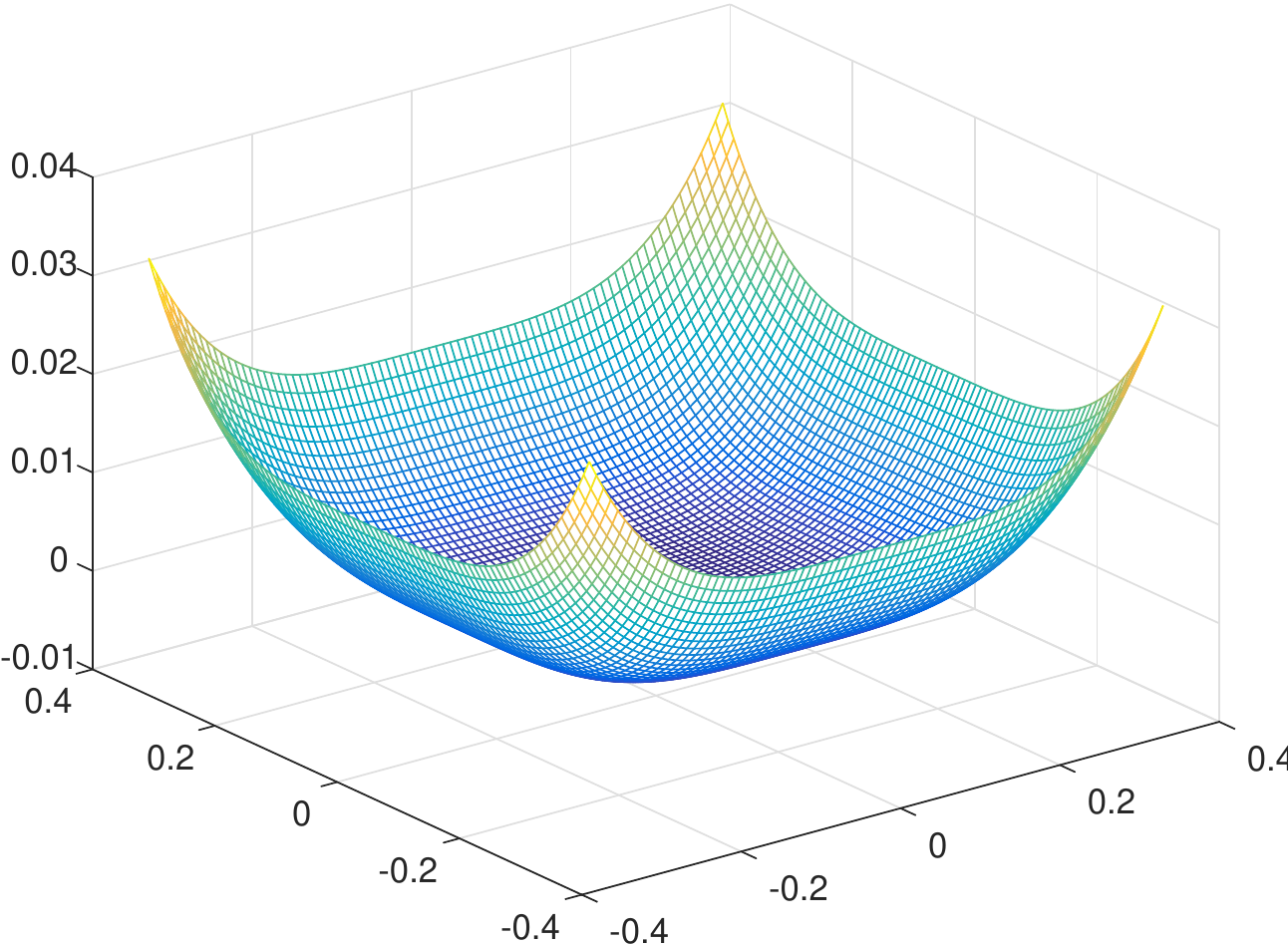}
\caption{$r_2$}
\label{fig:r_2}
\end{subfigure}

\begin{subfigure}{0.5\textwidth}
\includegraphics[width=\textwidth]{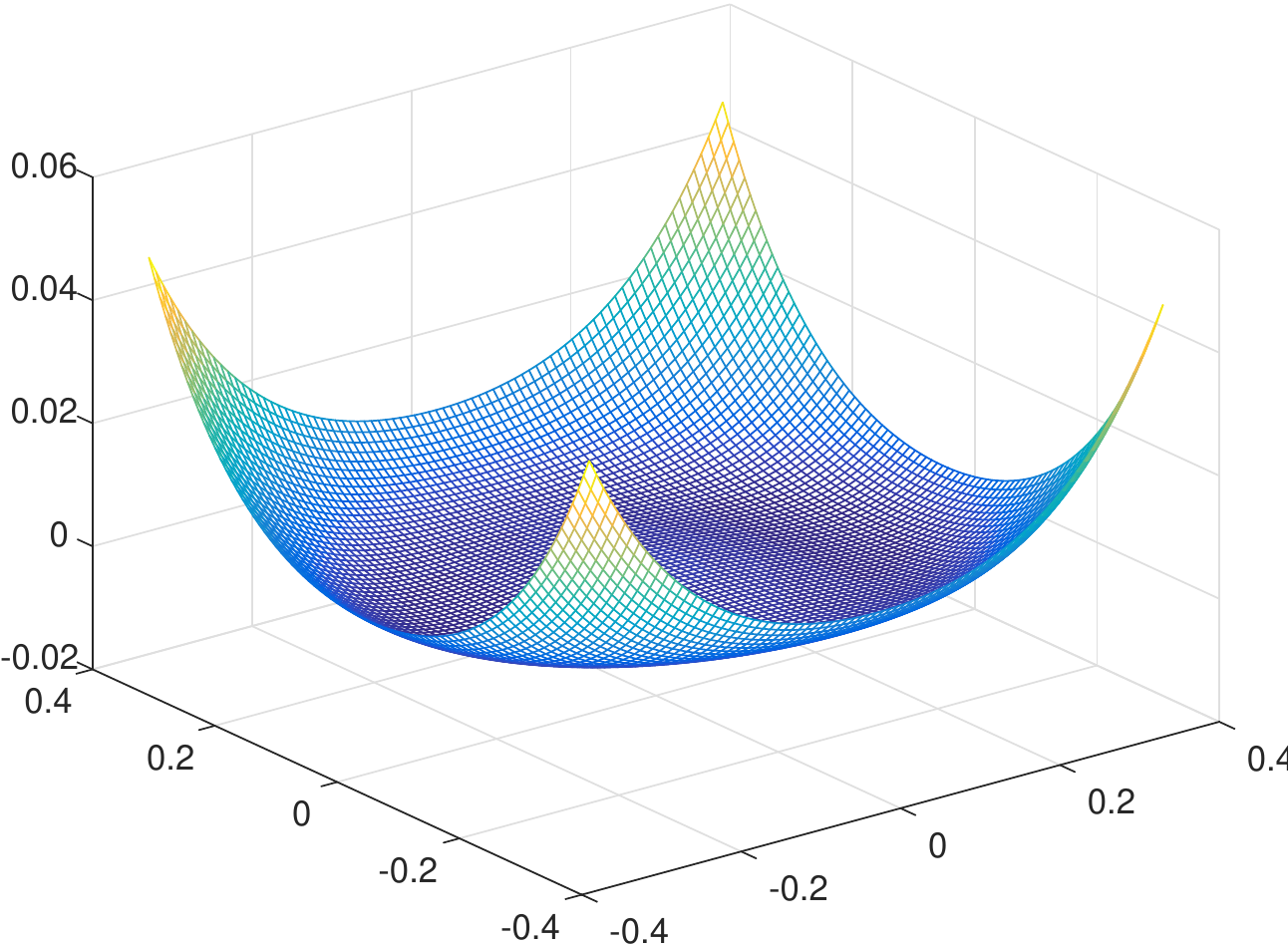}
\caption{$r_3$}
\label{fig:r_3}
\end{subfigure}
\begin{subfigure}{0.5\textwidth}
\includegraphics[width=\textwidth]{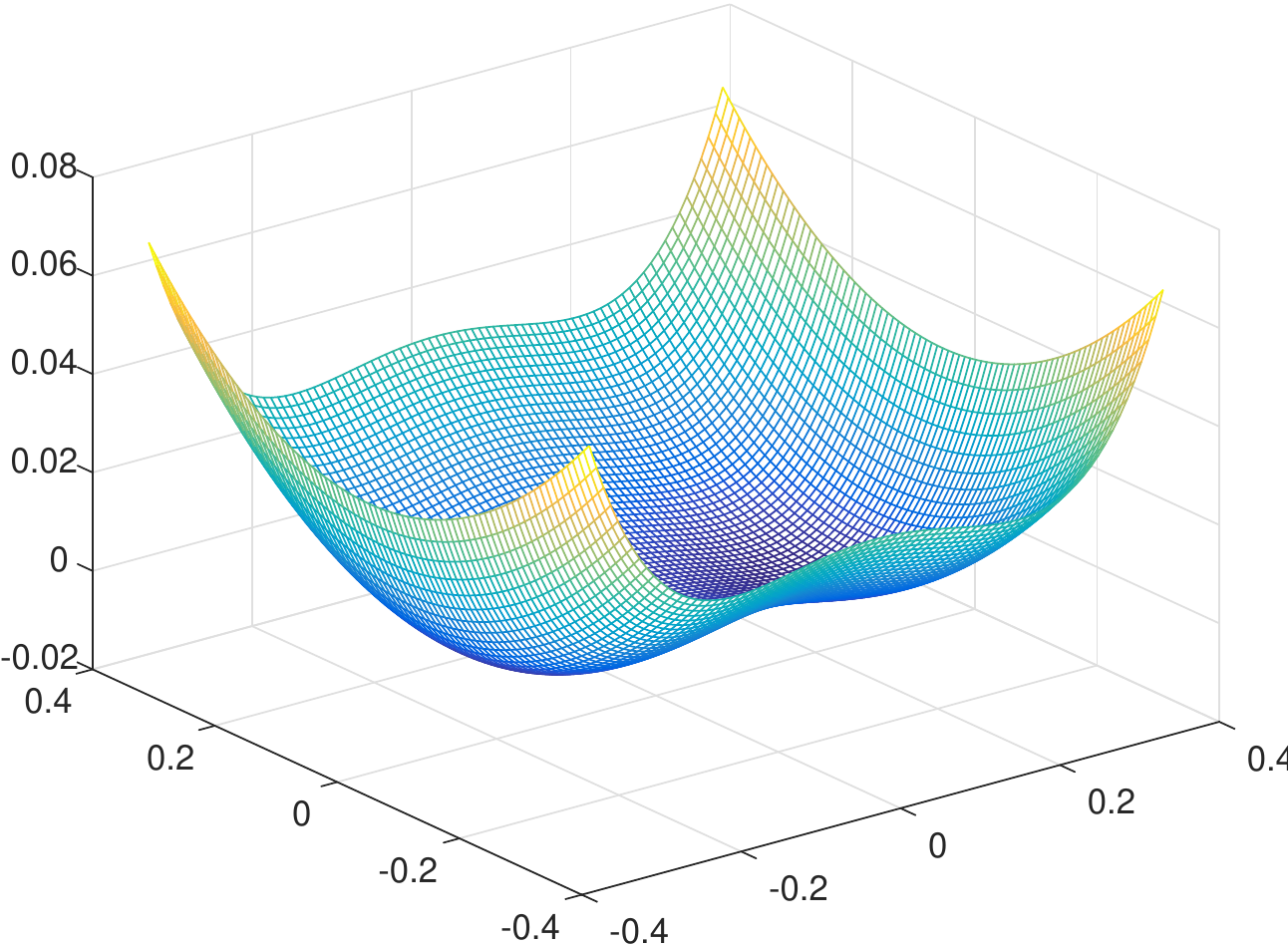}
\caption{$r_4$}
\label{fig:r_4}
\end{subfigure}

\caption{The functions $r_k$ for $k = 1 \ , \ \ldots \ , \ 4$.} 
\label{fig:r}
\end{figure*}

For our test we took our data set $X$ to be a random set of $41^2$ data points with fill distance $0.036$ in the region \[[-0.4 , 0.4]^2 \subset \mathbb{R}^2 \ . \]
Using the sampled data points, we approximated a function, and compared our approximation to the actual function values. We performed this comparison on a uniform mesh of $81^2$ points with fill distance $h = 0.01$ in the region \[[-0.36 , 0.36]^2\subset \mathbb{R}^2 \ . \]
Note that we avoided testing points close to the edge, wishing to avoid approximation errors resulting from partial neighbourhood near the edges.
For the quasi-interpolation we used the MLS approximation with weight function 
\[ \omega \left( \frac{\Vert y - x \Vert}{h} \right) = e^{-\frac{\Vert y - x \Vert^2}{40 h^2}} \ . \]
We approximated the functions $f_k = g + (r_k)_+$ for $1 \leq k \leq 4$ where 
\[ g(x,y) = e^{x+y} \]
and 
\begin{itemize}
\item $r_1(x,y) = x^2 + y^2 - \frac{1}{5}^2$
\item $r_2(x,y) = x^4 + y^4 - \frac{1}{5}^4$
\item $r_3(x,y) = \left(\left(x-\frac{1}{10}\right)^2+y^2-\frac{1}{5}^2\right) \cdot \left(\left(x+\frac{1}{10}\right)^2+y^2-\frac{1}{5}^2\right)$
\item $r_4(x,y) = 4 \cdot x^4+ \frac{y^2 - x^2}{4}$
\end{itemize}
See \cref{fig:r} for the graphs of these functions.

Unless otherwise specified, we have set the value of $M$, the maximal total degree of the approximating polynomials, to $4$.

\begin{figure*}
\begin{subfigure}{0.3\textwidth}
\includegraphics[width=\textwidth]{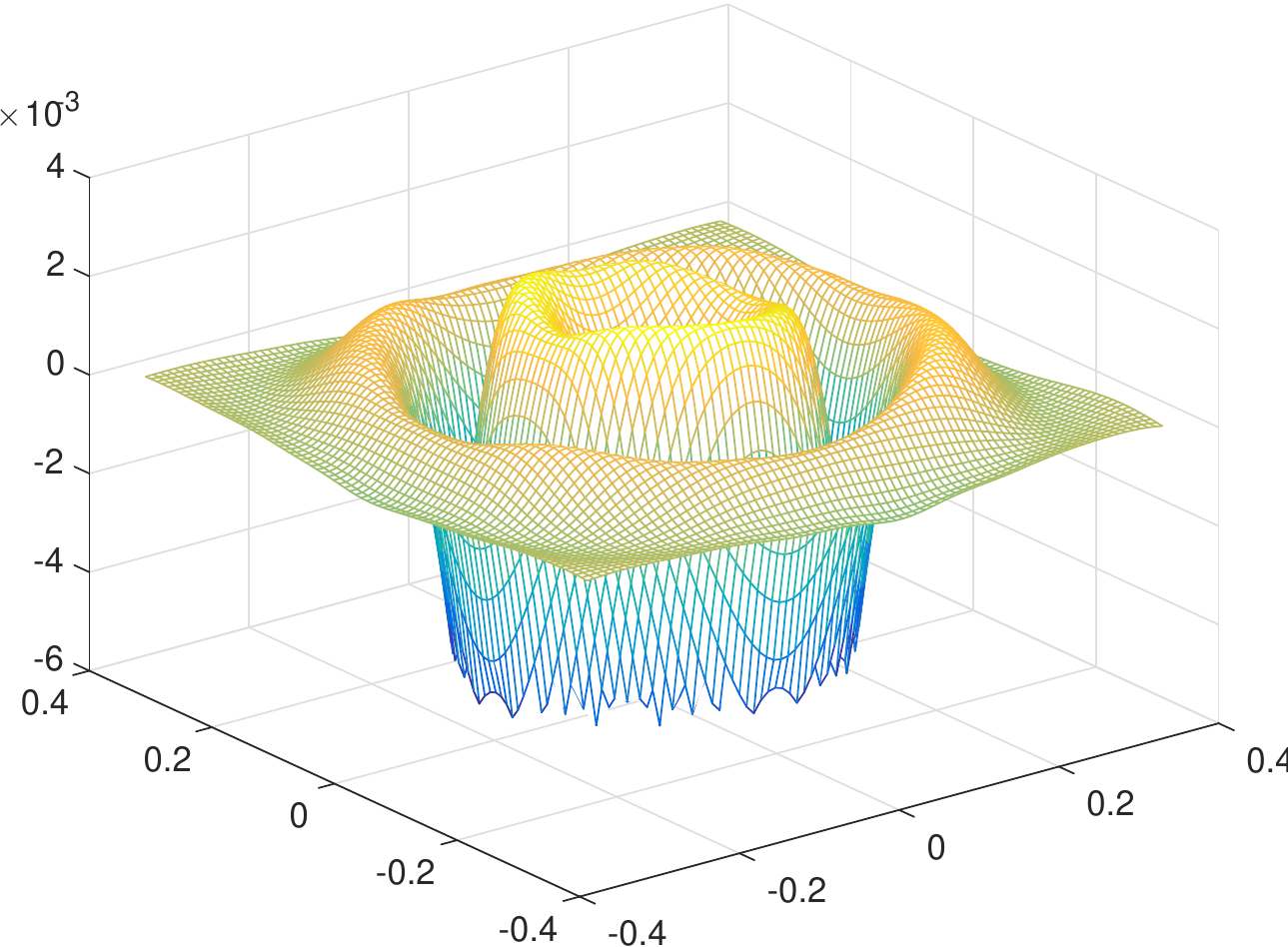}
\caption{$Ef_1$}
\label{fig:Ef_1}
\end{subfigure}
\begin{subfigure}{0.3\textwidth}
\includegraphics[width=\textwidth]{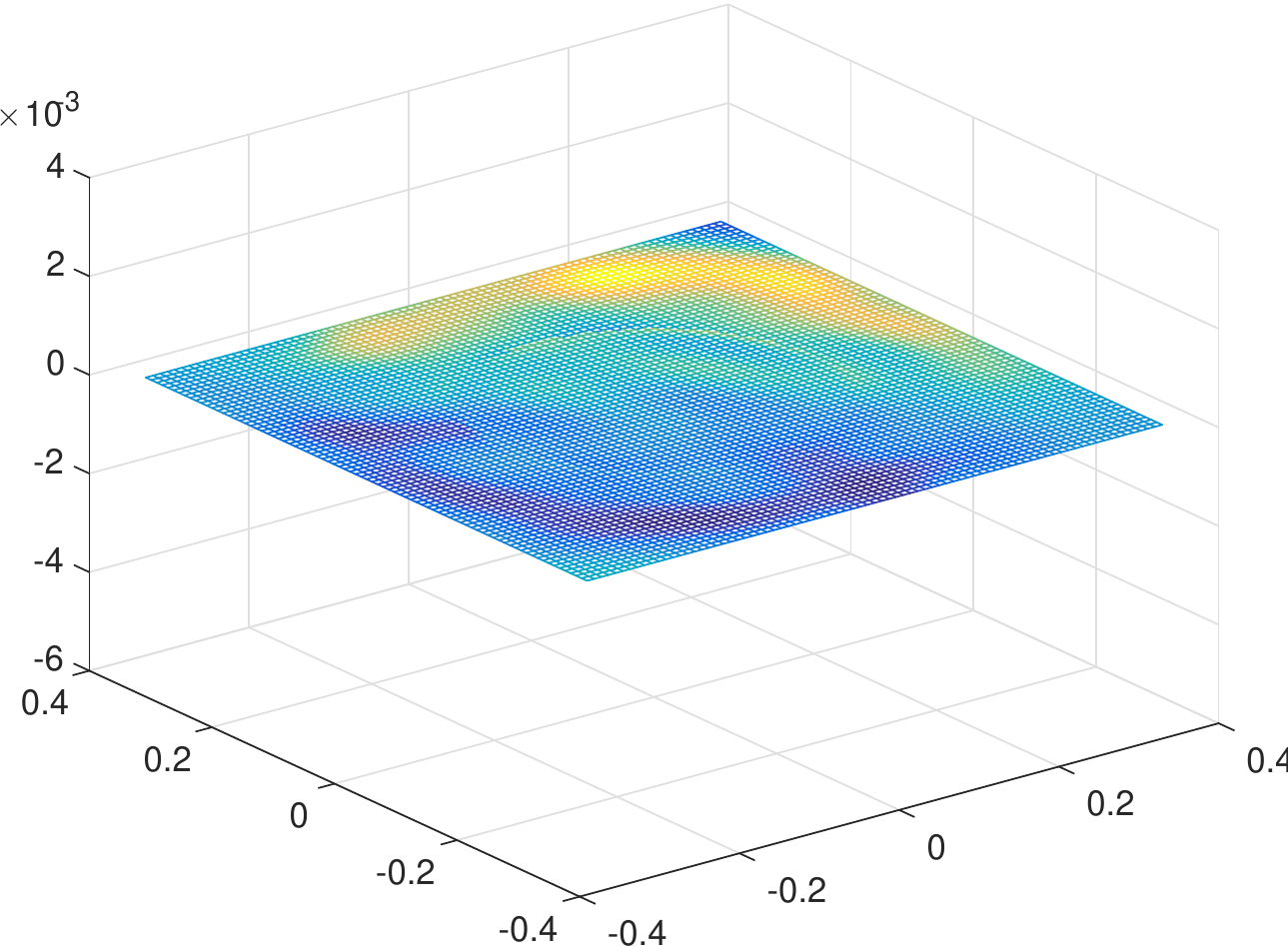}
\caption{$\widehat{E}_1 f_1$}
\label{fig:E_1 f_1}
\end{subfigure}
\begin{subfigure}{0.3\textwidth}
\includegraphics[width=\textwidth]{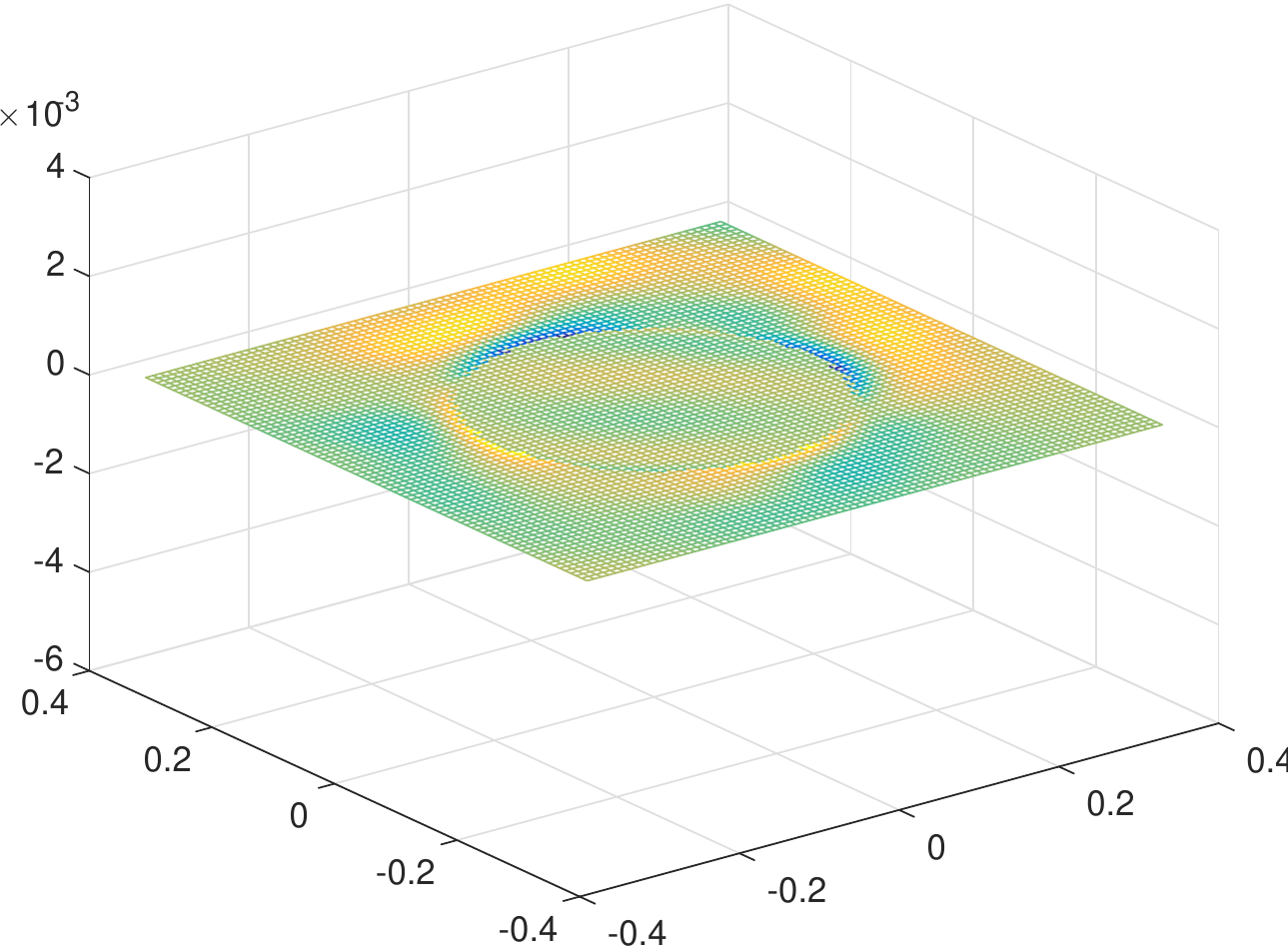}
\caption{$\widehat{E}_2 f_1$}
\label{fig:E_2 f_1}
\end{subfigure}
\caption{Comparison between the original MLS errors and the errors of our corrected approximations for $f_1$. }
\label{fig:errs1}
\end{figure*}

\begin{figure*}
\begin{subfigure}{0.3\textwidth}
\includegraphics[width=\textwidth]{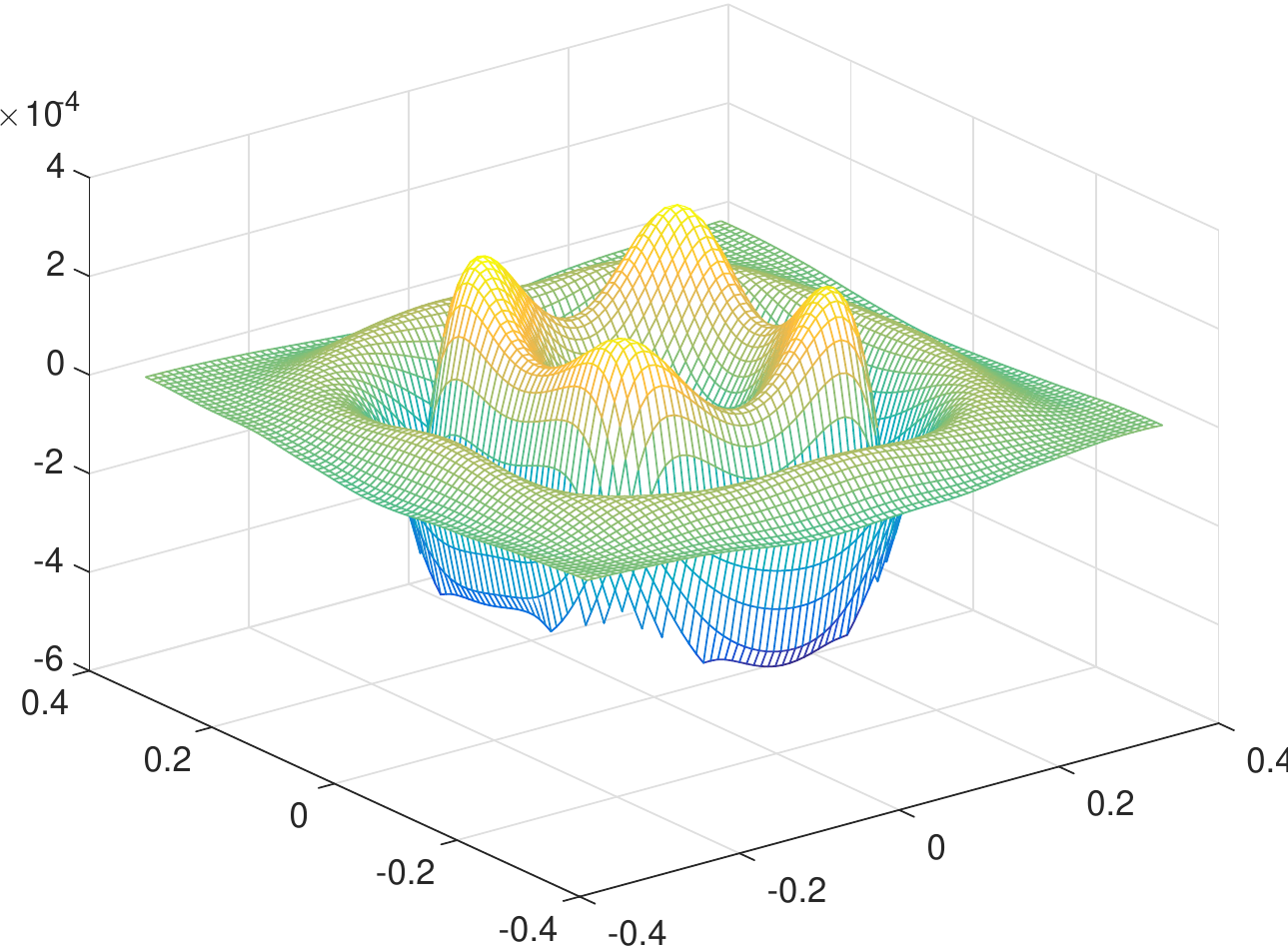}
\caption{$Ef_2$}
\label{fig:Ef_2}
\end{subfigure}
\begin{subfigure}{0.3\textwidth}
\includegraphics[width=\textwidth]{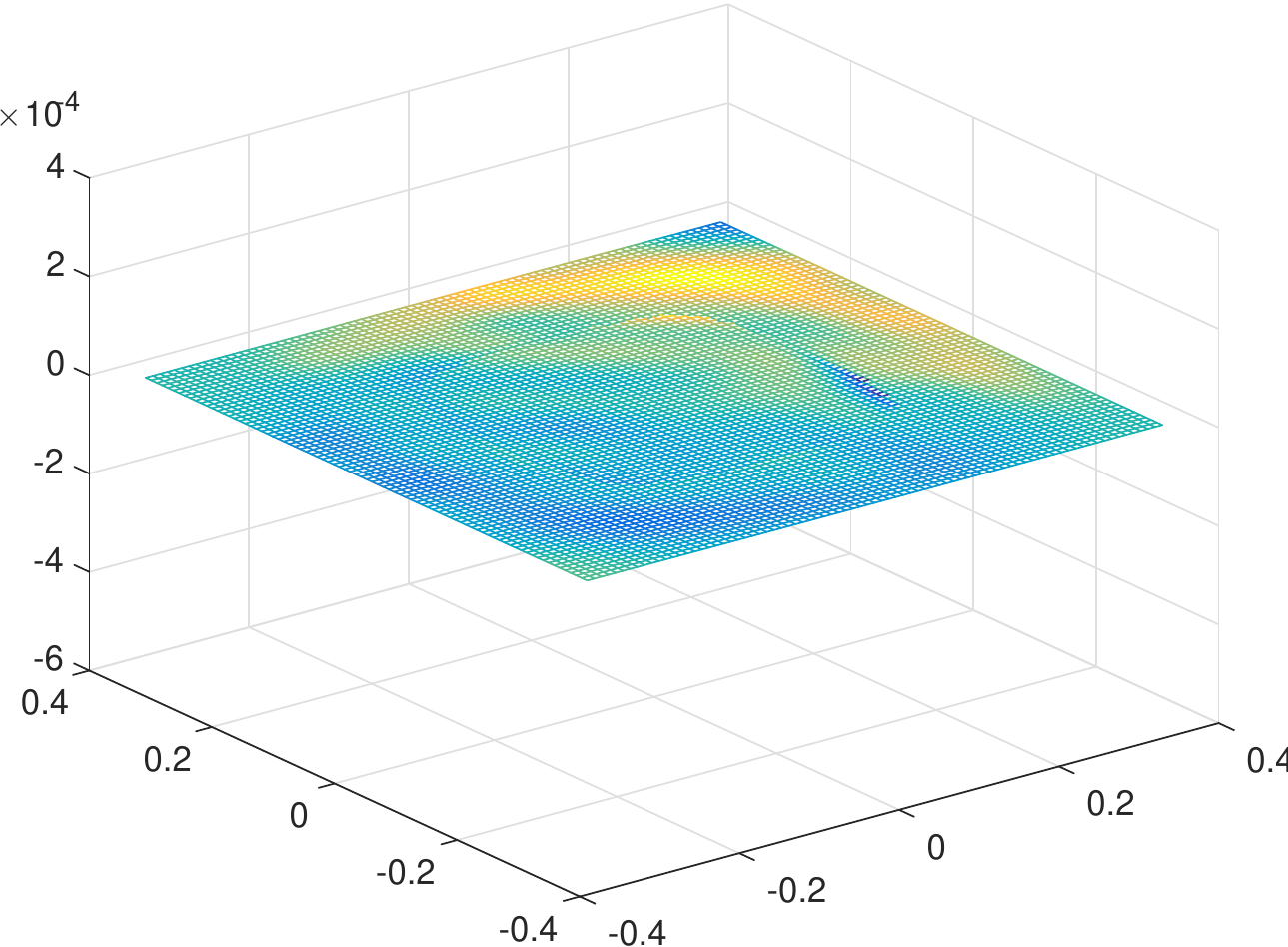}
\caption{$\widehat{E}_1 f_2$}
\label{fig:E_1 f_2}
\end{subfigure}
\begin{subfigure}{0.3\textwidth}
\includegraphics[width=\textwidth]{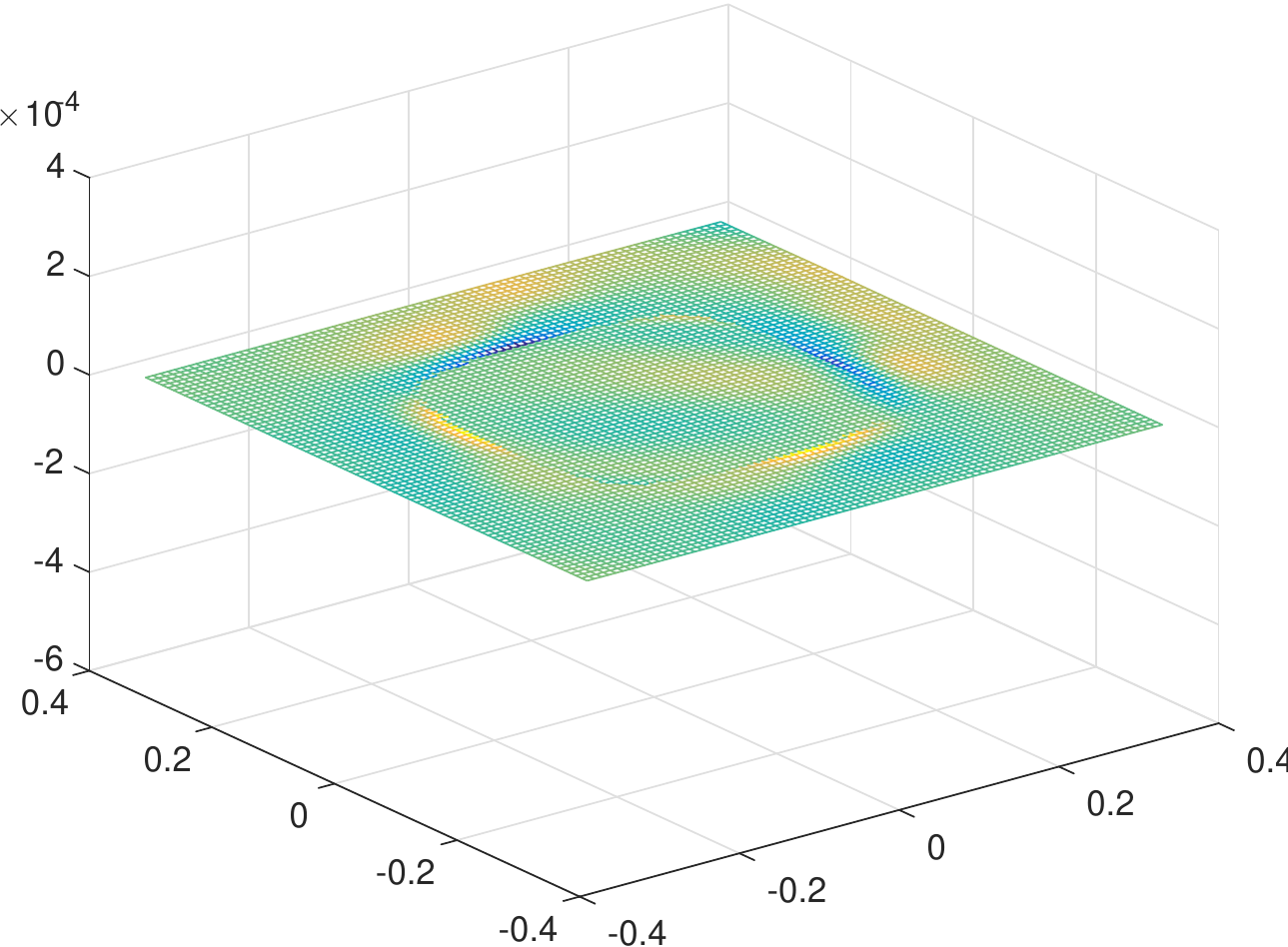}
\caption{$\widehat{E}_2 f_2$}
\label{fig:E_2 f_2}
\end{subfigure}
\caption{Comparison between the original MLS errors and the errors of our corrected approximations for $f_2$. }
\label{fig:errs2}
\end{figure*}

\begin{figure*}
\begin{subfigure}{0.3\textwidth}
\includegraphics[width=\textwidth]{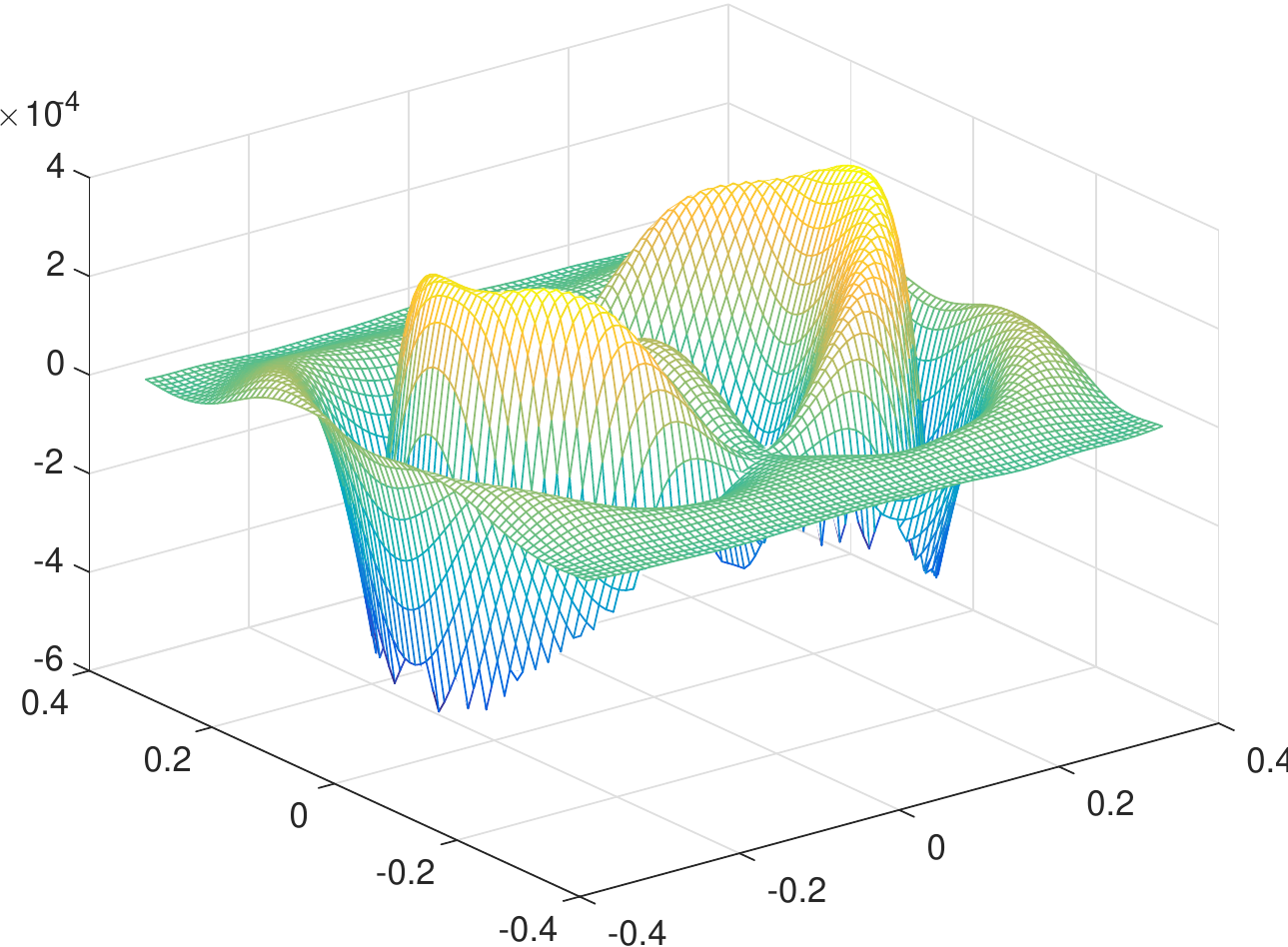}
\caption{$Ef_3$}
\label{fig:Ef_3}
\end{subfigure}
\begin{subfigure}{0.3\textwidth}
\includegraphics[width=\textwidth]{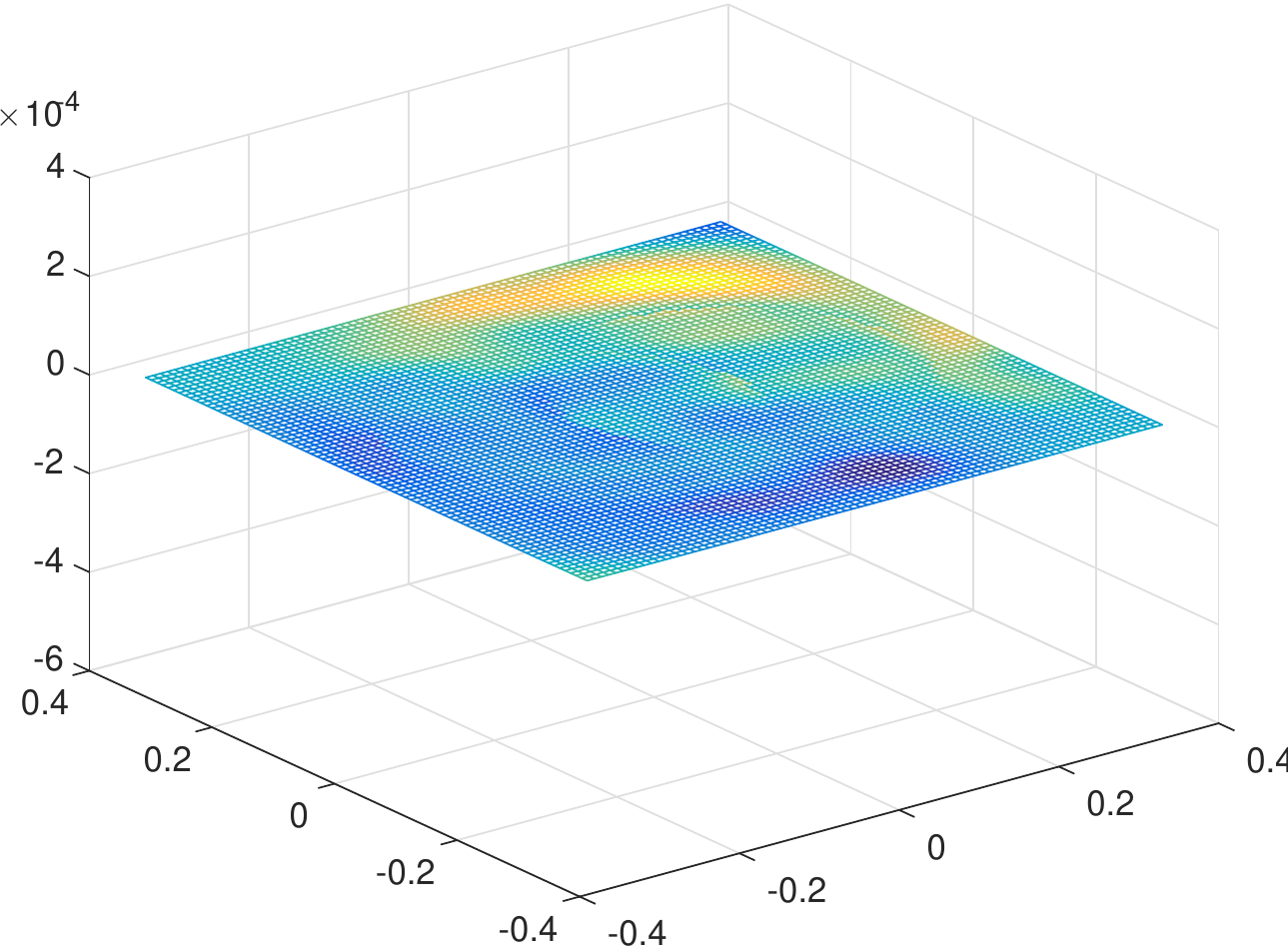}
\caption{$\widehat{E}_1 f_3$}
\label{fig:E_1 f_3}
\end{subfigure}
\begin{subfigure}{0.3\textwidth}
\includegraphics[width=\textwidth]{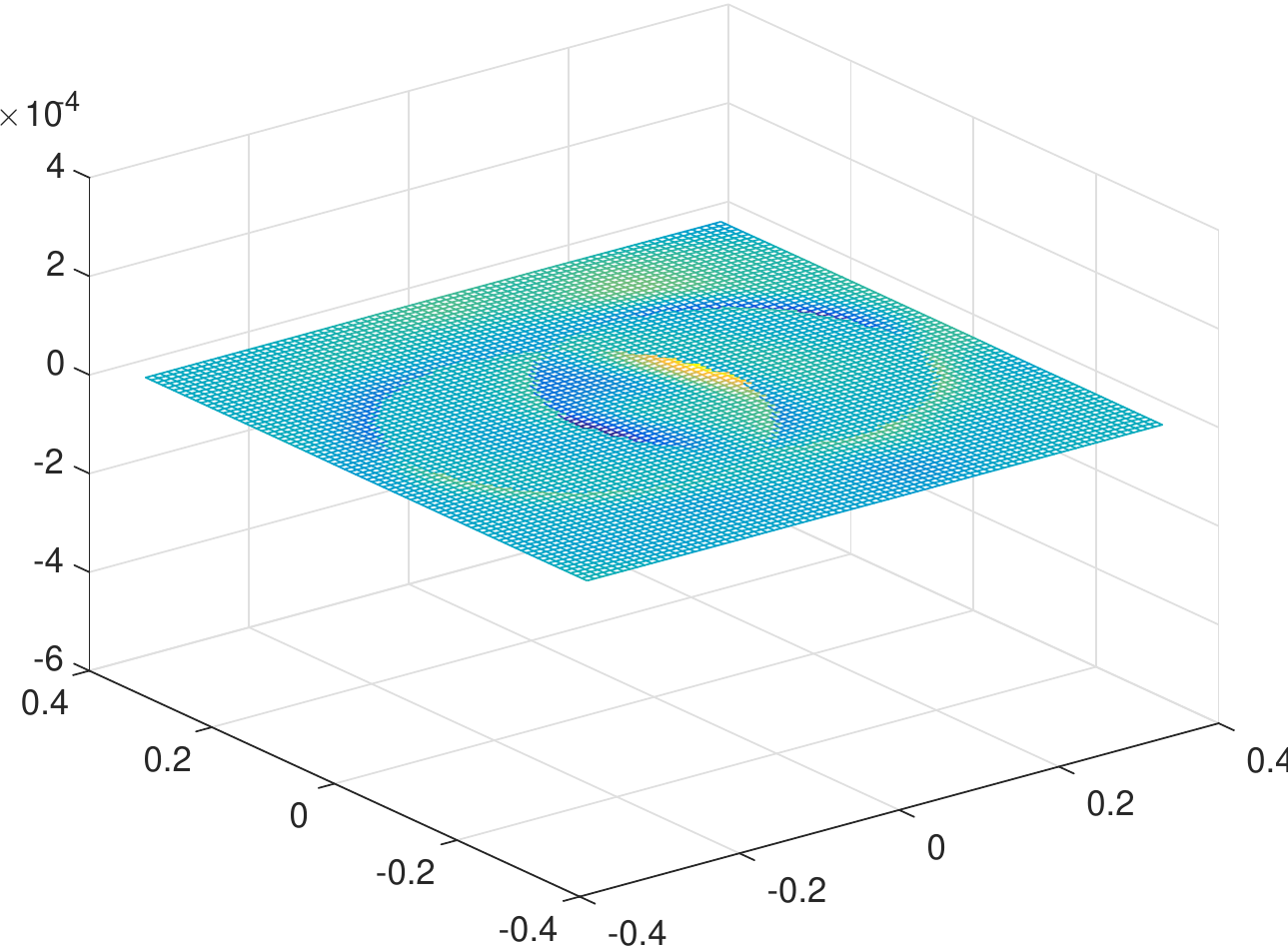}
\caption{$\widehat{E}_2 f_3$}
\label{fig:E_2 f_3}
\end{subfigure}
\caption{Comparison between the original MLS errors and the errors of our corrected approximations for $f_3$. }
\label{fig:errs3}
\end{figure*}

\begin{figure*}
\begin{subfigure}{0.3\textwidth}
\includegraphics[width=\textwidth]{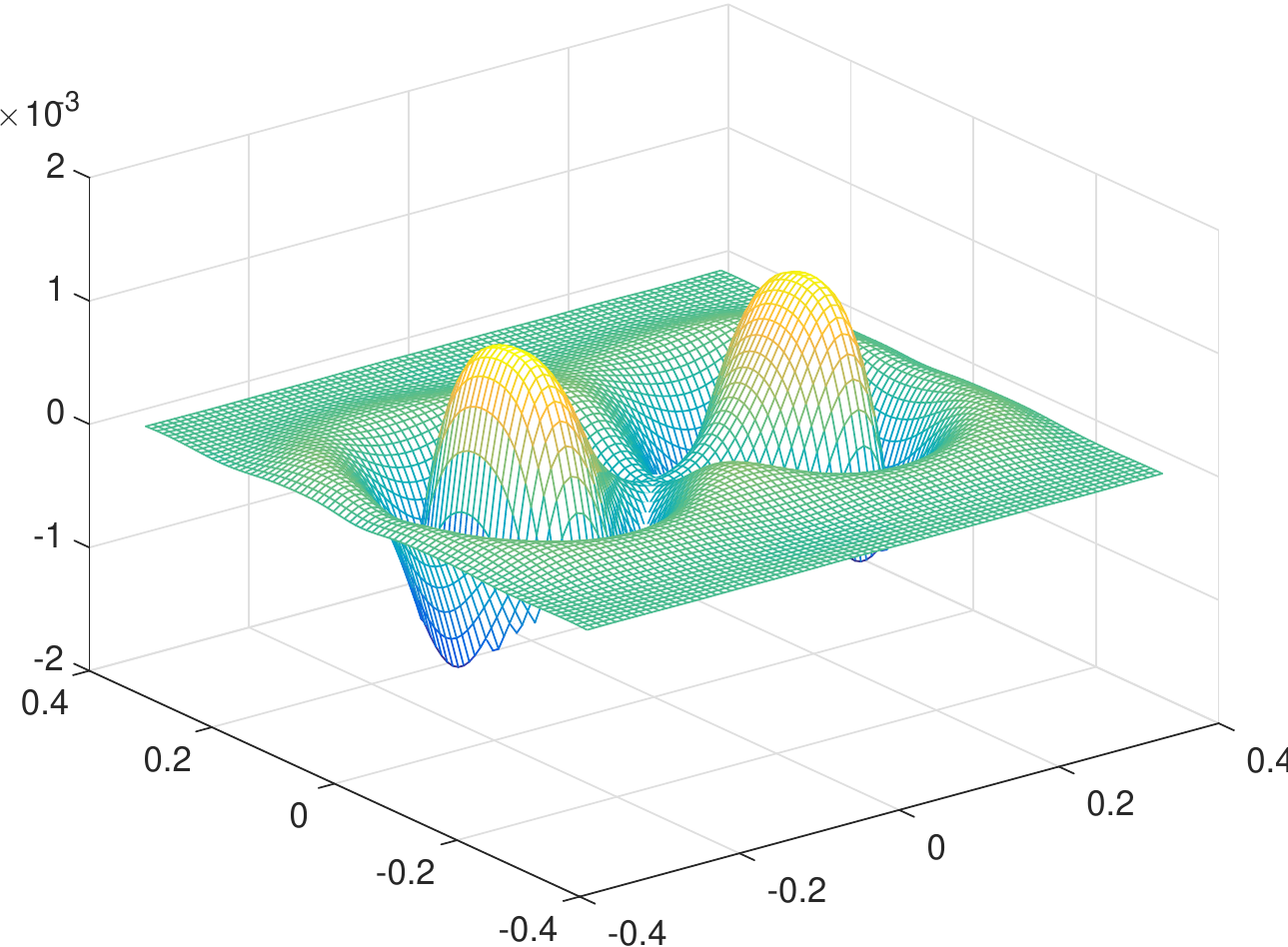}
\caption{$Ef_4$}
\label{fig:Ef_4}
\end{subfigure}
\begin{subfigure}{0.3\textwidth}
\includegraphics[width=\textwidth]{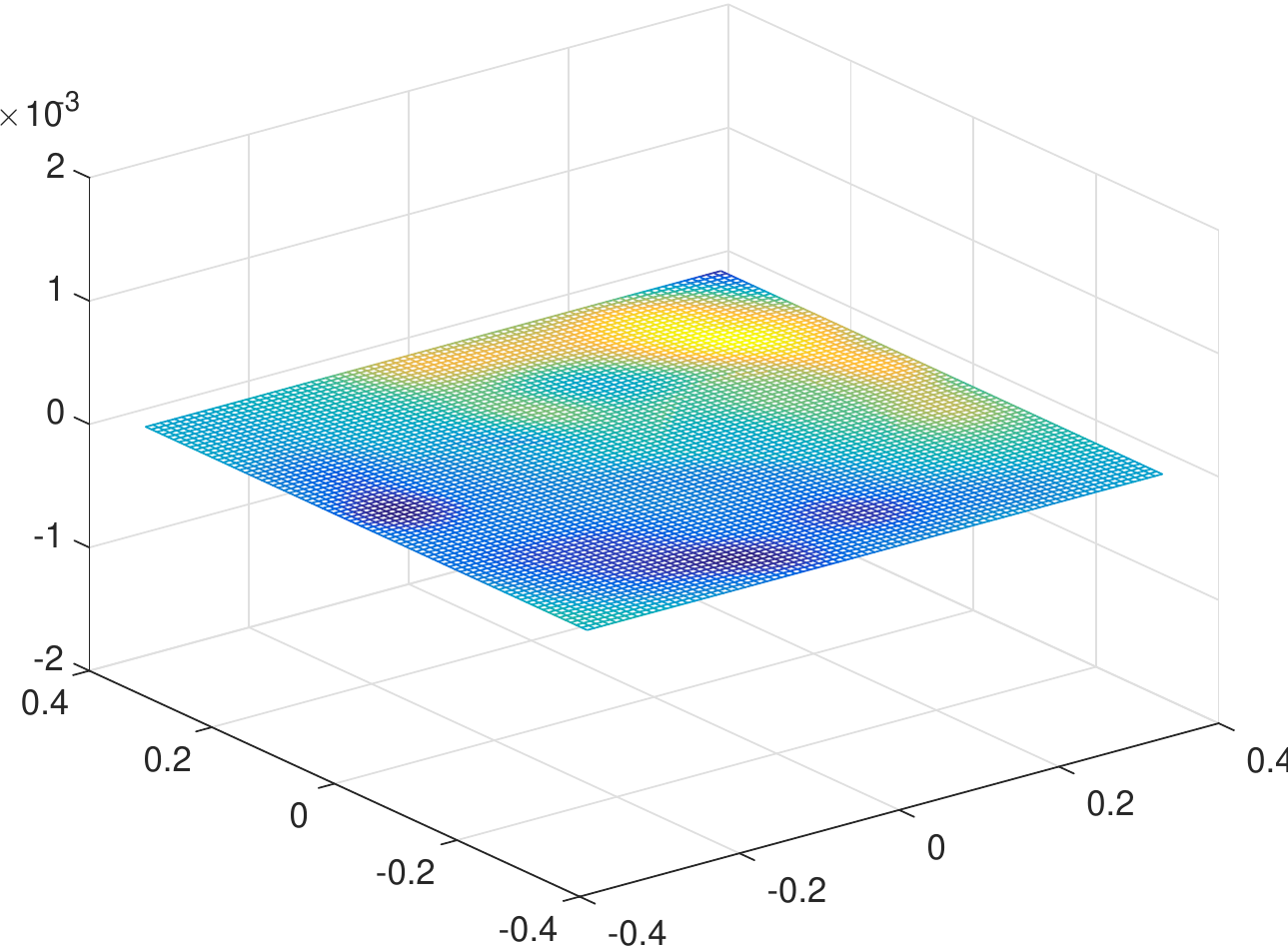}
\caption{$\widehat{E}_1 f_4$}
\label{fig:E_1 f_4}
\end{subfigure}
\begin{subfigure}{0.3\textwidth}
\includegraphics[width=\textwidth]{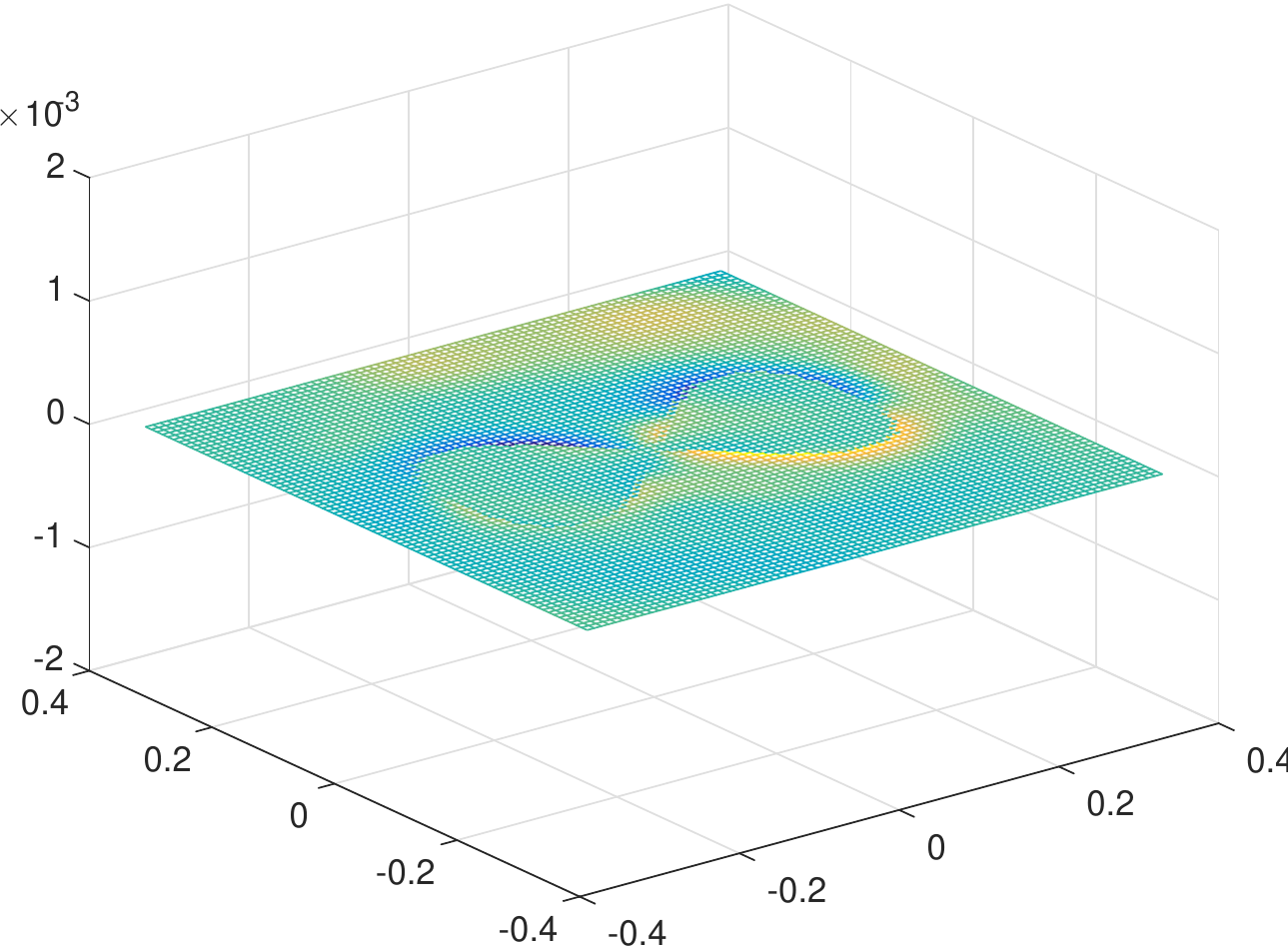}
\caption{$\widehat{E}_2 f_4$}
\label{fig:E_2 f_4}
\end{subfigure}
\caption{Comparison between the original MLS errors and the errors of our corrected approximations for $f_4$. }
\label{fig:errs4}
\end{figure*}

In \cref{fig:errs1,fig:errs2,fig:errs3,fig:errs4}, one can see a comparison between the errors of the original MLS approximation, $Ef$, and our corrected approximations, $\widehat{E}_1 f$ and $\widehat{E}_2 f$, based upon $p^y_1$ and $p^y_2$ respectively.

Note that while the errors of the original approximation were distinctly higher in the vicinity of the singularities, the errors of the corrected approximation based upon $p^y_1$ have the same order near the singularity and far from it. Also, the errors of the corrected approximation based upon $p^y_2$ are much reduced near the singularity.

\begin{figure*}
\begin{subfigure}{0.5\textwidth}
\includegraphics[width=\textwidth]{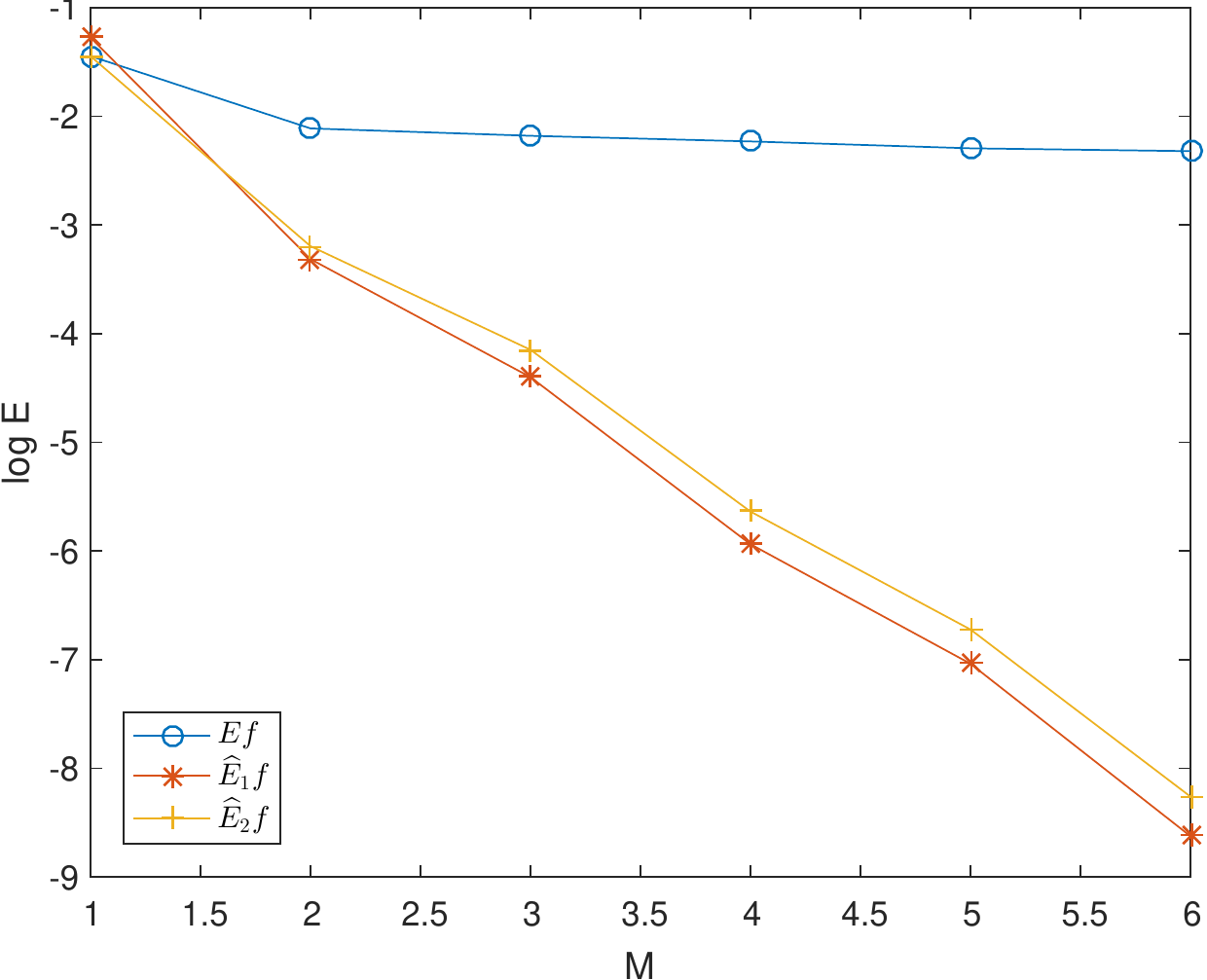}
\caption{$f_1$}
\label{fig:M_test_f1}
\end{subfigure}
\begin{subfigure}{0.5\textwidth}
\includegraphics[width=\textwidth]{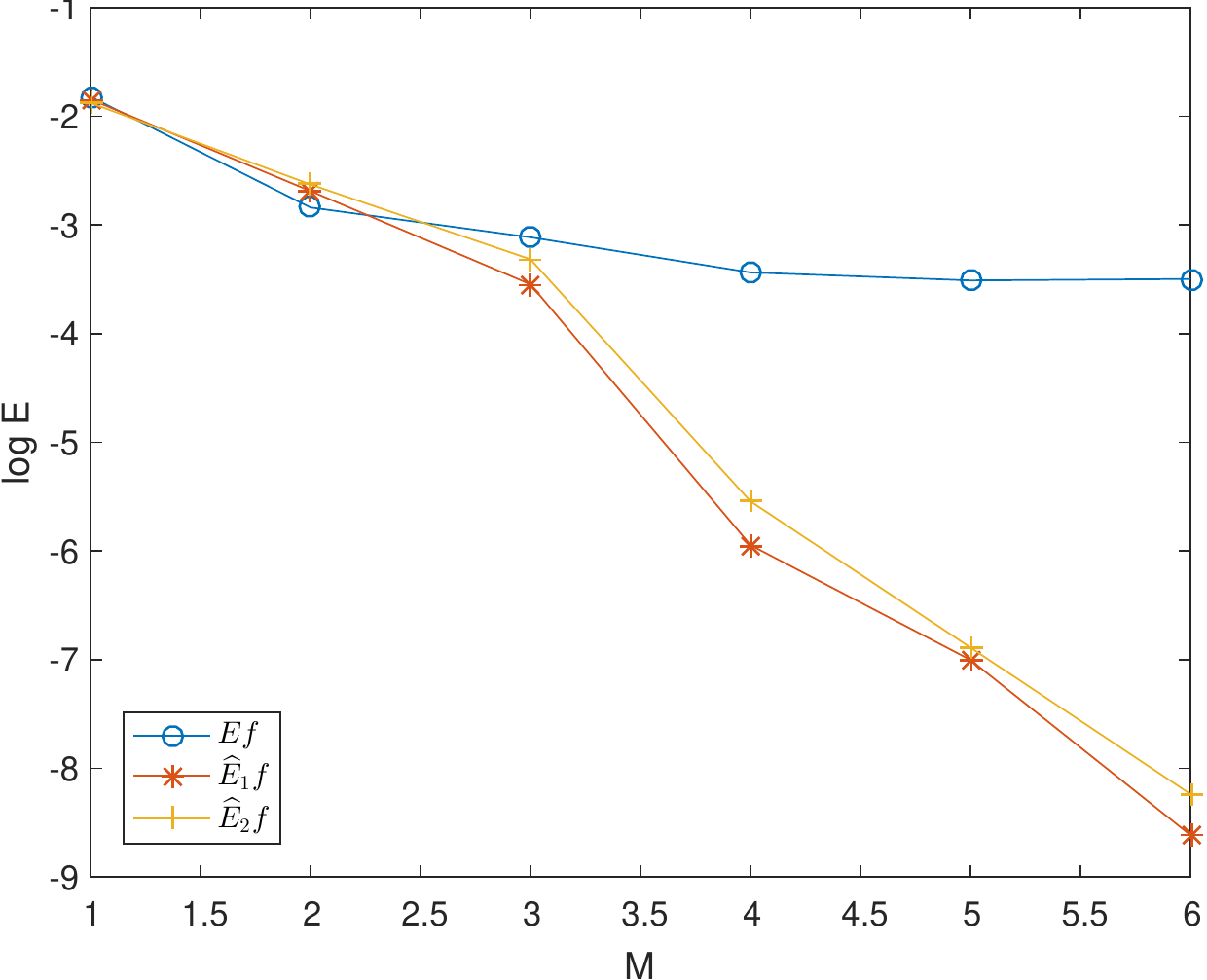}
\caption{$f_2$}
\label{fig:M_test_f2}
\end{subfigure}

\begin{subfigure}{0.5\textwidth}
\includegraphics[width=\textwidth]{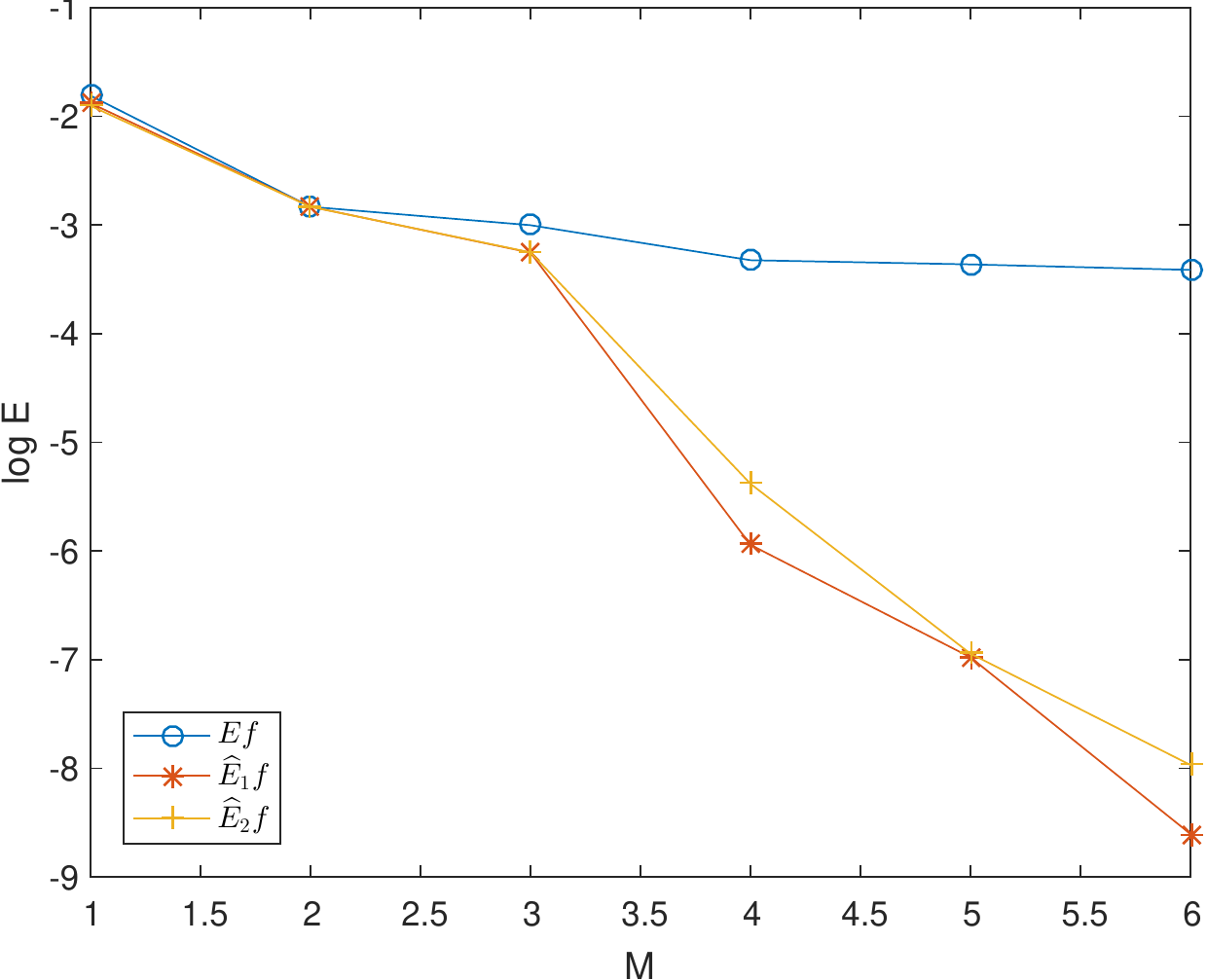}
\caption{$f_3$}
\label{fig:M_test_f3}
\end{subfigure}
\begin{subfigure}{0.5\textwidth}
\includegraphics[width=\textwidth]{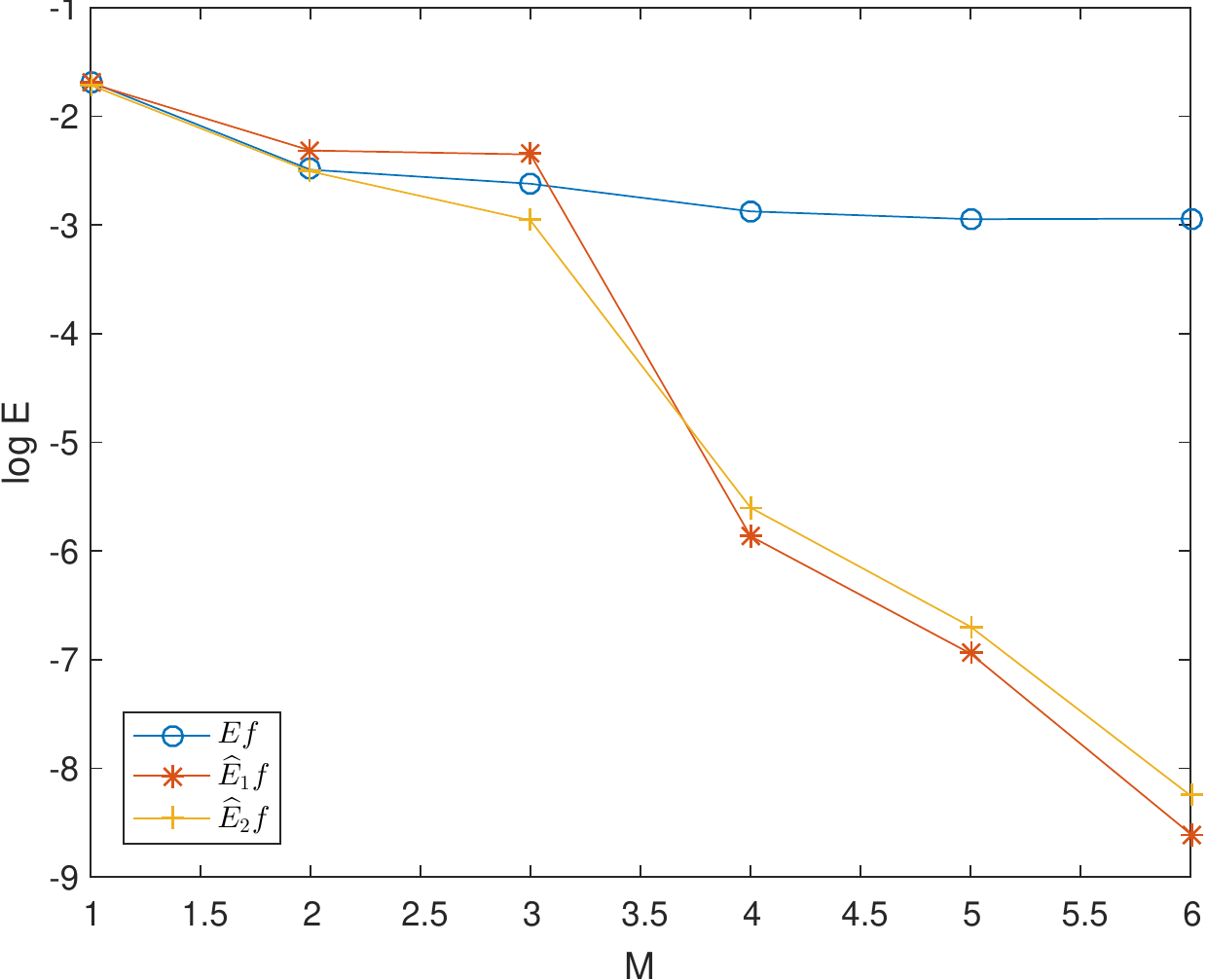}
\caption{$f_4$}
\label{fig:M_test_f4}
\end{subfigure}

\caption{The errors of the original MLS approximation and the corrected approximations for varying $M$ values for $f_1$, $f_2$, $f_3$ and $f_4$.}
\label{fig:M_test}
\end{figure*}

In \cref{tab:changeM} and \cref{fig:M_test}, there is a comparison of the maximal errors on the entire domain for varying $M$. 

\begin{table}[width=\textwidth]
\centering
\caption{Comparison of the maximal errors for varying $M$.}
\label{tab:changeM}{
\begin{tabular}{c|ccc|ccc}
\hline\noalign{\smallskip}
$M$ & \bf $E(f_1)$ &  \bf $\widehat{E}_1(f_1)$ & $\widehat{E}_2(f_1)$ & $E(f_2)$ & $\widehat{E}_1(f_2)$ & $\widehat{E}_2(f_2)$ \\ 
\noalign{\smallskip}\hline\noalign{\smallskip}
$1$ & $h^{0.99}$ &  $h^{0.87}$ &  $h^{0.99}$ & $h^{1.32}$ &  $h^{1.34}$ &  $h^{1.36}$ \\ 
$2$ & $h^{1.45}$ &  $h^{2.27}$ &  $h^{2.19}$ & $h^{2.06}$ &  $h^{1.95}$ &  $h^{1.91}$ \\ 
$3$ & $h^{1.49}$ &  $h^{3.02}$ &  $h^{2.84}$ & $h^{2.26}$ &  $h^{2.58}$ &  $h^{2.41}$ \\ 
$4$ & $h^{1.53}$ &  $h^{4.07}$ &  $h^{3.86}$ & $h^{2.50}$ &  $h^{4.32}$ &  $h^{4.03}$ \\ 
$5$ & $h^{1.57}$ &  $h^{4.82}$ &  $h^{4.61}$ & $h^{2.55}$ &  $h^{5.09}$ &  $h^{5.01}$ \\ 
$6$ & $h^{1.59}$ &  $h^{5.91}$ &  $h^{5.66}$ & $h^{2.54}$ &  $h^{6.26}$ &  $h^{5.99}$ \\ 
\noalign{\smallskip}\hline
\end{tabular}

\begin{tabular}{c|ccc|ccc}
\hline\noalign{\smallskip}
$M$ & \bf $E(f_3)$ & $\widehat{E}_1(f_3)$ & $\widehat{E}_2(f_3)$ & $E(f_4)$ & $\widehat{E}_1(f_4)$ & $\widehat{E}_2 (f_4)$ \\
\noalign{\smallskip}\hline\noalign{\smallskip}
$1$ & $h^{1.28}$ &  $h^{1.33}$ &  $h^{1.35}$ & $h^{1.18}$ &  $h^{1.19}$ &  $h^{1.20}$ \\ 
$2$ & $h^{2.00}$ &  $h^{2.00}$ &  $h^{2.00}$ & $h^{1.74}$ &  $h^{1.62}$ &  $h^{1.75}$ \\ 
$3$ & $h^{2.12}$ &  $h^{2.30}$ &  $h^{2.30}$ & $h^{1.83}$ &  $h^{1.64}$ &  $h^{2.07}$ \\ 
$4$ & $h^{2.35}$ &  $h^{4.21}$ &  $h^{3.81}$ & $h^{2.01}$ &  $h^{4.10}$ &  $h^{3.92}$ \\ 
$5$ & $h^{2.38}$ &  $h^{4.94}$ &  $h^{4.92}$ & $h^{2.06}$ &  $h^{4.86}$ &  $h^{4.68}$ \\ 
$6$ & $h^{2.42}$ &  $h^{6.09}$ &  $h^{5.64}$ & $h^{2.06}$ &  $h^{6.02}$ &  $h^{5.77}$ \\ 
\noalign{\smallskip}\hline
\end{tabular}}
\end{table}

Here too we can see that the errors of the original approximation attain a maximum that is independent of $M$, while the maximal errors of the corrected approximations decrease as $M$ increases.

\begin{figure}

\begin{subfigure}{0.5\textwidth}
\includegraphics[width=\textwidth]{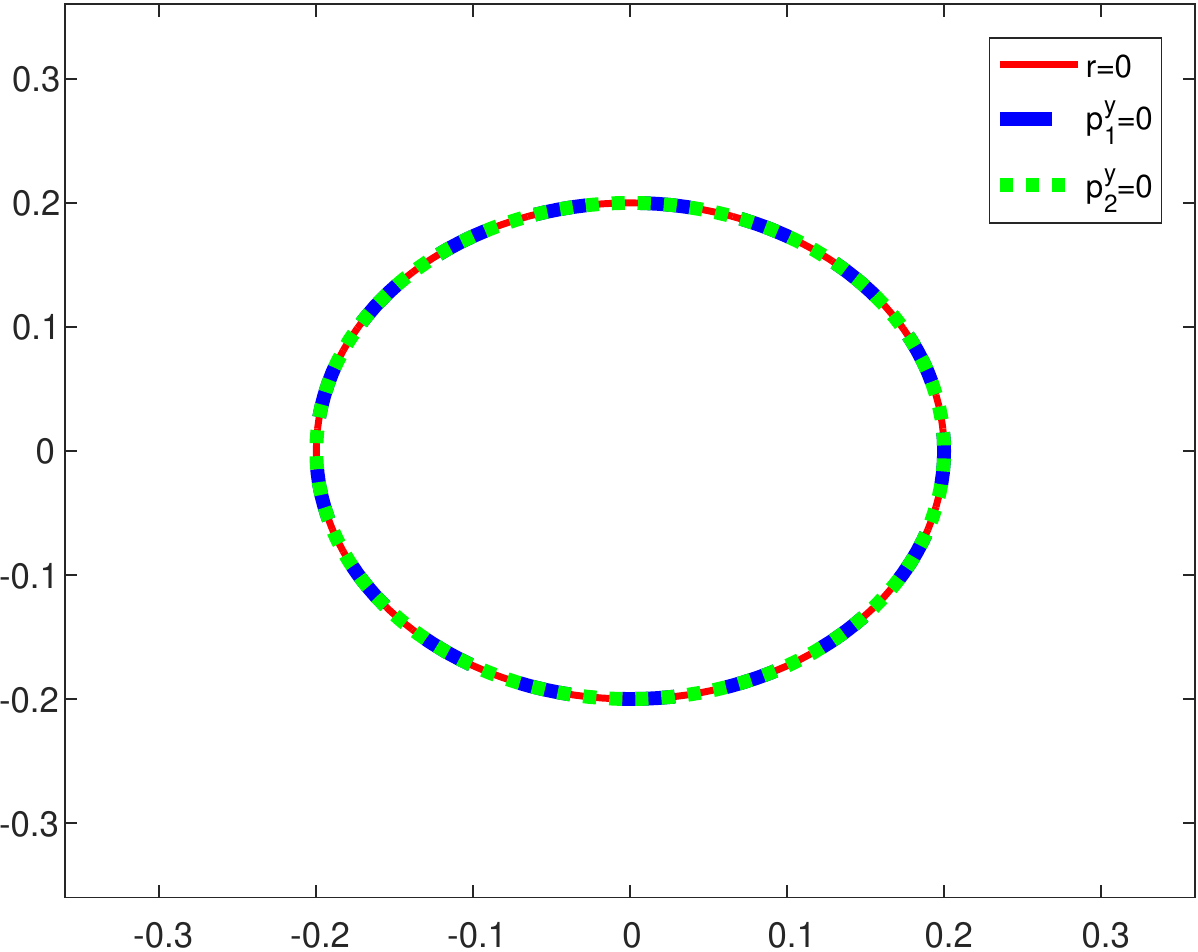}
\caption{$\left\{ r_1 = 0 \right\}$}
\label{fig:curve_test_f1}
\end{subfigure}
\begin{subfigure}{0.5\textwidth}
\includegraphics[width=\textwidth]{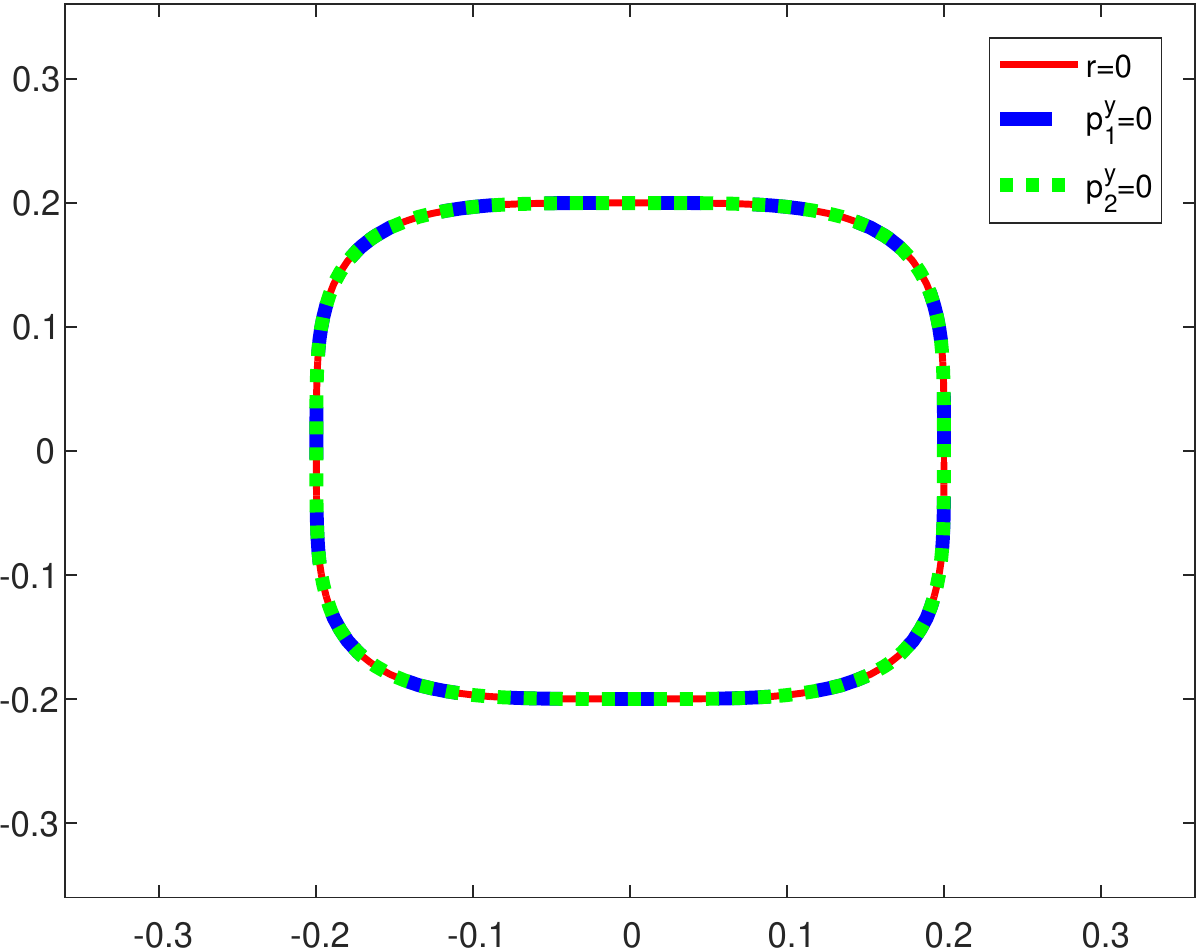}
\caption{$\left\{ r_2 = 0 \right\}$}
\label{fig:curve_test_f2}
\end{subfigure}

\begin{subfigure}{0.5\textwidth}
\includegraphics[width=\textwidth]{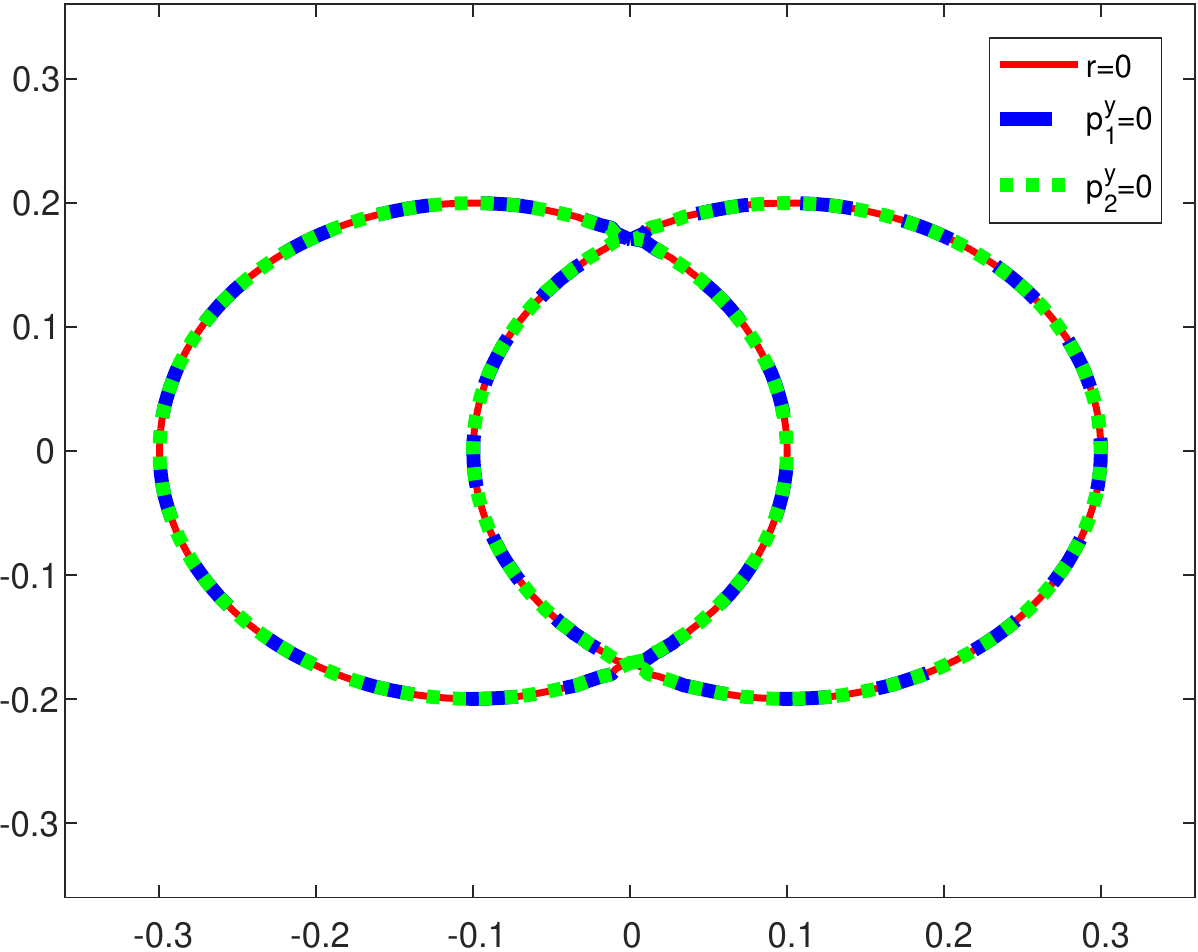}
\caption{$\left\{ r_3 = 0 \right\}$}
\label{fig:curve_test_f3}
\end{subfigure}
\begin{subfigure}{0.5\textwidth}
\includegraphics[width=\textwidth]{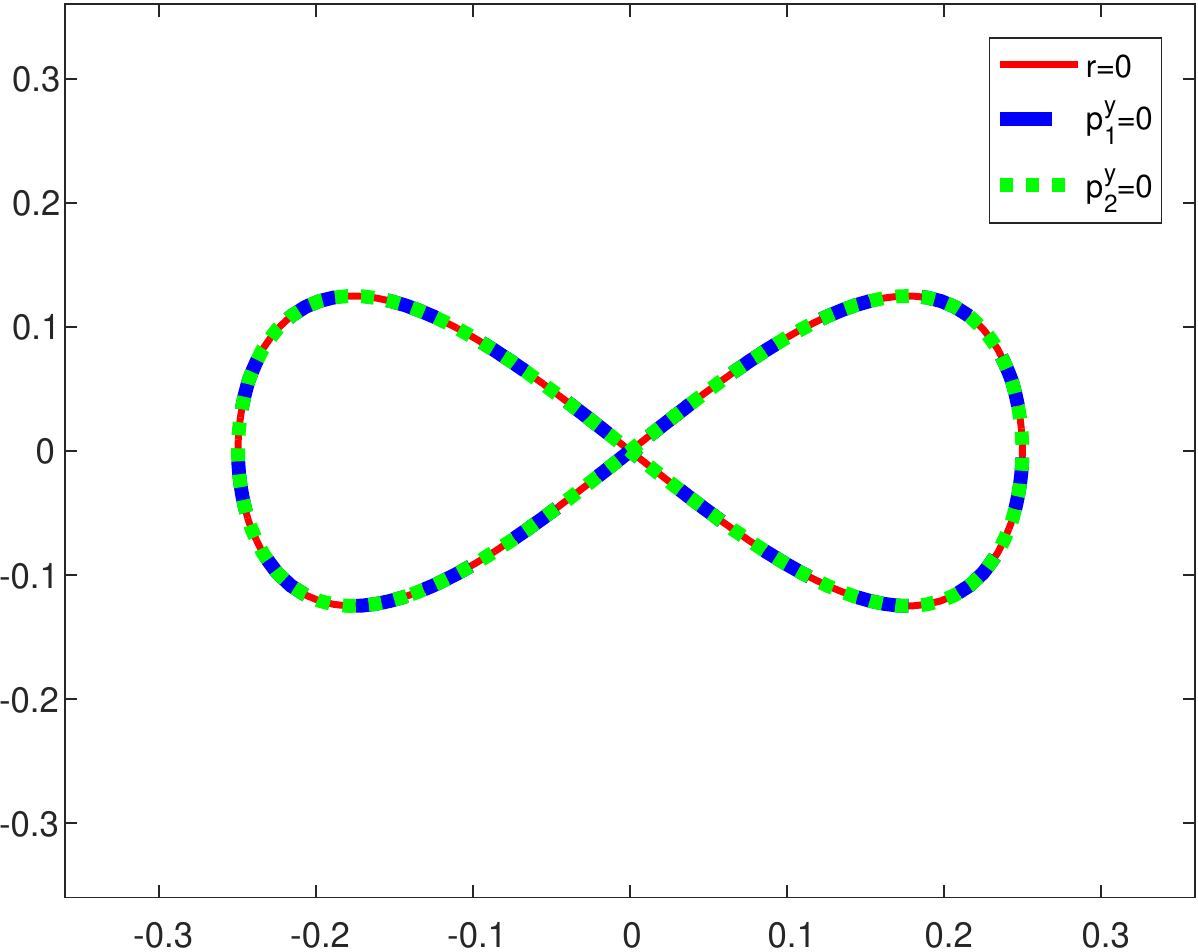}
\caption{$\left\{ r_4 = 0 \right\}$}
\label{fig:curve_test_f4}
\end{subfigure}

\caption{Comparison between the curves $\left\{ r_k=0 \right\} $ (red) and their approximations $\left\{ p^y_1 = 0 \right\}$ (blue dashed) and $\left\{ p^y_2 = 0 \right\}$ (green dotted) for $k = 1 \ , \ \ldots \ , \ 4$. }
\label{fig:curves}
\end{figure}

In \cref{fig:curves}, we drew the curves $\left\{ r_k = 0 \right\}$ on which the function $f_k$ has singularities. We also drew the curves $\left\{ p^y_1(y) = 0 \right\}$ and $\left\{ p^y_2(y) = 0 \right\}$, for each $r_k$, these are our approximations of the curve $\left\{ r_k = 0 \right\}$. We have used the MATLAB \texttt{contour} command to draw these curves.

In these figures we see that the curves $\left\{ r_k = 0 \right\}$ are indistinguishable from their approximations based upon $p^y_1$ and $p^y_2$.

\begin{figure}
\centering
\includegraphics[width=0.5\textwidth]{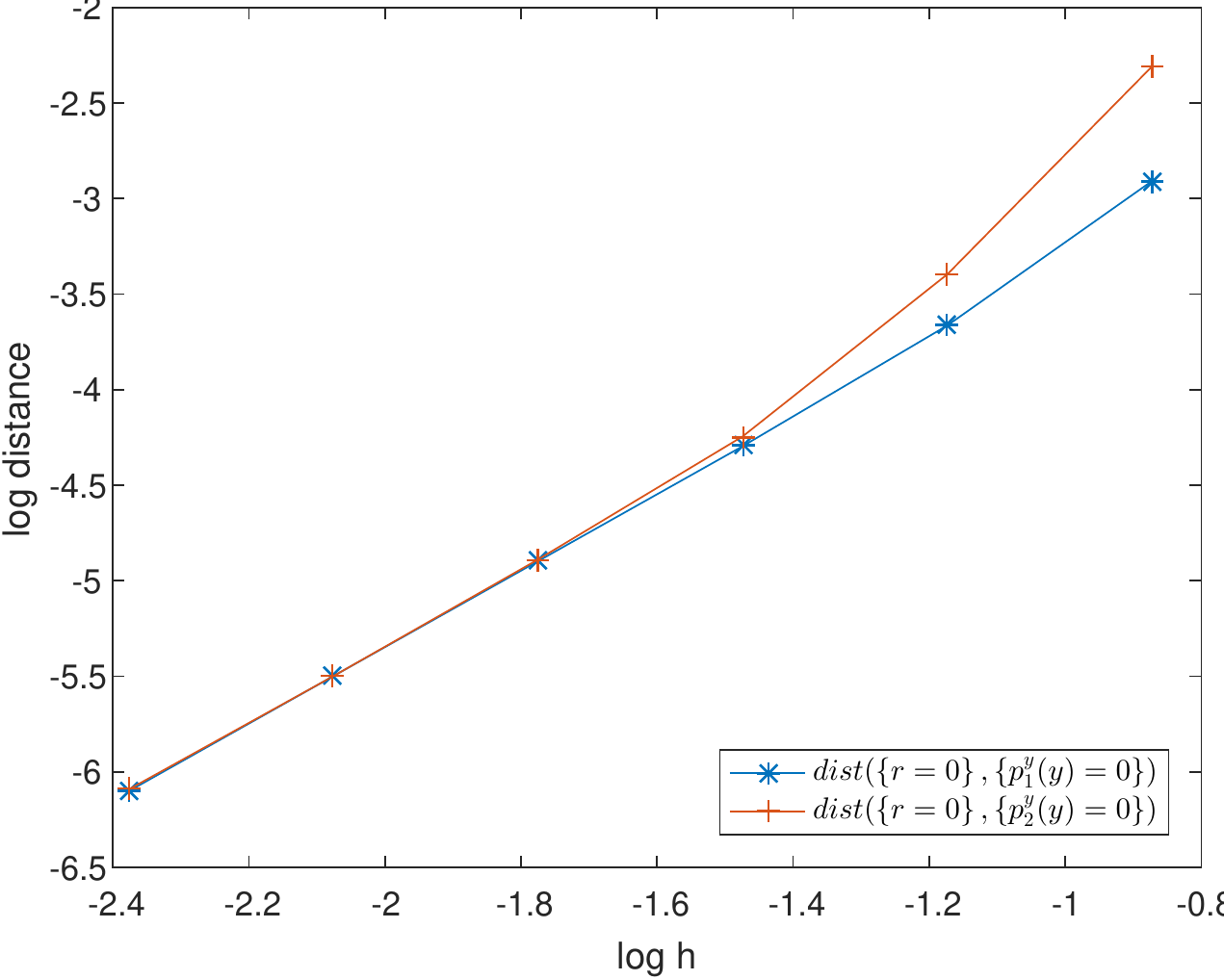}
\caption{The Hausdorff distance between the curve $\left\{ r_1=0 \right\}$ and our approximations of this curve $\left\{ p^y_1=0 \right\}$ and $\left\{ p^y_2 = 0 \right\}$ for varying values of $h$.}
\label{fig:curve_test}
\end{figure}

In \cref{fig:curve_test}, we outline the Hausdorff distance between the singularity curve $\left\{ r_1=0 \right\}$ and our approximations of this curve, $\left\{ p^y_1 = 0 \right\}$ and $\left\{ p^y_2 = 0 \right\}$ for varying values of fill-distance $h$. Note that as we took varying $h$ values we did not change the number of data points, but only the size of the region in which we test the procedure.

\begin{figure*}
\begin{subfigure}{0.5\textwidth}
\includegraphics[width=\textwidth]{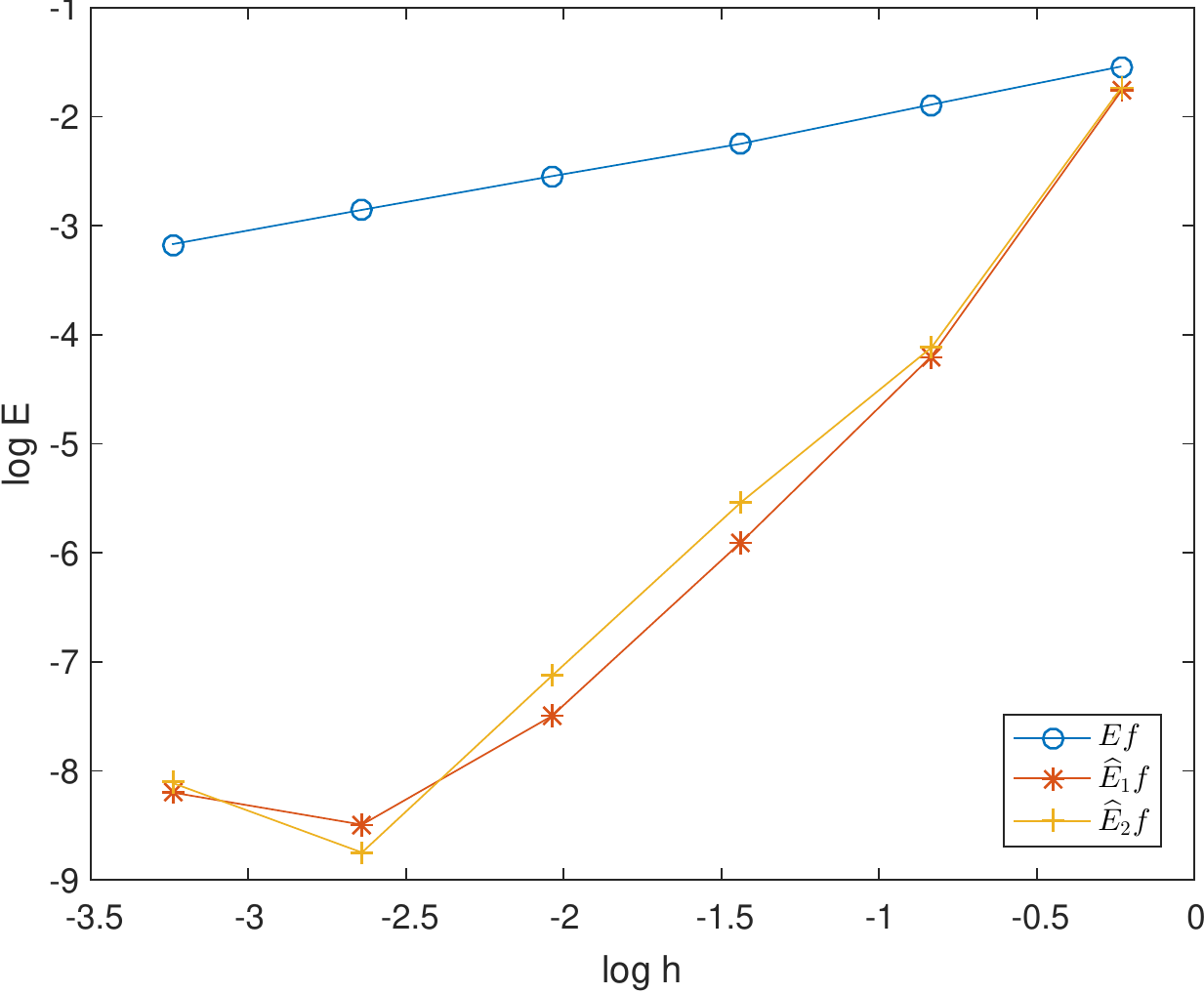}
\caption{$f_1$}
\label{fig:h_test_f1}
\end{subfigure}
\begin{subfigure}{0.5\textwidth}
\includegraphics[width=\textwidth]{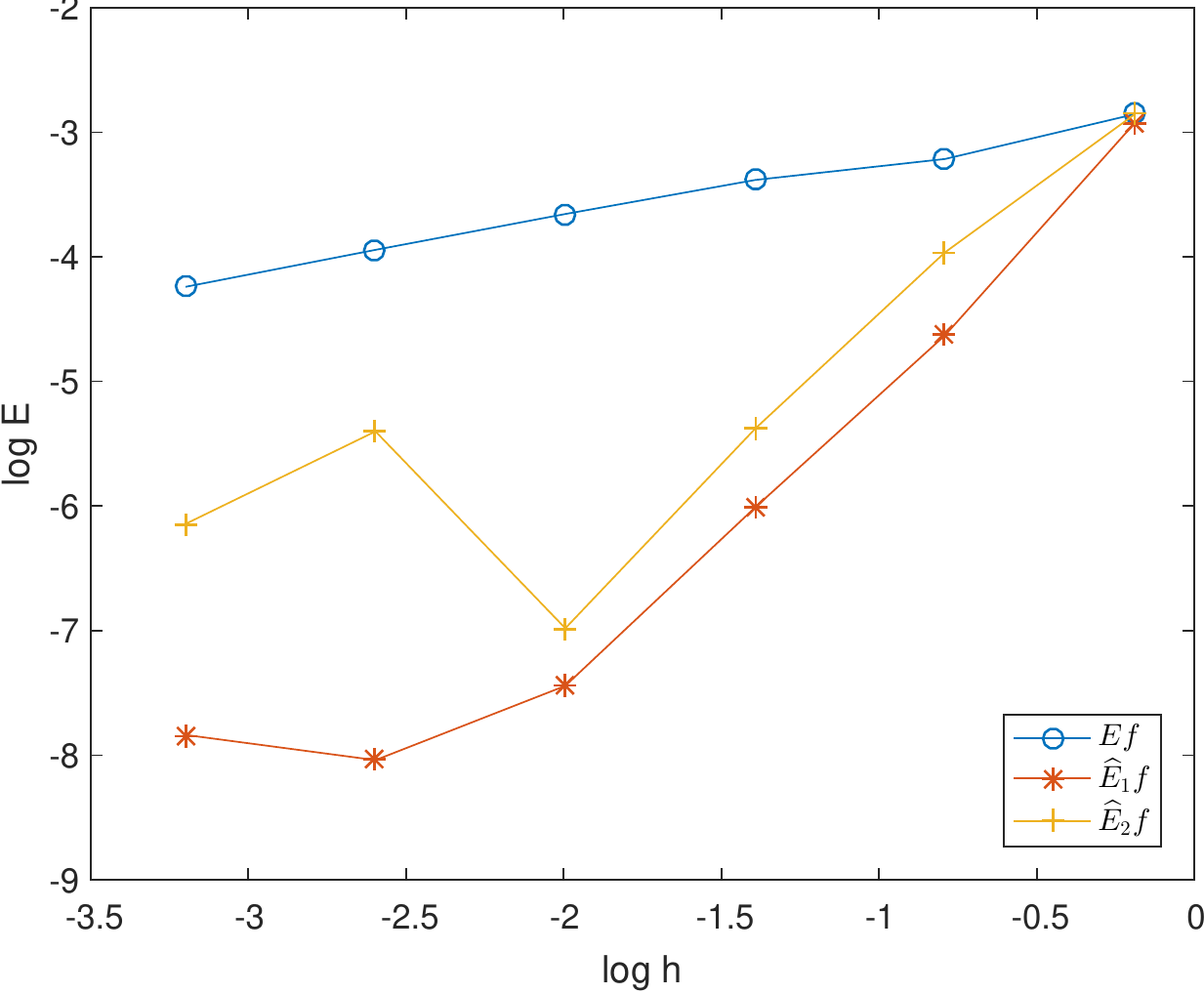}
\caption{$f_2$}
\label{fig:h_test_f2}
\end{subfigure}

\begin{subfigure}{0.5\textwidth}
\includegraphics[width=\textwidth]{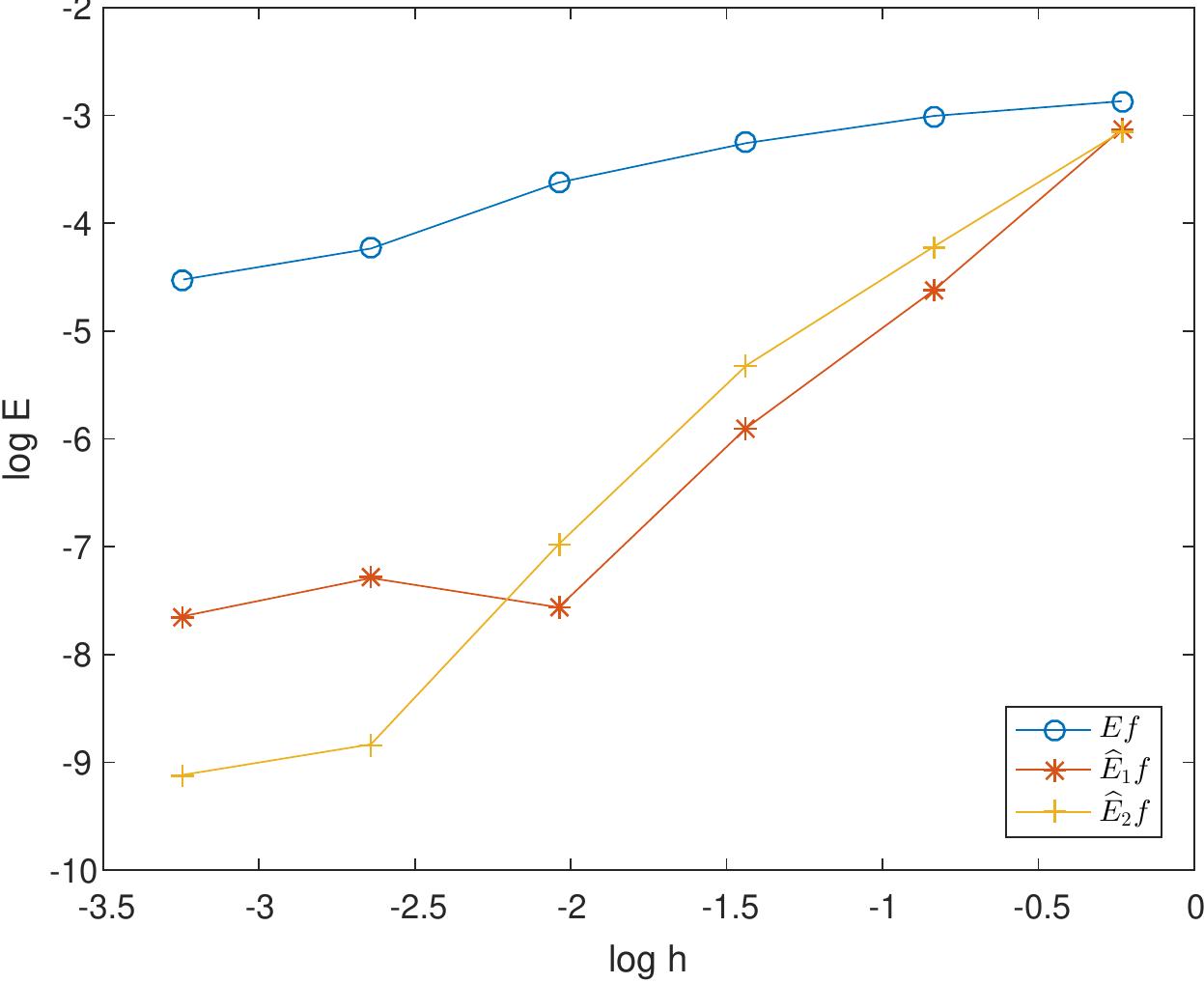}
\caption{$f_3$}
\label{fig:h_test_f3}
\end{subfigure}
\begin{subfigure}{0.5\textwidth}
\includegraphics[width=\textwidth]{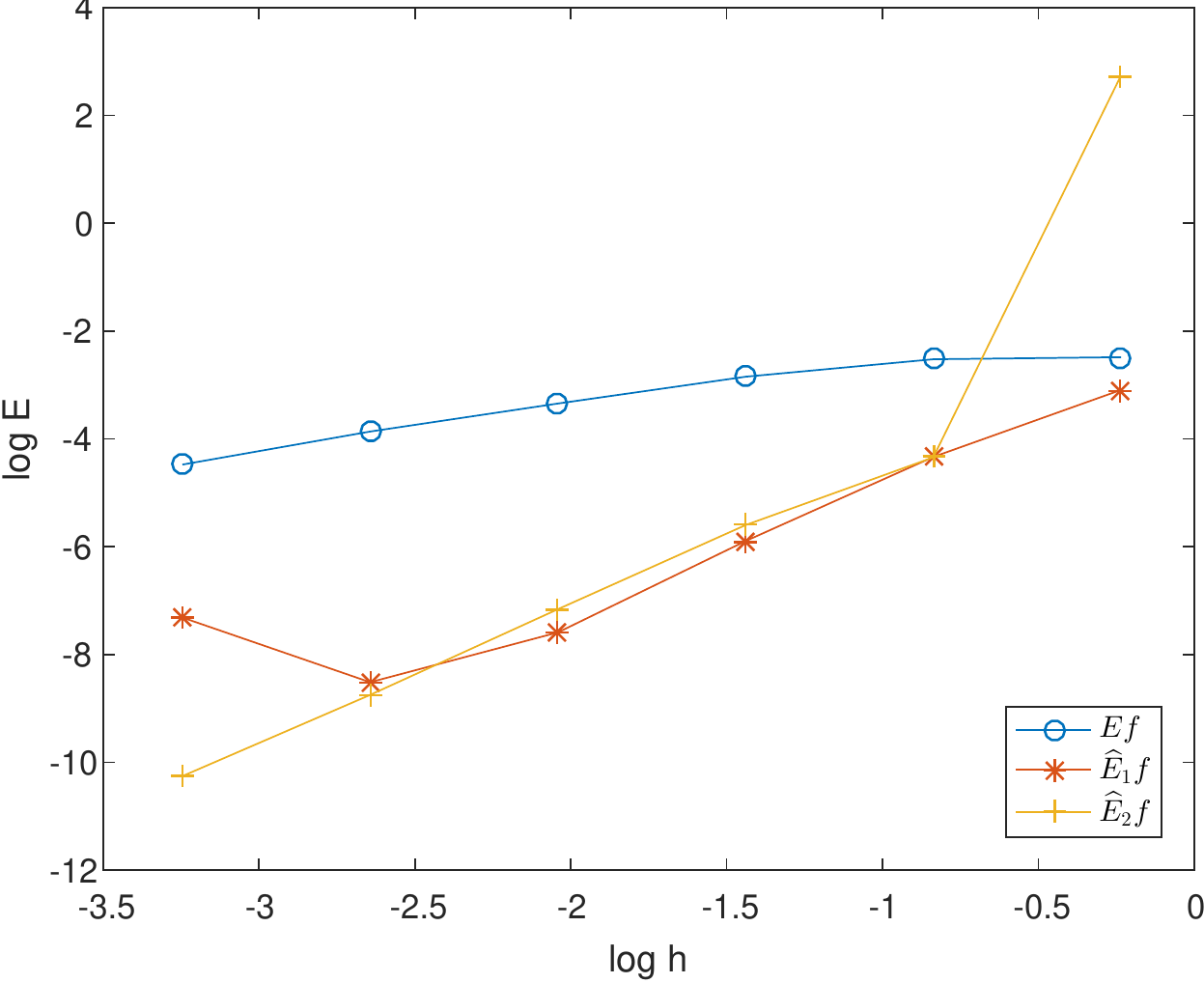}
\caption{$f_4$}
\label{fig:h_test_f4}
\end{subfigure}

\caption{The errors of the original MLS approximation and the corrected approximations for varying $h$ values for $f_1$, $f_2$, $f_3$ and $f_4$.}
\label{fig:h_test}
\end{figure*}

In \cref{fig:h_test}, we compare the original approximation errors to the corrected approximations errors for varying $h$ values.

Here also it is clear that the corrected approximation errors are much more affected by the decreasing $h$ than the original approximation errors.

\begin{figure*}
\begin{subfigure}{0.5\textwidth}
\includegraphics[width=\textwidth]{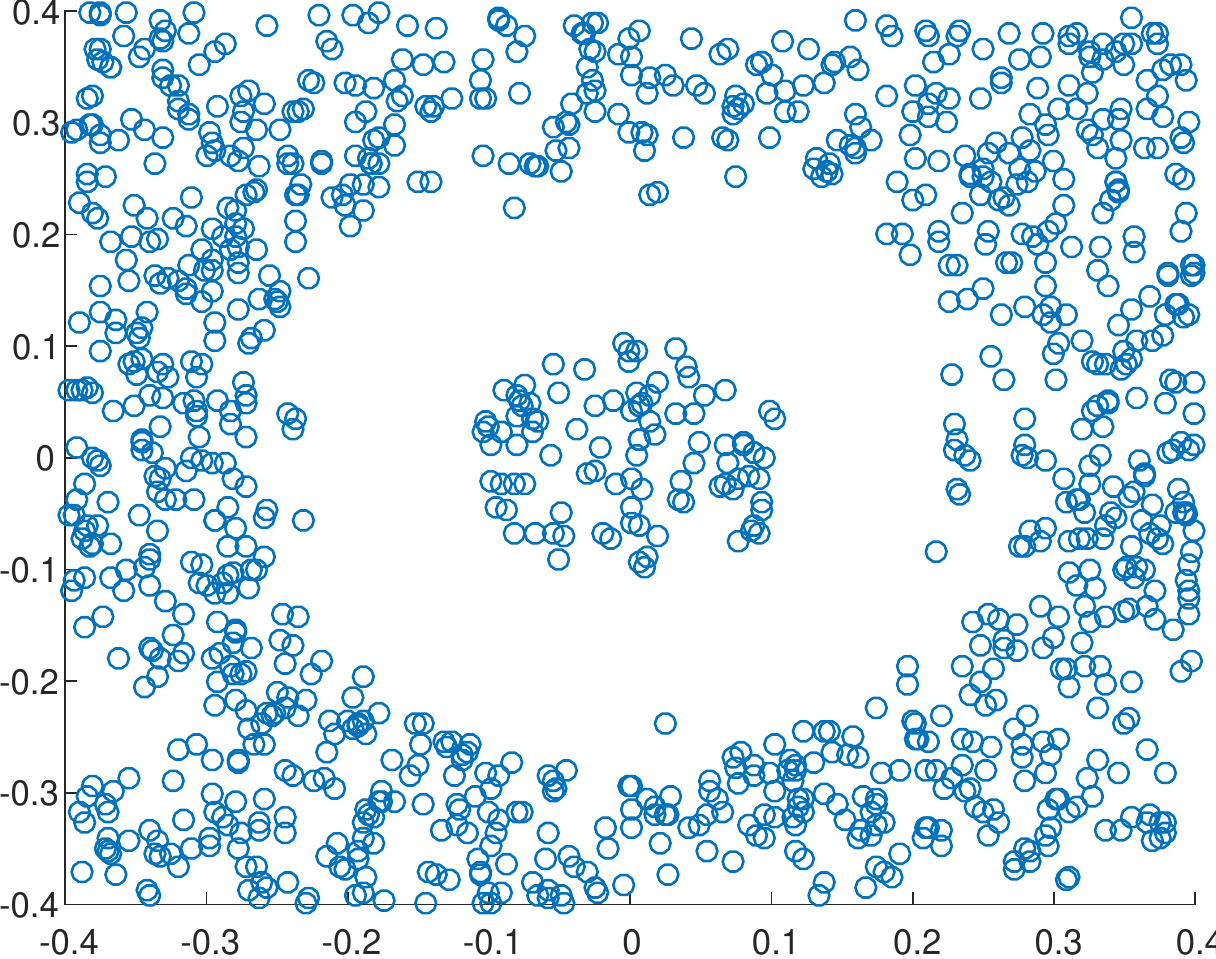}
\caption{$X \setminus \mathcal{S}$}
\label{fig:r1_X-S}
\end{subfigure}
\begin{subfigure}{0.5\textwidth}
\includegraphics[width=\textwidth]{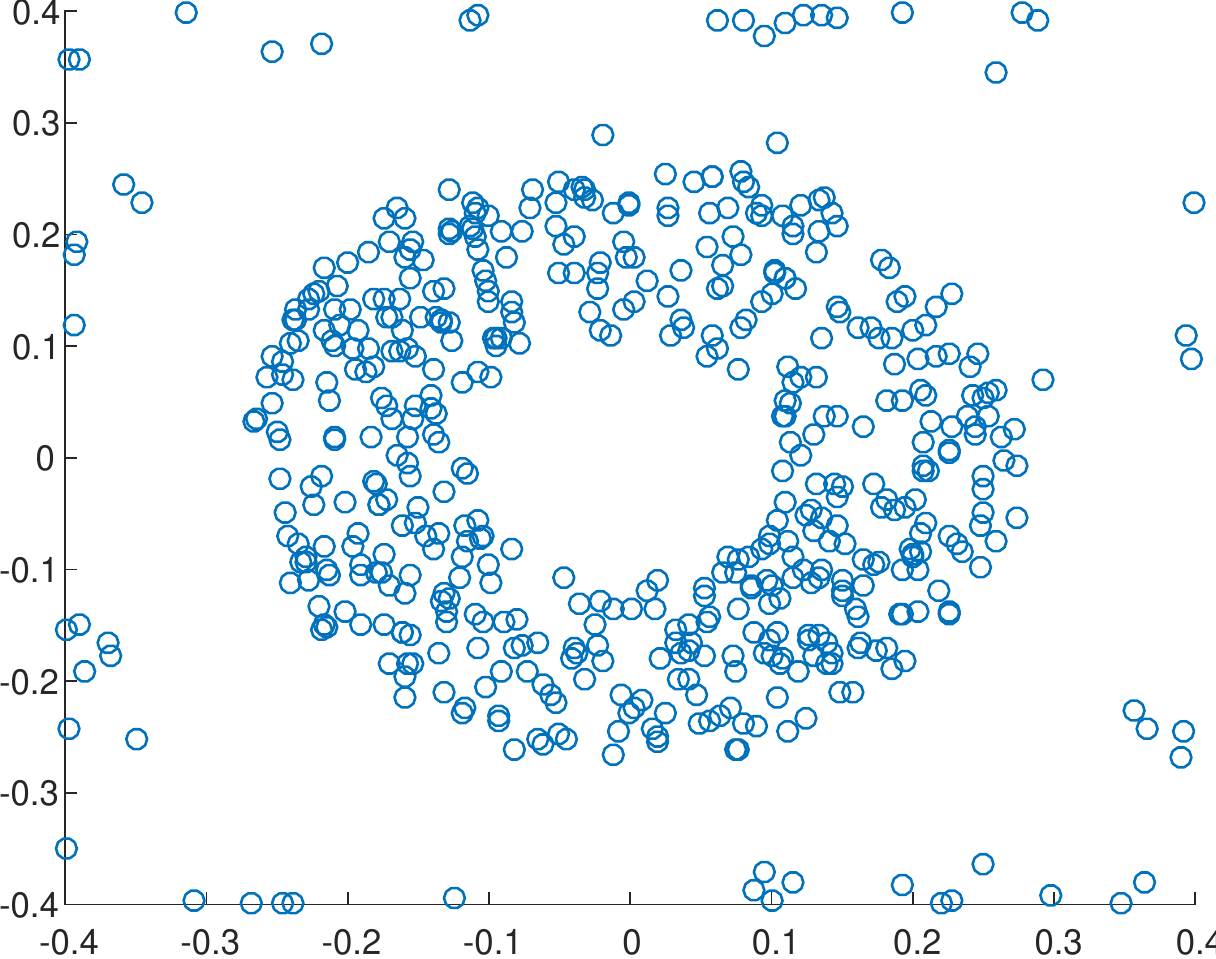}
\caption{$\mathcal{S}$}
\label{fig:r1_S}
\end{subfigure}

\begin{subfigure}{0.5\textwidth}
\includegraphics[width=\textwidth]{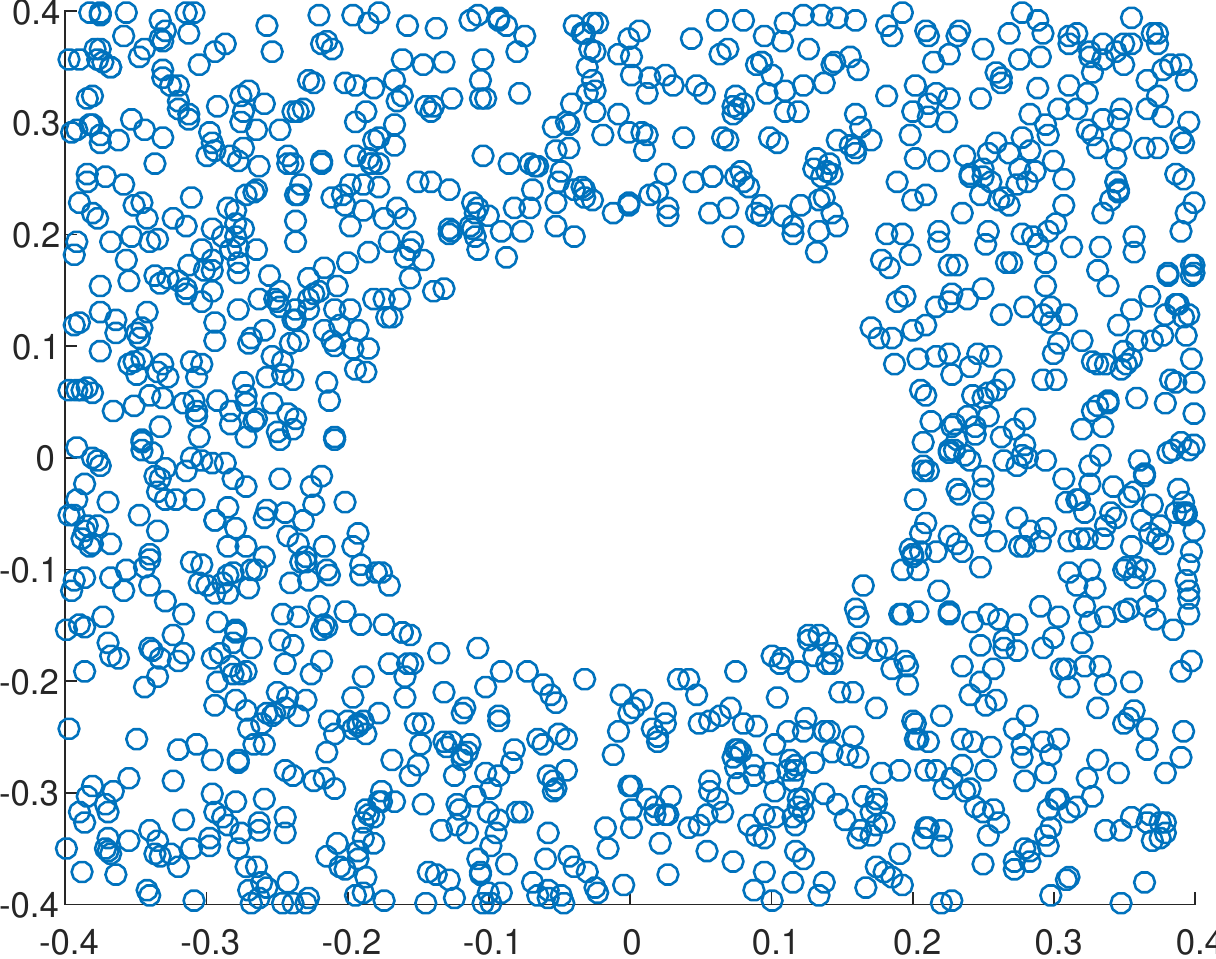}
\caption{$\mathcal{P}$}
\label{fig:r1_P}
\end{subfigure}
\begin{subfigure}{0.5\textwidth}
\includegraphics[width=\textwidth]{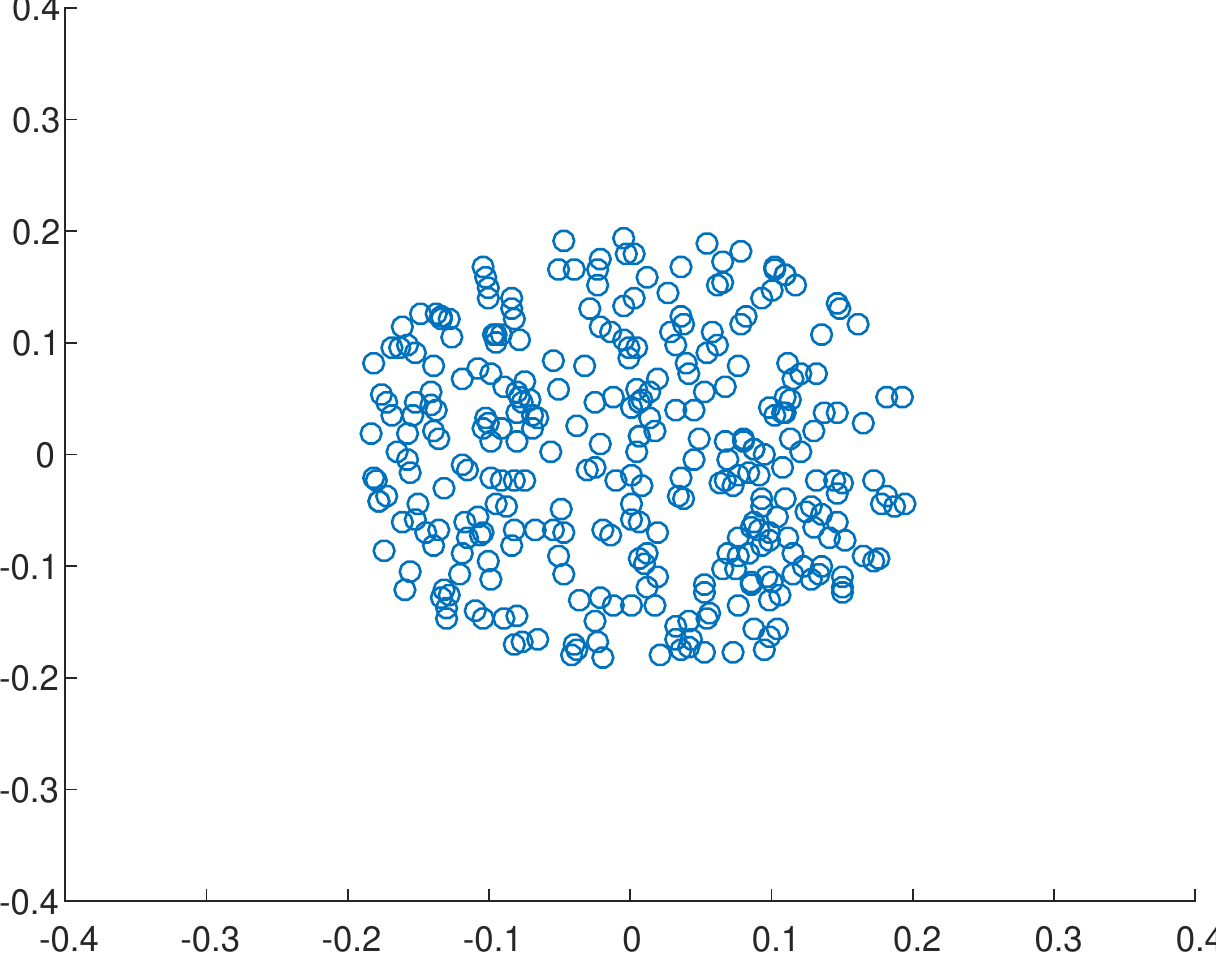}
\caption{$X \setminus \mathcal{P}$}
\label{fig:r1_X-P}
\end{subfigure}

\caption{The results of the partitioning with respect to the signs of $r_1$.}
\label{fig:signs_f1}
\end{figure*}

\begin{figure*}
\begin{subfigure}{0.5\textwidth}
\includegraphics[width=\textwidth]{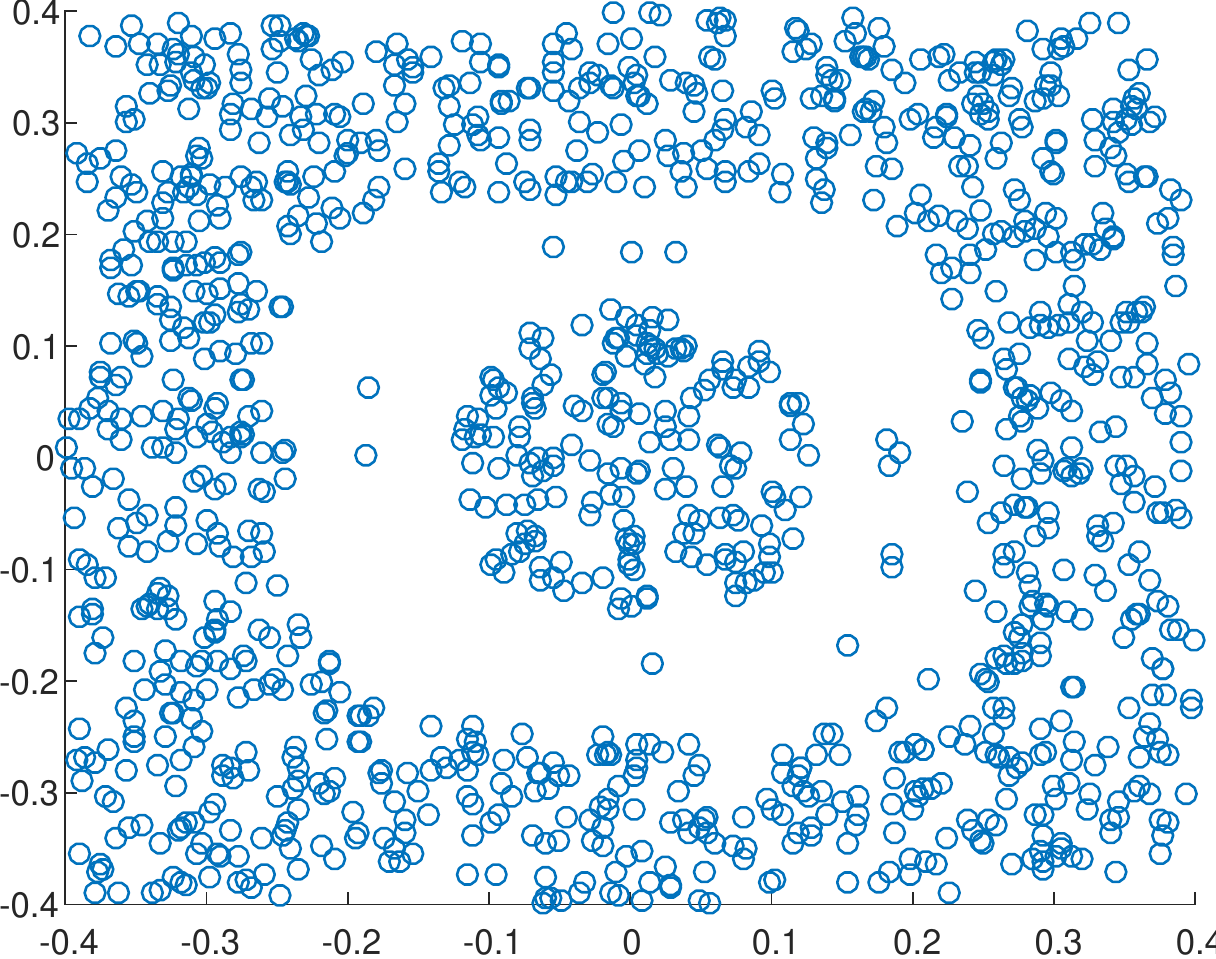}
\caption{$X \setminus \mathcal{S}$}
\label{fig:r2_X-S}
\end{subfigure}
\begin{subfigure}{0.5\textwidth}
\includegraphics[width=\textwidth]{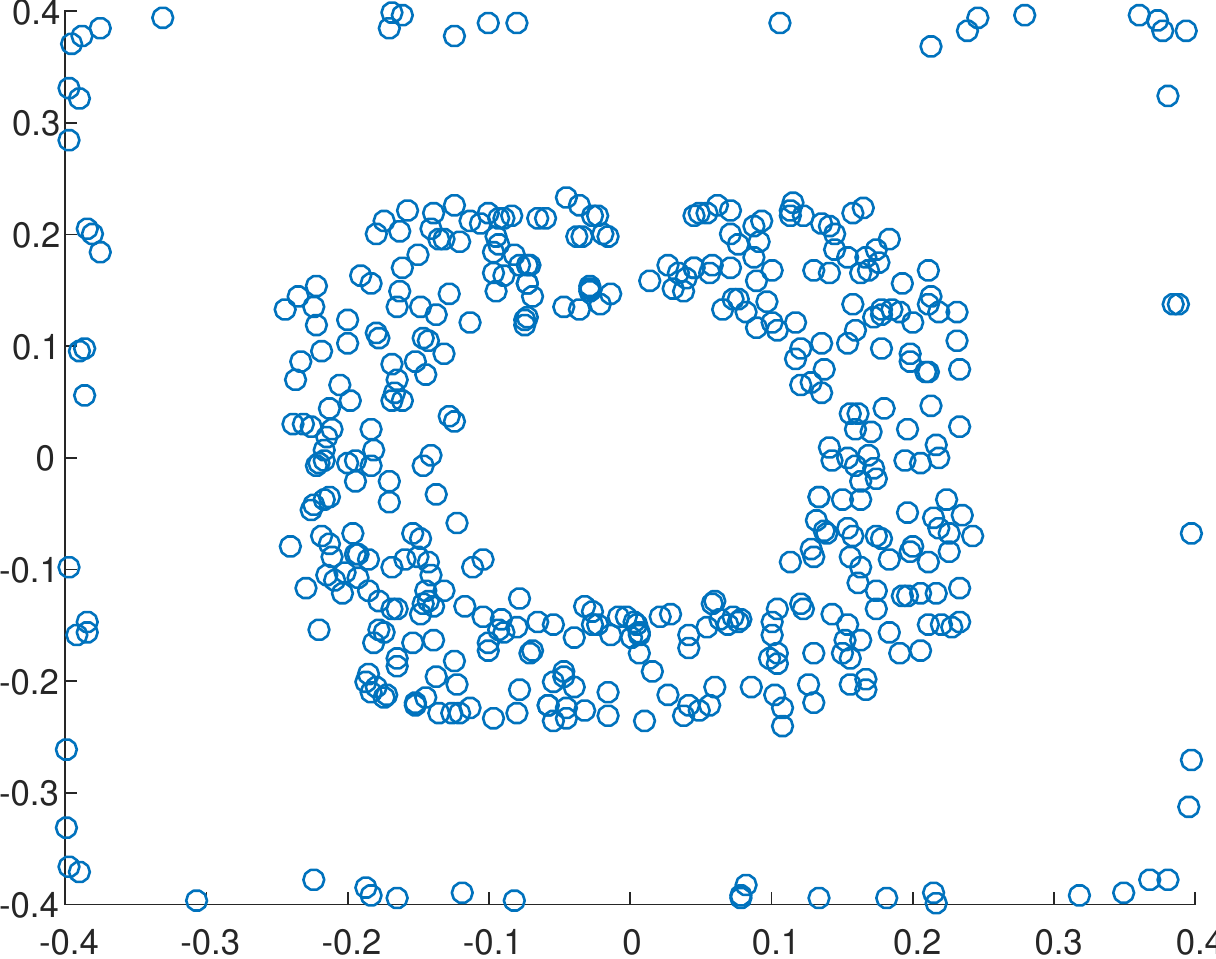}
\caption{$\mathcal{S}$}
\label{fig:r2_S}
\end{subfigure}

\begin{subfigure}{0.5\textwidth}
\includegraphics[width=\textwidth]{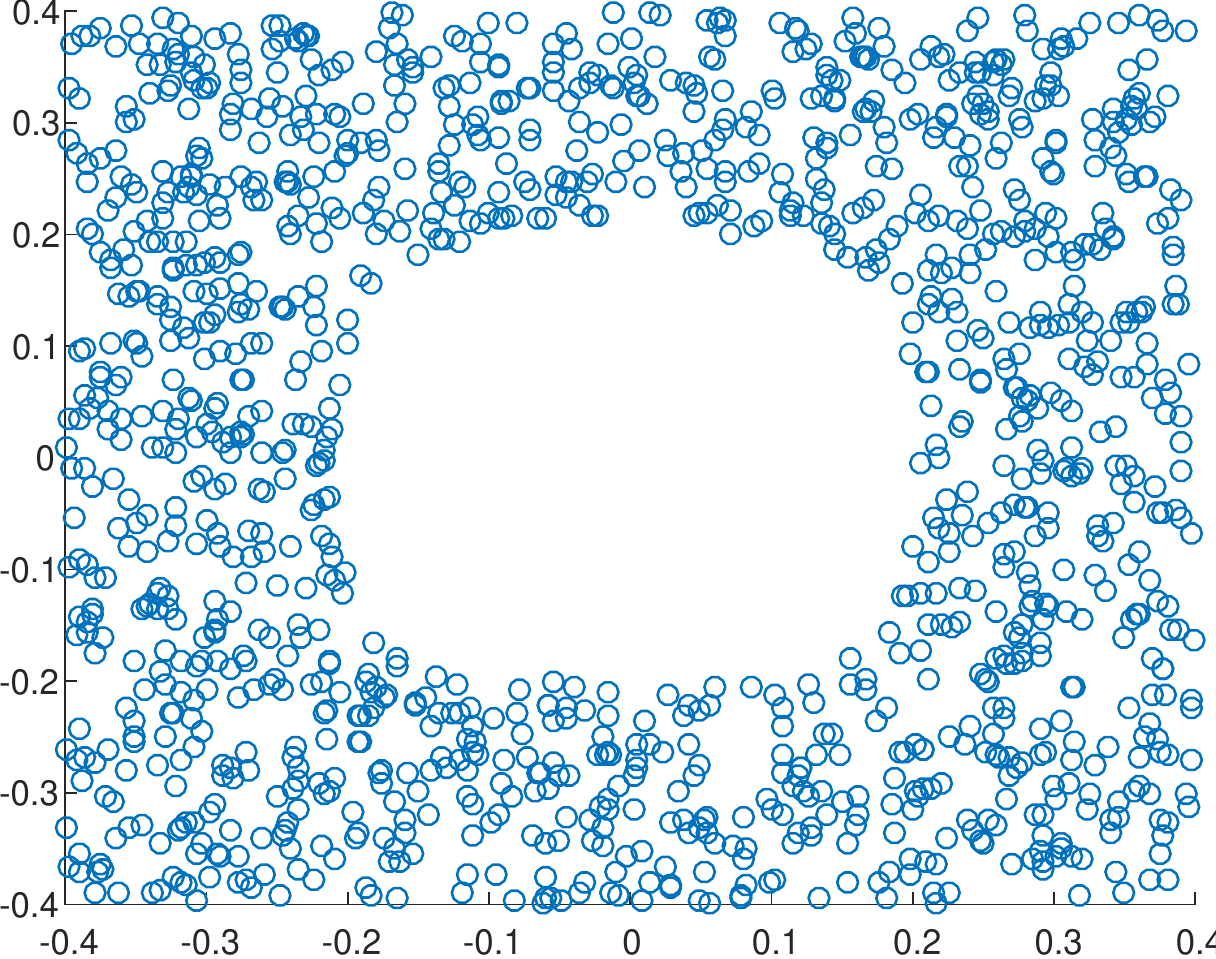}
\caption{$\mathcal{P}$}
\label{fig:r2_P}
\end{subfigure}
\begin{subfigure}{0.5\textwidth}
\includegraphics[width=\textwidth]{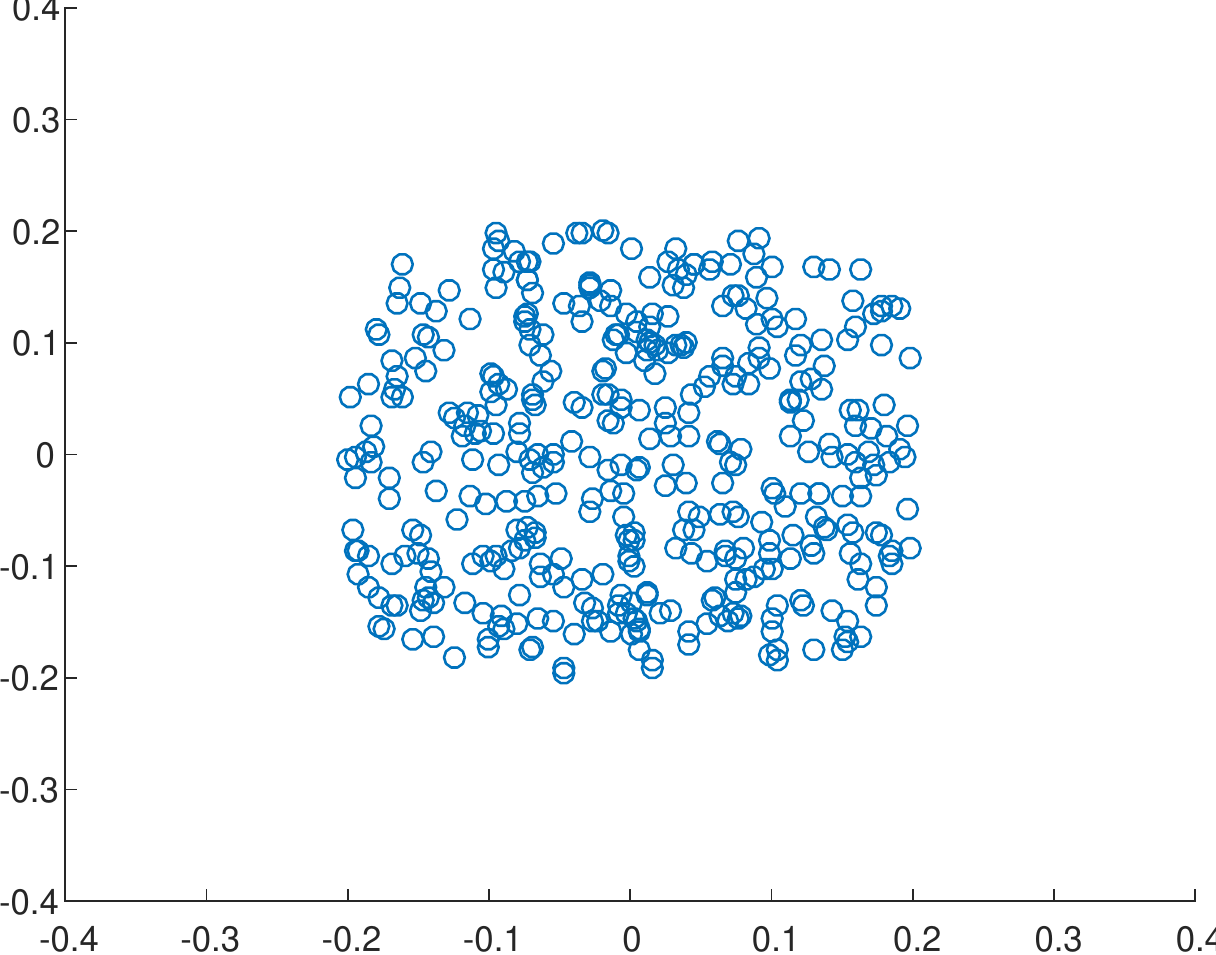}
\caption{$X \setminus \mathcal{P}$}
\label{fig:r2_X-P}
\end{subfigure}

\caption{The results of the partitioning with respect to the signs of $r_2$.}
\label{fig:signs_f2}
\end{figure*}

\begin{figure*}
\begin{subfigure}{0.5\textwidth}
\includegraphics[width=\textwidth]{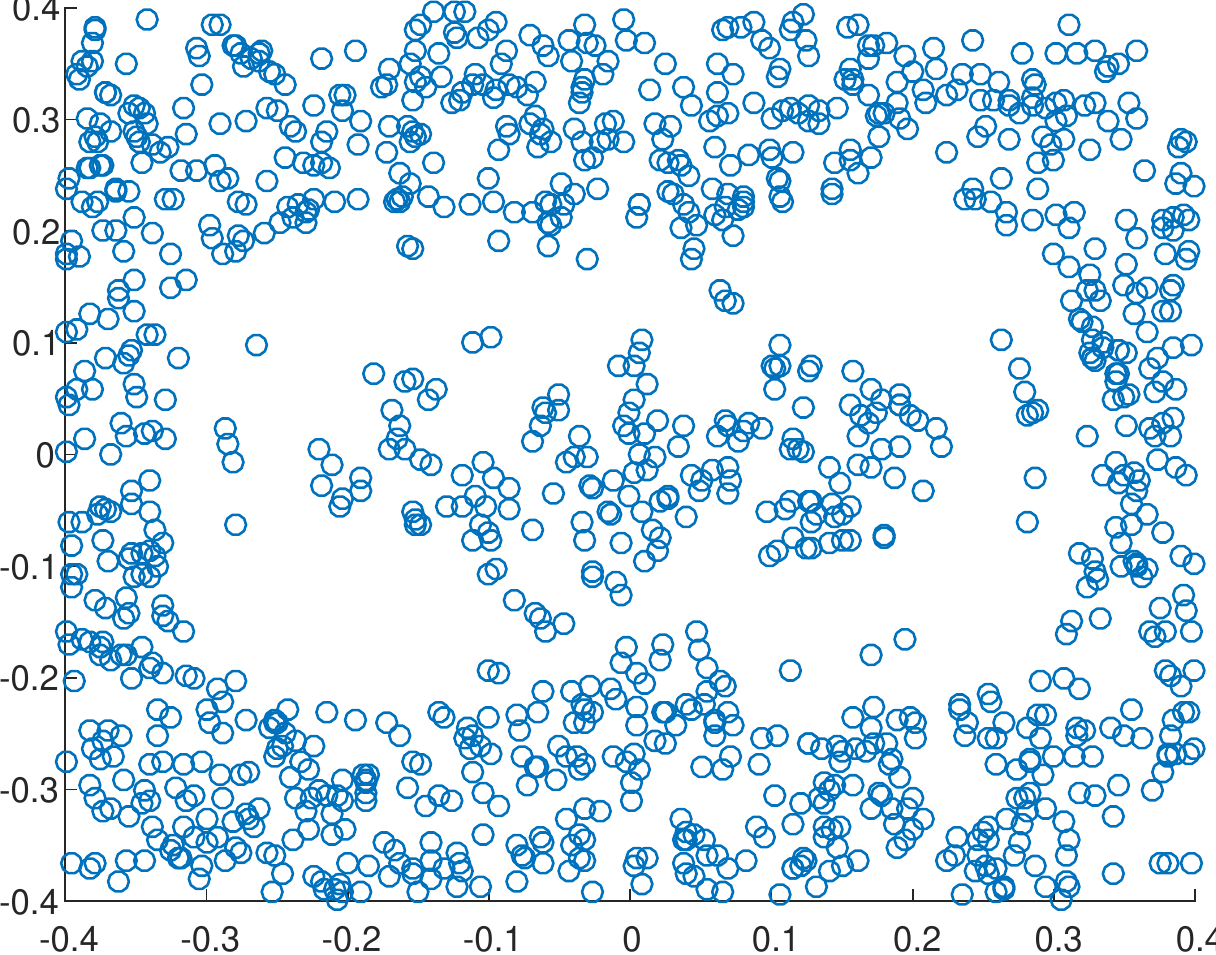}
\caption{$X \setminus \mathcal{S}$}
\label{fig:r3_X-S}
\end{subfigure}
\begin{subfigure}{0.5\textwidth}
\includegraphics[width=\textwidth]{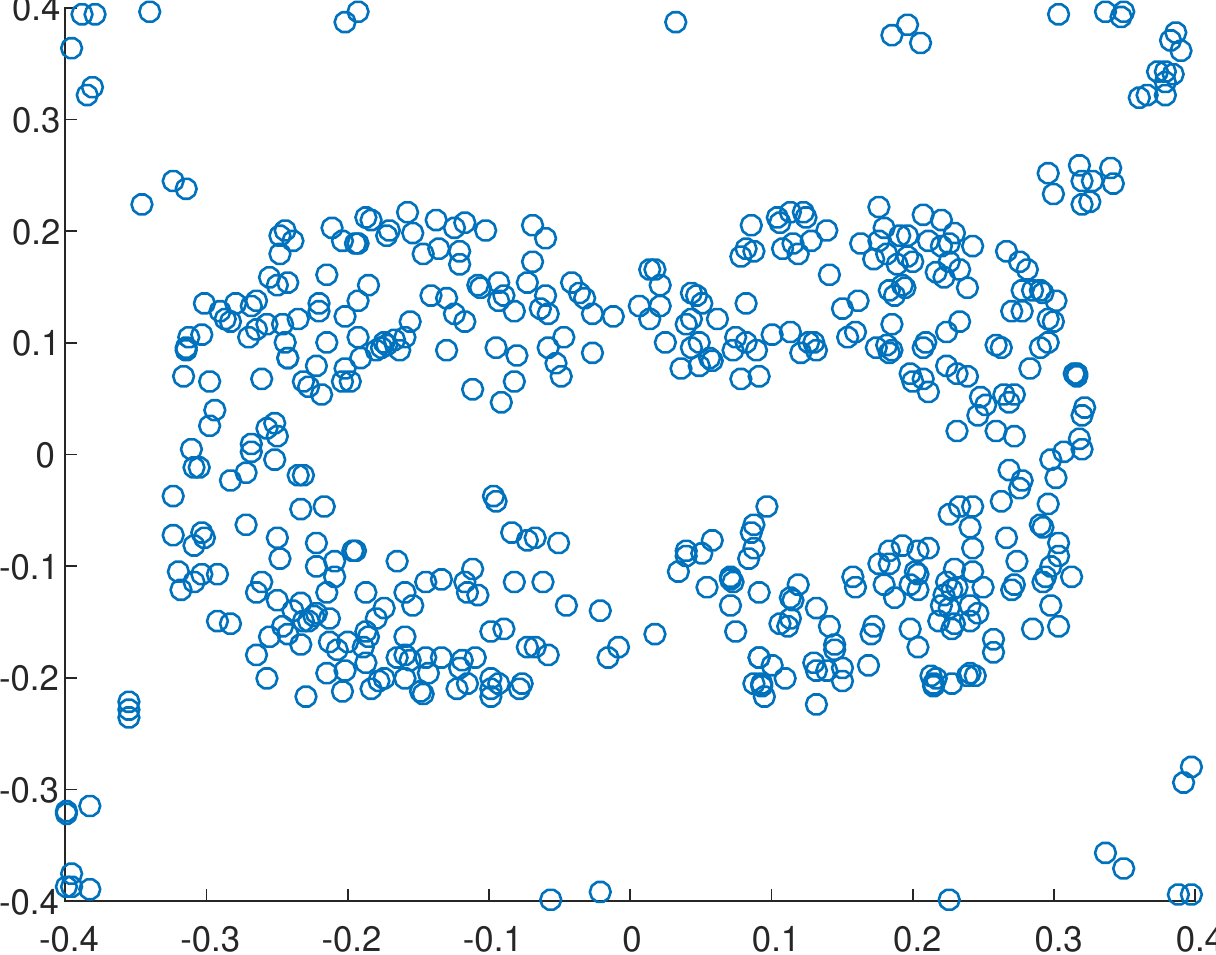}
\caption{$\mathcal{S}$}
\label{fig:r3_S}
\end{subfigure}

\begin{subfigure}{0.5\textwidth}
\includegraphics[width=\textwidth]{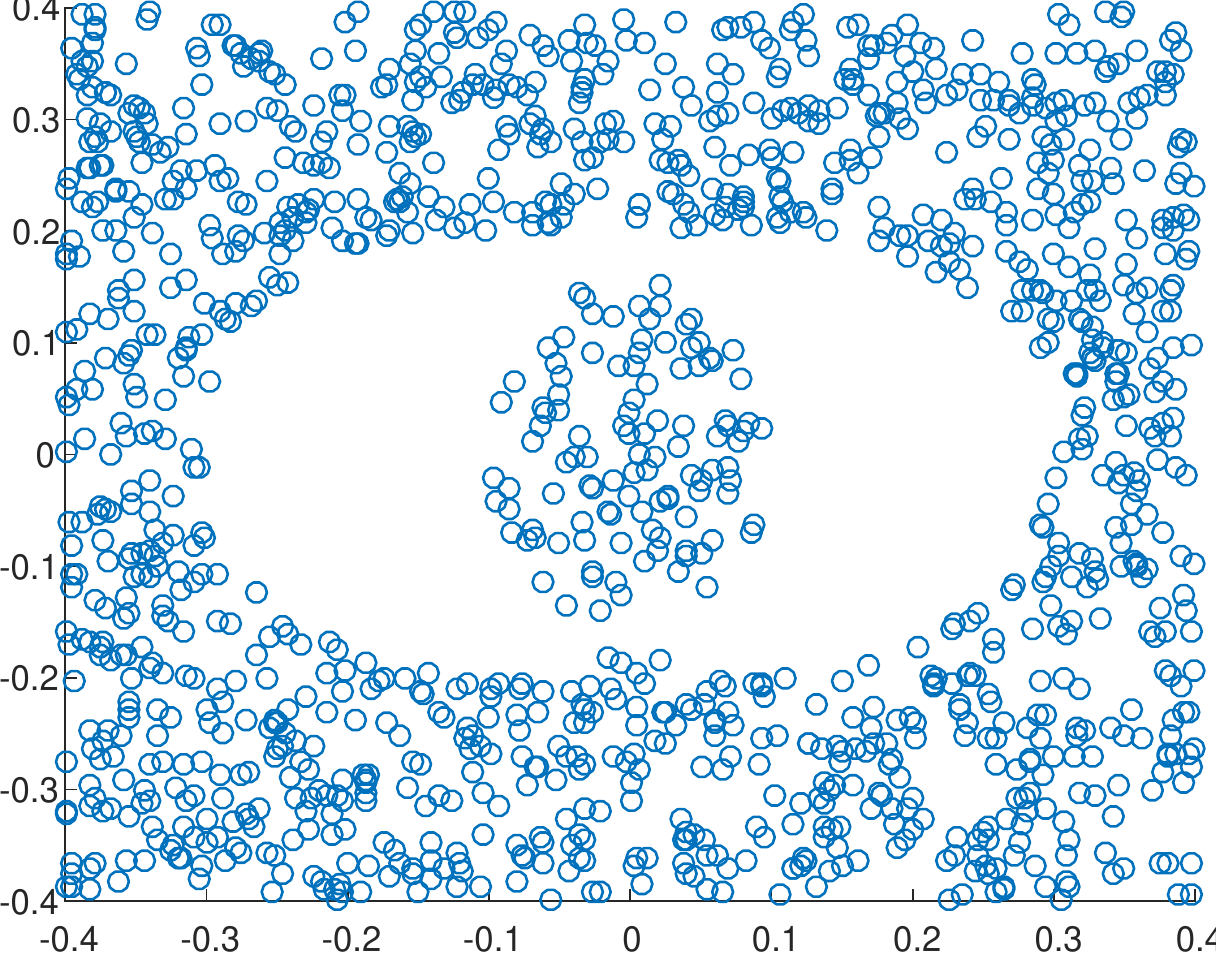}
\caption{$\mathcal{P}$}
\label{fig:r3_P}
\end{subfigure}
\begin{subfigure}{0.5\textwidth}
\includegraphics[width=\textwidth]{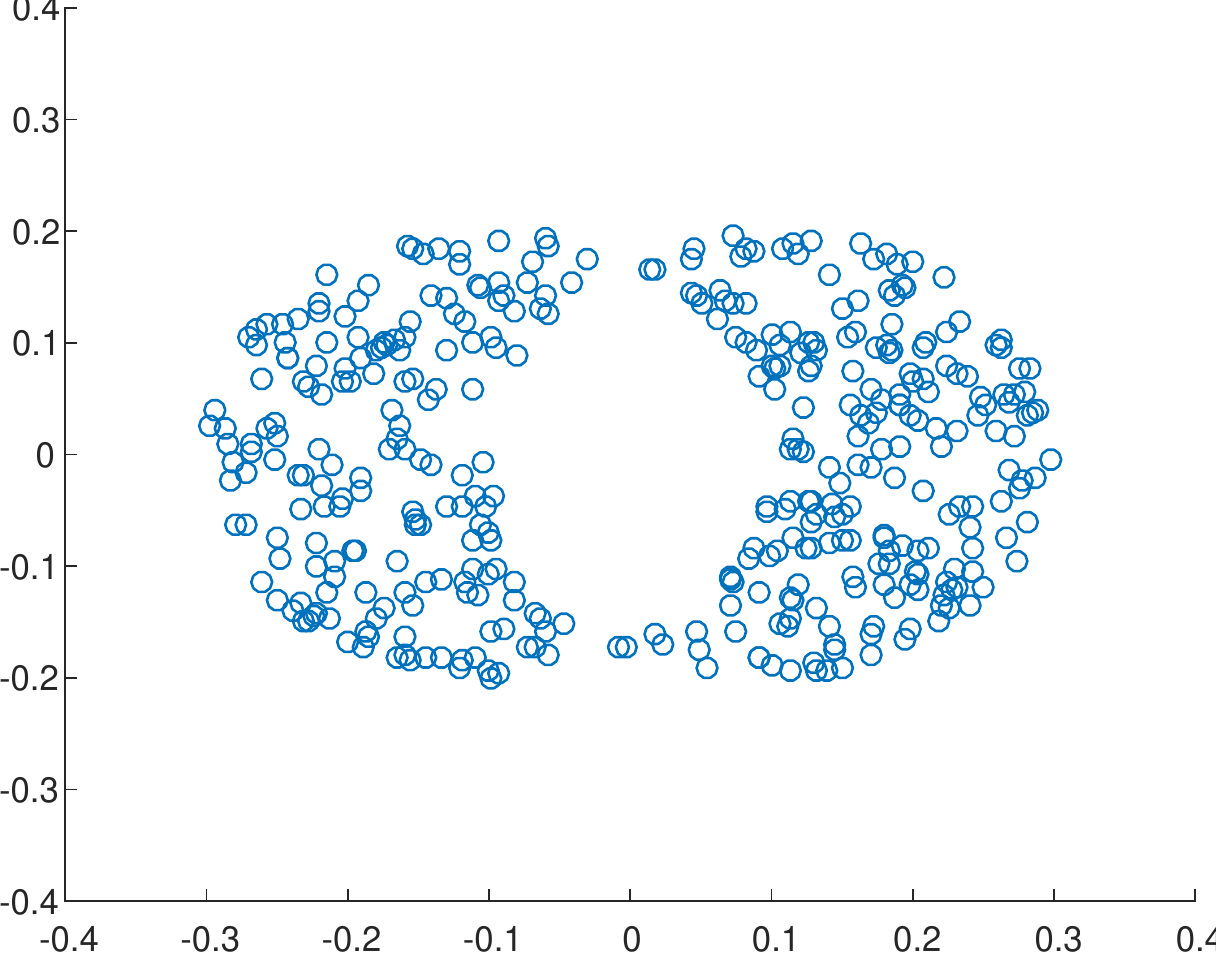}
\caption{$X \setminus \mathcal{P}$}
\label{fig:r3_X-P}
\end{subfigure}

\caption{The results of the partitioning with respect to the signs of $r_3$.}
\label{fig:signs_f3}
\end{figure*}

\begin{figure*}
\begin{subfigure}{0.5\textwidth}
\includegraphics[width=\textwidth]{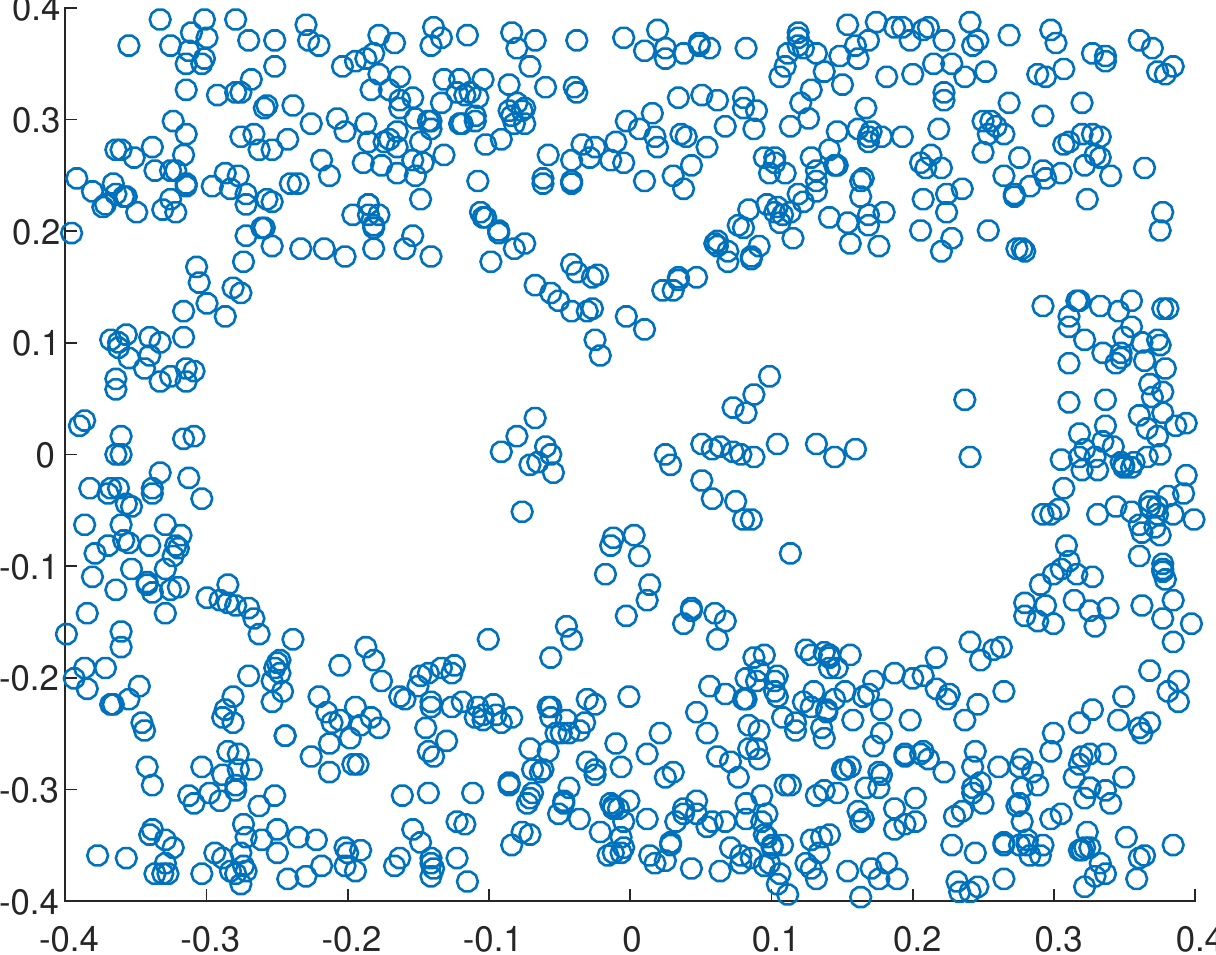}
\caption{$X \setminus \mathcal{S}$}
\label{fig:r4_X-S}
\end{subfigure}
\begin{subfigure}{0.5\textwidth}
\includegraphics[width=\textwidth]{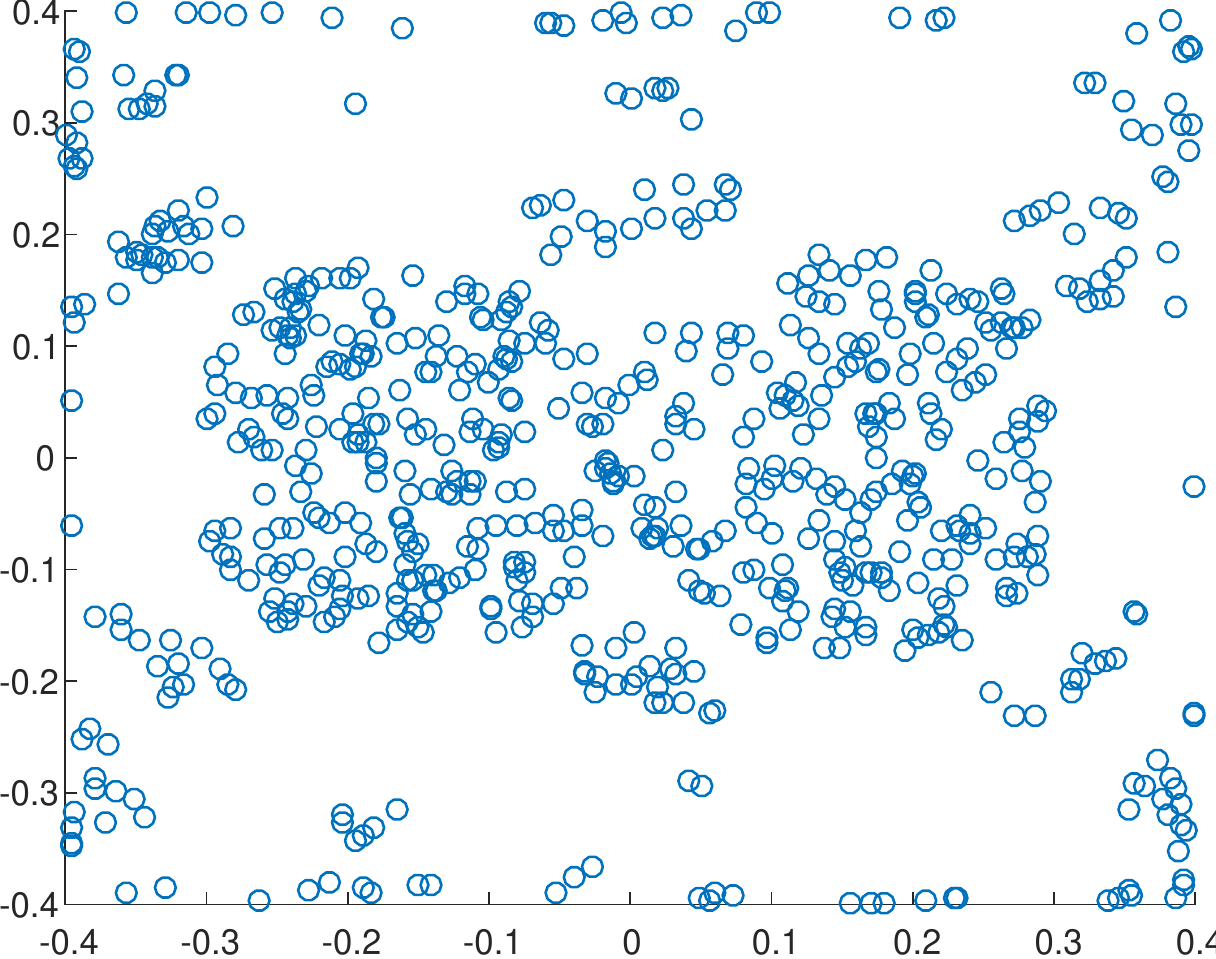}
\caption{$\mathcal{S}$}
\label{fig:r4_S}
\end{subfigure}

\begin{subfigure}{0.5\textwidth}
\includegraphics[width=\textwidth]{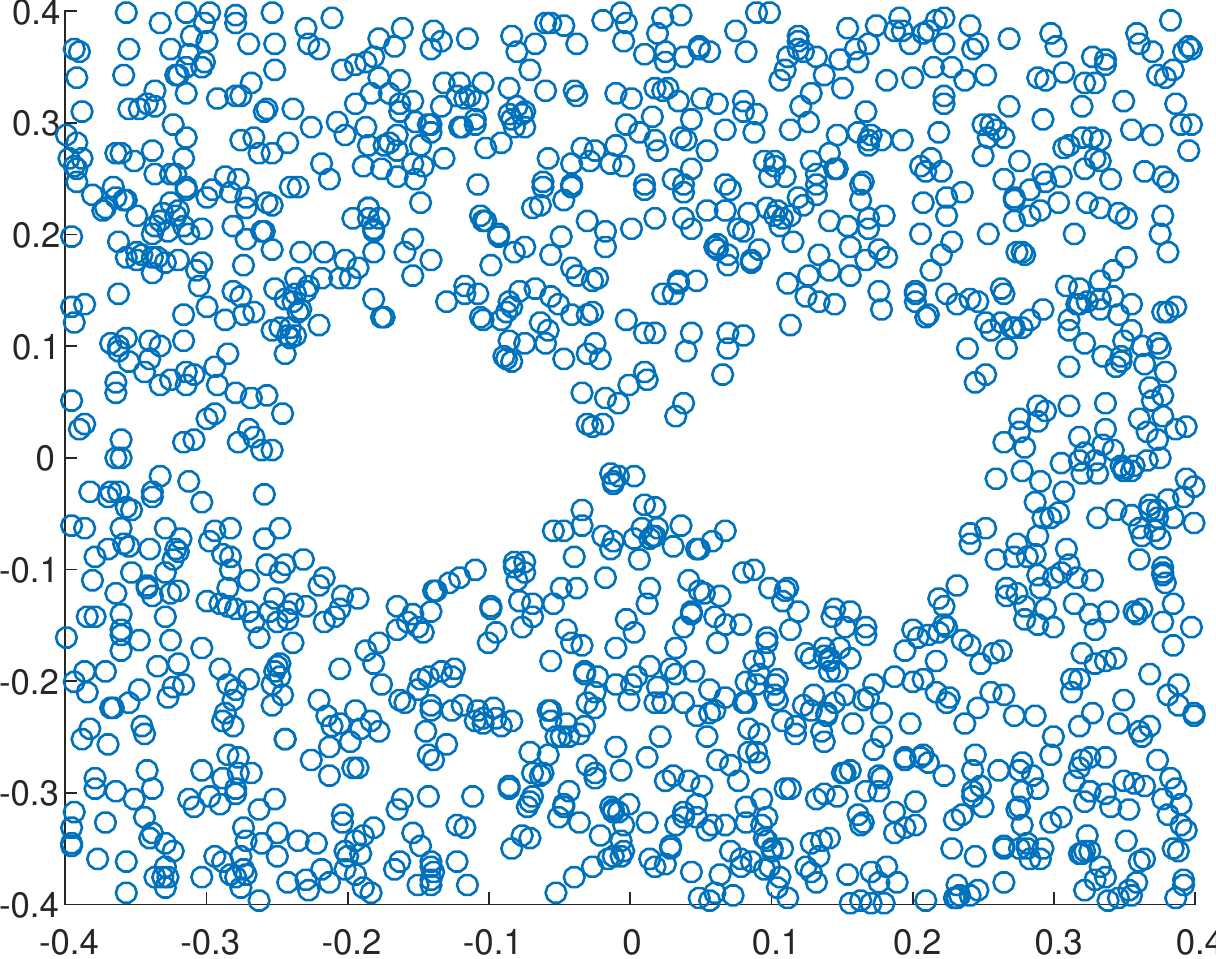}
\caption{$\mathcal{P}$}
\label{fig:r4_P}
\end{subfigure}
\begin{subfigure}{0.5\textwidth}
\includegraphics[width=\textwidth]{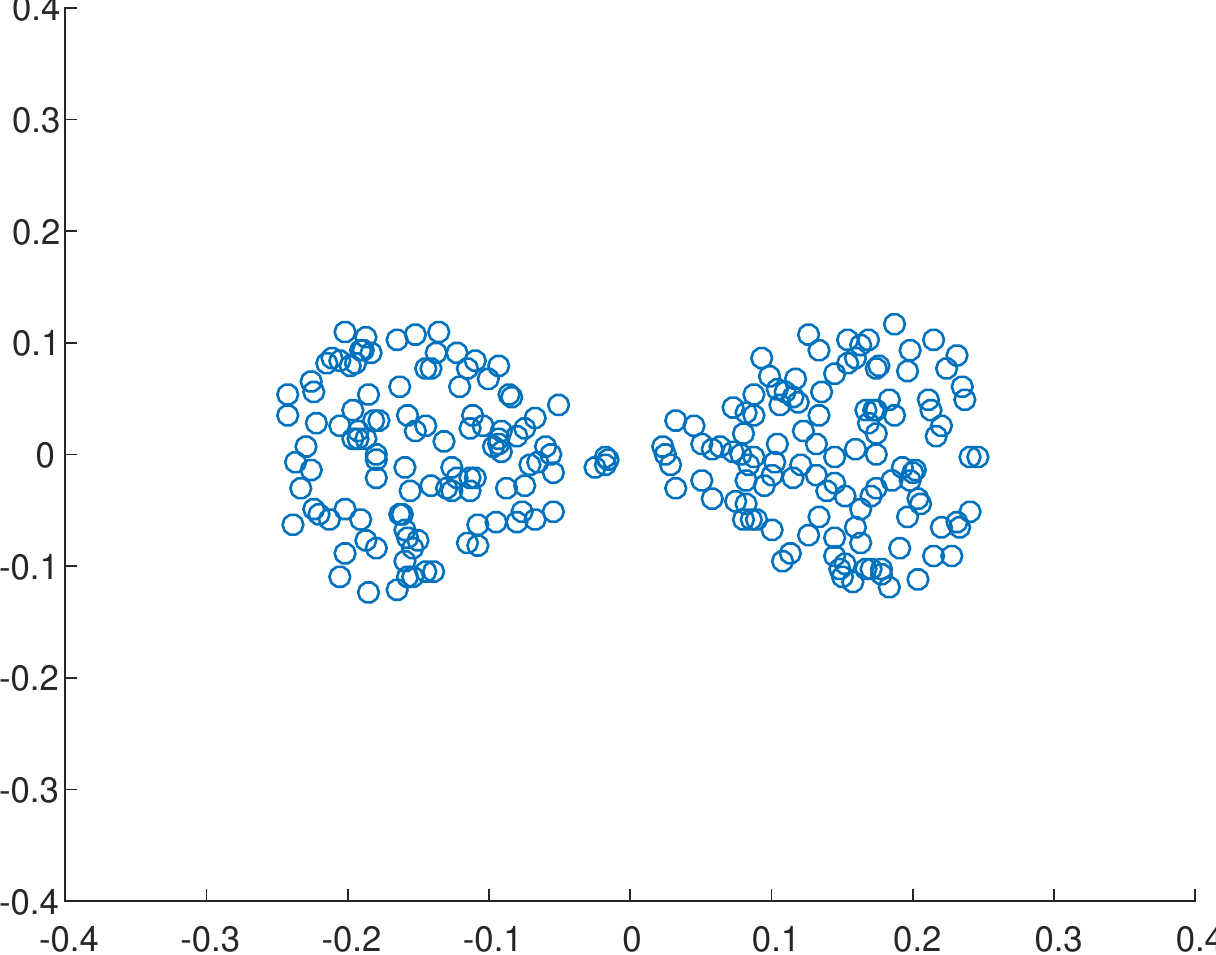}
\caption{$X \setminus \mathcal{P}$}
\label{fig:r4_X-P}
\end{subfigure}

\caption{The results of the partitioning with respect to the signs of $r_4$.}
\label{fig:signs_f4}
\end{figure*}

We also tested our partitioning with respect to the sign of $r_k$ algorithm (see \cref{sec:signs}) on our data set for $1 \leq k \leq 4$. In \cref{fig:signs_f1,fig:signs_f2,fig:signs_f3,fig:signs_f4} we show the steps of the partitioning algorithm.

Note how well the algorithm succeeds in partitioning the data points without any knowledge of the function $r$. 
\bibliography{singular}{}

\begin{thebibliography}{10}

\bibitem{Arandiga}
F.~Arandiga, A.~Cohen, R.~Donat, and N.~Dyn.
\newblock Interpolation and approximation of piecewise smooth functions.
\newblock {\em SIAM Journal on Numerical Analysis}, 43(1):41--57, 2005.

\bibitem{Gelb2009}
R.~Archibald, A.~Gelb, R.~Saxena, and D.~Xiu.
\newblock Discontinuity detection in multivariate space for stochastic
  simulations.
\newblock {\em Journal of Computational Physics}, 228(7):2676--2689, 2009.

\bibitem{Gelb2008}
R.~Archibald, A.~Gelb, and J.~Yoon.
\newblock Determining the locations and discontinuities in the derivatives of
  functions.
\newblock {\em Applied Numerical Mathematics}, 58(5):577--592, 2008.

\bibitem{batenkov2012complete}
D.~Batenkov.
\newblock Complete algebraic reconstruction of piecewise-smooth functions from
  fourier data.
\newblock {\em arXiv preprint arXiv:1211.0680}, 2012.

\bibitem{Batenkov2011algebraic}
D.~Batenkov, N.~Sarig, and Y.~Yomdin.
\newblock Algebraic reconstruction of piecewise-smooth functions from integral
  measurements.
\newblock {\em arXiv preprint arXiv:1103.3969}, 2011.

\bibitem{Batenkov2012algebraic}
D.~Batenkov and Y.~Yomdin.
\newblock Algebraic fourier reconstruction of piecewise smooth functions.
\newblock {\em Mathematics of Computation}, 81(277):277--318, 2012.

\bibitem{bosMLS}
L.~Bos and K.~Salkauskas.
\newblock Moving least-squares are backus-gilbert optimal.
\newblock {\em Journal of Approximation Theory}, 59(3):267--275, 1989.

\bibitem{Harten}
A.~Harten.
\newblock {ENO} schemes with subcell resolution.
\newblock {\em Journal of Computational Physics}, 83(1):148 -- 184, 1989.

\bibitem{levin1998approximation}
D.~Levin.
\newblock The approximation power of moving least-squares.
\newblock {\em Mathematics of Computation of the American Mathematical
  Society}, 67(224):1517--1531, 1998.

\bibitem{Lipman}
Y.~Lipman and D.~Levin.
\newblock Approximating piecewise-smooth functions.
\newblock {\em IMA Journal of Numerical Analysis}, 30(4):1159--1183, 2010.

\bibitem{markakis2014high}
C.~Markakis and L.~Barack.
\newblock High-order difference and pseudospectral methods for discontinuous
  problems.
\newblock {\em arXiv preprint arXiv:1406.4865}, 2014.

\bibitem{Acoustics}
P.~Morse and K.~Ingard.
\newblock {\em Theoretical Acoustics}.
\newblock International series in pure and applied physics. Princeton
  University Press, 1986.

\bibitem{LevelSet}
S.~Osher and R.~Fedkiw.
\newblock {\em Level Set Methods and Dynamic Implicit Surfaces}.
\newblock Applied Mathematical Sciences. Springer, 2003.

\bibitem{Plaskota2009power}
L.~Plaskota and G.~W. Wasilkowski.
\newblock The power of adaptive algorithms for functions with singularities.
\newblock {\em Journal of fixed point theory and applications}, 6(2):227--248,
  2009.

\bibitem{Plaskota2013adaptive}
L.~Plaskota, G.~W. Wasilkowski, and Y.~Zhao.
\newblock An adaptive algorithm for weighted approximation of singular
  functions over r.
\newblock {\em SIAM Journal on Numerical Analysis}, 51(3):1470--1493, 2013.

\bibitem{wendland2004scattered}
H.~Wendland.
\newblock {\em Scattered data approximation}, volume~17.
\newblock Cambridge university press, 2004.

\end{thebibliography}

\end{document}